\documentclass[11pt]{article}
\usepackage{smile}

\usepackage{fullpage}
\usepackage{lscape}
\usepackage{bigints}
\usepackage{framed}
\usepackage{mdframed}
\usepackage{enumerate}
\usepackage[inline]{enumitem}
\usepackage[T1]{fontenc}
\usepackage{moresize}
\usepackage{bm}
\usepackage{bbm}
\usepackage{dsfont}
\usepackage{amsmath}
\usepackage{amssymb}
\usepackage{amsthm}
\usepackage{amsfonts}
\usepackage{stmaryrd}
\usepackage{array}
\usepackage{mathrsfs}
\usepackage{mathtools} 
\usepackage{extarrows}
\usepackage{stackrel}
\usepackage{relsize,exscale}
\usepackage{scalerel}
\usepackage[nodisplayskipstretch]{setspace}
\usepackage{color}
\usepackage[usenames,dvipsnames]{xcolor}
\usepackage{cancel}
\usepackage{soul}
\usepackage{undertilde}
\usepackage{xfrac}
\usepackage{siunitx}
\usepackage{graphicx}
\usepackage{float}
\usepackage{rotating}
\usepackage{subcaption}
\usepackage{overpic}
\usepackage[all]{xy}
\usepackage{tikz}
\usetikzlibrary{arrows,matrix,positioning,calc,automata,patterns}
\usepackage{booktabs}
\usepackage{dcolumn}
\usepackage{multirow}
\usepackage{diagbox}
\usepackage{tabularx}
\usepackage{verbatim}
\usepackage{listings}
\usepackage[ruled,vlined]{algorithm2e}
\usepackage{fancyvrb}
\usepackage{hyperref}
\usepackage[round]{natbib}
\usepackage{sectsty}

\hypersetup{
    bookmarks=true,         
    unicode=false,          
    pdftoolbar=true,        
    pdfmenubar=true,        
    pdffitwindow=false,     
    pdfstartview={FitH},    
    pdftitle={My title},    
    pdfauthor={Author},     
    pdfsubject={Subject},   
    pdfcreator={Creator},   
    pdfproducer={Producer}, 
    pdfkeywords={key1, key2}, 
    pdfnewwindow=true,      
    colorlinks=true,        
    linkcolor=blue,         
    citecolor=blue,         
    filecolor=blue,         
    urlcolor=cyan           
}

\usepackage{stackengine}
\stackMath
\newcommand\tenq[2][1]{%
\def\useanchorwidth{T}%
\ifnum#1>1%
\stackunder[0pt]{\tenq[\numexpr#1-1\relax]{#2}}{\!\scriptscriptstyle\thicksim}%
\else%
\stackunder[1pt]{#2}{\!\scriptstyle\thicksim}%
\fi%
}

\makeatletter
\DeclareRobustCommand\widecheck[1]{{\mathpalette\@widecheck{#1}}}
\def\@widecheck#1#2{%
    \setbox\z@\hbox{\m@th$#1#2$}%
    \setbox\tw@\hbox{\m@th$#1%
       \widehat{%
          \vrule\@width\z@\@height\ht\z@
          \vrule\@height\z@\@width\wd\z@}$}%
    \dp\tw@-\ht\z@
    \@tempdima\ht\z@ \advance\@tempdima2\ht\tw@ \divide\@tempdima\thr@@
    \setbox\tw@\hbox{%
       \raise\@tempdima\hbox{\scalebox{1}[-1]{\lower\@tempdima\box
\tw@}}}%
    {\ooalign{\box\tw@ \cr \box\z@}}}
\makeatother

\def\given{\,|\,}

\def\Biggiven{\,\Big{|}\,}
\def\tr{\mathop{\text{tr}}\kern.2ex}

\def\P{{\mathrm P}}

\def\E{{\mathrm E}}

\def\d{{\mathrm d}}

\newcommand{\zahl}[1]{\llbracket #1\rrbracket}
\newcommand\yestag{\addtocounter{equation}{1}\tag{\theequation}}
\newcolumntype{L}[1]{>{\raggedright\let\newline\\\arraybackslash\hspace{0pt}}m{#1}}
\newcolumntype{C}[1]{>{  \centering\let\newline\\\arraybackslash\hspace{0pt}}m{#1}}
\newcolumntype{R}[1]{>{ \raggedleft\let\newline\\\arraybackslash\hspace{0pt}}m{#1}}
\newcolumntype{d}[1]{D{.}{.}{#1}}
\newcolumntype{H}{>{\setbox0=\hbox\bgroup}c<{\egroup}@{}}
\newcolumntype{Z}{>{\setbox0=\hbox\bgroup}c<{\egroup}@{\hspace*{-\tabcolsep}}}
\newcolumntype{b}{X}
\newcolumntype{s}{>{\hsize=.5\hsize}X}

\numberwithin{equation}{section}

\newtheorem{theorem}{Theorem}[section]
\newtheorem{lemma}{Lemma}[section]
\newtheorem{proposition}{Proposition}[section]
\newtheorem{assumption}{Assumption}[section]
\newtheorem{corollary}{Corollary}[section]
\newtheorem{innercustomgeneric}{\customgenericname}
\providecommand{\customgenericname}{}
\newcommand{\newcustomtheorem}[2]{%
  \newenvironment{#1}[1]
  {%
   \renewcommand\customgenericname{#2}%
   \renewcommand\theinnercustomgeneric{##1}%
   \innercustomgeneric
  }
  {\endinnercustomgeneric}
}
\newcustomtheorem{customdefinition}{Definition}
\newcustomtheorem{customdefinitions}{Definitions}
\newcustomtheorem{customtheorem}{Theorem}
\newcustomtheorem{customassumption}{Assumption}
\newcustomtheorem{customlemma}{Lemma}
\newcustomtheorem{customexample}{Example}
\theoremstyle{definition}

\usepackage{enumitem}
\makeatletter
\newcommand{\mylabel}[2]{#2\def\@currentlabel{#2}\label{#1}}
\makeatother

\graphicspath{{./fig3/}}

\allowdisplaybreaks

\begin{document}

\setlength{\abovedisplayskip}{5pt}
\setlength{\belowdisplayskip}{5pt}
\setlength{\abovedisplayshortskip}{5pt}
\setlength{\belowdisplayshortskip}{5pt}
\hypersetup{colorlinks,breaklinks,urlcolor=blue,linkcolor=blue}

\title{\LARGE On Rosenbaum's Rank-based Matching Estimator}

\author{Matias D. Cattaneo\thanks{Department of Operations Research and Financial Engineering, Princeton University, Princeton, NJ 08544, USA; e-mail: {\tt cattaneo@princeton.edu}} \and
	    Fang Han\thanks{Department of Statistics, University of Washington, Seattle, WA 98195, USA; e-mail: {\tt fanghan@uw.edu}}\and
	    Zhexiao Lin\thanks{Department of Statistics, University of California, Berkeley, CA 94720, USA; e-mail: {\tt zhexiaolin@berkeley.edu}}
}

\date{\today}

\maketitle

\vspace{-1em}

\begin{abstract}
In two influential contributions, \citet{rosenbaum2005exact,rosenbaum2010design} advocated for using the distances between component-wise ranks, instead of the original data values, to measure covariate similarity when constructing matching estimators of average treatment effects. While the intuitive benefits of using covariate ranks for matching estimation are apparent, there is no theoretical understanding of such procedures in the literature. We fill this gap by demonstrating that Rosenbaum's rank-based matching estimator, when coupled with a regression adjustment, enjoys the properties of double robustness and semiparametric efficiency without the need to enforce restrictive covariate moment assumptions. Our theoretical findings further emphasize the statistical virtues of employing ranks for estimation and inference, more broadly aligning with the insights put forth by Peter Bickel in his 2004 Rietz lecture \citep{bickel2004}.
\end{abstract}


\section{Introduction}

Consider the problem of estimating the average treatment effect (ATE),
\[
\tau:=\E\Big[Y(1)-Y(0)\Big],
\]
based on an observational study encompassing $n$ observations of a binary treatment $D\in\{0,1\}$, some measured pre-treatment covariates $X\in\mathbb{R}^d$, and an outcome $Y=Y(D)\in\mathbb{R}$ that is realized from the two potential outcomes $(Y(0),Y(1))$. Among the techniques employed to estimate $\tau$, nearest neighbor (NN) matching stands as one of the most widely adopted and comprehensible approaches \citep[see][and references therein]{stuart2010matching}. These estimators aim to impute the missing potential outcome of each unit in one treatment group by finding units from the opposite treatment group whose covariate profile closely resemble that of the unit with the missing potential outcome. The quantification of covariate similarity relies on the user-specified {\it distance metric} between the (distribution of the) covariates of each unit.

\cite{abadie2006large} laid out the mathematical groundwork to study NN matching estimators employing the Euclidean distance metric, while \cite{abadie2011bias} established a root-$n$ central limit theorem for a bias-corrected version of those NN matching estimators. More recently, \cite{lin2021estimation} established connections between NN estimators and augmented inverse probability weighted (AIPW) methods \citep{robins1994estimation, scharfstein1999adjusting}, thereby establishing double robustness and semiparametric efficiency theory for NN matching estimators when the number of matches diverges to infinity with the sample size.

The aforementioned theoretical work, however, centers around the Euclidean distance metric for determining NN matches. This approach may exhibit sensitivity to alterations in scale and to the presence of extreme outliers or heavy-tailed distributions. Indeed, all the existing theoretical results on matching require a compact supported covariates, which is theoretically hard to alleviate, if not impossible. On the other hand, in practice, distance metrics are often derived from a ``standardized'' representation of the data, and the selection of a distance metric is an important factor in causal inference because various metrics can lead to different conclusions \citep[Chapter 9]{rosenbaum2010design}. 

This paper focuses on a particular standardization approach that identifies NNs by measuring the Euclidean distance between the component-wise ranks of the covariates $X$, as proposed for the celebrated Rosenbaum's rank-based matching estimator \citep[Chapter 9.3]{rosenbaum2010design}. The concept of rank-based standardization is straightforward to interpret, easy to implement, and computationally efficient, while also being scale-invariant and insensitive to heavy-tailed distributions. Furthermore, due to their data adaptivity, rank-based methods are often used in treatment effect settings such as for analysis of experiments \citep{rosenbaum2010design}, regression discontinuity plots \citep{calonico2015optimal}, and binscatter regressions \citep{cattaneo2023binscatter}.

The challenge in formally studying rank-based methods lies on the theoretical side, as the transformation of the covariates into their ranks disrupts the independence structure of the original data, and thus complicates the subsequent statistical analysis. Our main theorem (Theorem \ref{thm:bc}) offers the {\it first} theoretical analysis of Rosenbaum's rank-based matching estimator, elucidating many of its appealing properties when combined with regression adjustments. In particular, our theory not only confirms Rosenbaum's intuition that the rank-based distance can limit the influence of outliers and heavy-tailed distributions \citep[page 210]{rosenbaum2010design}, but also demonstrates that the rank-based matching estimator can be doubly robust and semiparametrically efficient, particularly without imposing restrictive moment assumptions on the distribution of the covariates. More broadly, 
they align with Peter Bickel's 2004 Rietz lecture advocating for ``standardization by ranks'' when performing statistical and machine learning related tasks \citep{bickel2004}.

Our paper further offers two technical contributions as a by-product of the analysis, which may be of independent interest. First, Theorem \ref{thm:bc} establishes consistency, asymptotic linearity, and semiparametric efficiency of Rosenbaum's Rank-based Matching Estimator under generic high-level conditions on the regression adjustment. The proof of that theorem relies on a careful combination of empirical process theory for rank-based statistics and tools established in \cite{lin2022regression} and \cite{lin2021estimation} for matching estimators involving a growing number of nearest neighbors, which are generalized herein to accommodate standardization/transformation functions, of which component-wise ranking is one particular example. Second, Theorem \ref{thm:series} presents novel mean square and uniform convergence rates for series estimators when the covariates are generated via possibly unknown functions of the original independent variables, of which component-wise ranking is one particular example. Those results are proven for general series estimators \citep{newey1997convergence,belloni2015some} with covariate-generated conditioning variables, and thus lead to suboptimal uniform approximation results, but we also discuss how they may be upgraded to deliver optimal uniform convergence rates for partition-based series estimators \citep{huang03-local-asympt-for-polyn-spline-regres,cattaneo2013optimal,cattaneo2020large,cattaneo2023binscatter}.

\vspace{0.2cm}

\noindent{\bf Paper organization.} The rest of this paper is organized as follows. In Section \ref{sec:setup} we describe the problem setup and introduce the studied rank-based matching estimator. Section \ref{sec:theory} establishes the main result of this paper under general regression adjustments. Section \ref{sec:regression} discusses primitive conditions for regression adjustment using series estimation. Section \ref{sec:extend} generalizes the results in Section \ref{sec:theory}. All proofs are in the appendix.

\section{Setup}\label{sec:setup}

We adopt the standard potential outcomes causal model for a binary treatment, where it is assumed that there are $n$ independent and identically (i.i.d.) distributed realizations $\{X_i,D_i,Y_i(0), Y_i(1)\}_{i=1}^n$, of a quadruple $(X,D,Y(0),Y(1))$. In practice, we are only able to observe a part of the data, i.e., $\{X_i,D_i,Y_i:=Y_i(D_i)\}_{i=1}^n$. The goal is to conduct estimation and inference for the population ATE, 
\[
\tau=\E[Y(1)-Y(0)],
\]
based only on the observed data.

Rosenbaum's rank-based matching approach estimates $\tau$ by plugging the component-wise ranks of $X_i$'s, instead of the original values, into the NN matching mechanism, with each unit matched to $M$ units in the opposite treatment group with replacement. It proceeds in three step as follows.

\vspace{0.2cm}

\noindent{\bf Step 1.}  Given a sample $\{(X_i, D_i, Y_i)\}_{i=1}^n$ with $X_i = (X_{i,1},\ldots,X_{i,d})^\top$, introduce the vector of marginal empirical cumulative distribution functions (CDFs), $\hat{\mF}_n(\cdot):\bR^d \to [0,1]^d$, such that for any input $x = (x_1,\ldots,x_d)^\top \in \bR^d$,
\begin{align*}
  \hat{\mF}_n(x) :=  (\hat{F}_{n,1}(x_1),\ldots,\hat{F}_{n,d}(x_d))^\top,~~{\rm with}~~\hat{F}_{n,k}(x_k) := \frac{1}{n}\sum_{i=1}^n \ind(X_{i,k} \le x_k),~k\in\zahl{d}.
\end{align*}
Here $\ind(\cdot)$ stands for the indicator function and $\zahl{d}:=\{1,2,\ldots,d\}$. For each $i\in\zahl{n}$, introduce $\hat U_i:=\hat{\mF}_n(X_i)\in\bR^d$. It is understood that each component of $n\hat U_i$ represents the corresponding rank of $X_{i,\cdot}$ among $X_{1,\cdot},\ldots,X_{n,\cdot}$. The vector of marginal population CDFs is $\mF(\cdot):\bR^d \to [0,1]^d$, such that for any input $x = (x_1,\ldots,x_d) \in \bR^d$,
\begin{align*}
  \mF(x) :=  (F_1(x_1),\ldots,F_d(x_d))^\top,~~{\rm with}~~F_k(x_k) := \P(X_{1,k} \le x_k),~k\in\zahl{d}.
\end{align*}
Define $U:=\mF(X)\in[0,1]^d$ and for each $i\in\zahl{n}$, $U_i:=\mF(X_i)\in[0,1]^d$.

\vspace{0.2cm}
\noindent{\bf Step 2.} To correct for the estimation bias from matching, regression adjustment is employed. Let $\hat{\mu}_0(\cdot)$ and $\hat{\mu}_1(\cdot)$ be mappings from $[0,1]^d$ to $\bR$ such that they separately estimate the conditional means of the outcomes,
\[
\mu_0(u) := \E [Y \given U=u,D=0]\qquad \text{and} \qquad \mu_1(u) := \E [Y \given U=u,D=1]. 
\]
We obtain $\hat{\mu}_0$ and $\hat{\mu}_1$ by regressing $Y_i$'s in either the treatment or control group on the corresponding $\{\hat{U}_i\}$'s, respectively.

\vspace{0.2cm}
\noindent{\bf Step 3.} Implement bias-corrected nearest neighbor matching on $(\hat U_i, D_i,Y_i)$'s. Specifically, let $\cJ(i)$ represent the index set of the $M$-NNs of $\hat{U}_i$ in $\{\hat{U}_j:D_j=1-D_i\}_{j=1}^n$, measured using the Euclidean metric $\|\cdot\|$ with {\it ties broken in arbitrary way}. The final rank-based matching estimator is
\begin{align}\label{eq:tau}
  \hat{\tau} := \frac{1}{n} \sum_{i=1}^n \Big[\hat{Y}_i(1) -\hat{Y}_i(0)\Big],
\end{align}
where, for $\omega\in\{0,1\}$,
\begin{align*}
    \hat{Y}_i(\omega) := \begin{cases}
        Y_i, & \mbox{ if } D_i=\omega,\\
        \frac{1}{M} \sum_{j \in \cJ(i)} (Y_j + \hat{\mu}_\omega(\hat{U}_i) - \hat{\mu}_\omega(\hat{U}_j)), & \mbox{ if } D_i=1-\omega   
    \end{cases}.
\end{align*}

\section{Main Result}\label{sec:theory}

Following the notation system of \cite{abadie2006large,abadie2011bias}, let $K(i)$ represent the number of matched times for each unit $i$, i.e.,
\[
  K(i) := \sum_{j=1, D_j = 1-D_i}^n \ind\big(i \in \cJ(j)\big).
\]
In this paper, however, $K(i)$ denotes the matched times according to {\it the rank-based distance}, not the original Euclidean distance. Following \cite{lin2022regression} and \cite{lin2021estimation}, the rank-based bias-corrected matching estimator in \eqref{eq:tau} can be represented as an AIPW estimator:
\begin{align}\label{eq:mbc}
	\hat{\tau} = \hat{\tau}^{\rm reg} + \frac{1}{n} \sum_{i=1}^n (2D_i-1)\Big(1 + \frac{K(i)}{M}\Big) \hat{R}_i,
\end{align}
where 
\[
\hat\tau^{\rm reg}:=\frac1n\sum_{i=1}^n[\hat\mu_1(\hat U_i)-\hat\mu_0(\hat U_i)]~~{\rm and}~~\hat R_i:=Y_i-\hat\mu_{D_i}(\hat U_i), \text{ for } i \in \zahl{n}.
\]
We employ this insight throughout our large sample distributional analysis of $\hat\tau$.

To establish our main theorem we impose some assumptions. The first assumption is posed to regulate basic features of the data generating distribution, in particular making the estimation problem identifiable. Conceptually, the main difference with prior literature is that the assumption concerns the triple $(U_i,D_i,Y_i)$'s, that is, $X_i$ is replaced by $U_i$, the scaled population rank.

\begin{assumption} \label{asp:dr}
	\phantomsection $ $
	\begin{enumerate}[itemsep=-.5ex,label=(\roman*)]
	    \item For almost all $u \in [0,1]^d$, $D$ is independent of $(Y(0),Y(1))$ conditional on $U=u$, and there exists a fixed constant $c > 0$ such that $c < e(u):=\P(D=1 \given U=u) < 1-c$.
	    \item $[(X_i,D_i,Y_i)]_{i=1}^n$ are i.i.d. following the joint distribution of $(X,D,Y)$.
	    \item $\E [(Y(\omega) - \mu_\omega(U))^2 \given U=u] $ is uniformly bounded for almost all $u \in [0,1]^d$ and $\omega =0,1$.
	    \item $\E [\mu^2_\omega(U)]$ is bounded for $\omega = 0,1$.
	\end{enumerate}
\end{assumption}

The second assumption ensures that the regression adjustment procedure is, at least, well-posited in the sense that the estimator $\hat{\mu}_\omega(x)$ is uniformly consistent for some well-behaved (possibly misspecified) conditional expectation. Let $\|\cdot\|_{\infty}$ be the $L^{\infty}$ function norm.

\begin{assumption}\label{asp:dr-p}
	\phantomsection $ $
	For $\omega = 0,1$, there exists a deterministic, possibly changing with $n$, continuous function $\bar{\mu}_\omega(\cdot):[0,1]^d \to \bR$ such that $\E [\bar{\mu}^2_\omega(U)]$ is uniformly bounded and the estimator $\hat{\mu}_\omega(x)$ satisfies $\Vert \hat{\mu}_\omega - \bar{\mu}_\omega \Vert_\infty = o_\P(1)$.
\end{assumption}

The next assumption regulates the population rank-transformed random vector $U$, which identifies the copula distribution for $X$ \citep{joe2014dependence}. This assumption requires $X$ to be continuous, but discrete components of $X$ can be easily handled by conditioning \citep{stuart2010matching}.

\begin{assumption}\label{asp:dr-p'}
	\phantomsection $ $
	The Lebesgue density of $U$ exists and is continuous over its support. 
\end{assumption}

The next three assumptions concern the case of consistent population regression functions $\mu_0(U)$ and $\mu_1(U)$, and will be used in contrast to Assumption \ref{asp:dr-p}, where the regression adjustment procedure is allowed to be inconsistent for the population ranked-based regression functions. More precisely, Assumption \ref{asp:dr-p} vis-\'a-vis Assumptions \ref{asp:dr-o}--\ref{asp:se2} are used for establishing the double robustness and semiparametric efficiency of $\hat\tau$, respectively. Using standard multi-index notation, let $\Lambda_k$ is the set of all $d$-dimensional vectors of nonnegative integers $t=(t_1,\ldots,t_d)$ such that $|t|=\sum_{i=1}^d t_i = k$ with $k$ any positive integer, and $\partial^t \mu_{\omega}$ denotes the corresponding partial derivative of $\mu_{\omega}$.

\begin{assumption}\label{asp:dr-o}
	\phantomsection $ $
	For $\omega =0,1$, $\mu_\omega$ is continuous and the estimator $\hat{\mu}_\omega(x)$ satisfies $\Vert \hat{\mu}_\omega - \mu_\omega \Vert_\infty = o_\P(1)$.
\end{assumption}

\begin{assumption} \label{asp:se1}
	\phantomsection $ $
	\begin{enumerate}[itemsep=-.5ex,label=(\roman*)]
		\item $\E [(Y(\omega) - \mu_\omega(U))^2 \given U=u] $ is uniformly bounded away from zero for almost all $u \in [0,1]^d$ and $\omega = 0,1$.
		\item There exists a constant $c>0$ such that $\E [\vert Y(\omega) - \mu_\omega(U) \vert ^{2+c} \given U=u]$ is uniformly bounded for almost all $u \in [0,1]^d$ and $\omega = 0,1$.
		\item\label{asp:se1-3} $\max_{t \in \Lambda_{\max\{\lfloor d/2 \rfloor,1\} + 1}} \Vert \partial^t \mu_{\omega} \Vert_\infty$ is bounded, where $\lfloor \cdot \rfloor$ stands for the floor function.
	\end{enumerate}
\end{assumption}

\begin{assumption} \label{asp:se2}
	\phantomsection $ $ 
	For $\omega = 0,1$, the estimator $\hat{\mu}_\omega(x)$ satisfies
	\[
	  \max_{t \in \Lambda_{\max\{\lfloor d/2 \rfloor,1\} + 1}} \Vert \partial^t \hat{\mu}_{\omega} \Vert_\infty = O_\P(1)~~~{\rm and}~~~ 
	  \max_{t \in \Lambda_\ell} \Vert \partial^t \hat{\mu}_{\omega} - \partial^t \mu_{\omega} \Vert_\infty = O_\P(n^{-\gamma_\ell}) ~~\mbox{\rm for all}~~ \ell \in \zahl{\max\{\lfloor d/2 \rfloor,1\}},
	\]
	with some constants $\gamma_\ell$'s satisfying $\gamma_\ell > \max\{\frac{1}{2} - \frac{\ell}{d},0\} $ for $\ell=1,2,\ldots,\max\{\lfloor d/2 \rfloor,1\}$.
\end{assumption}

Several remarks on the above assumptions align with our conceptual discussion. First, as the quantile transformation preserves all information, Assumption \ref{asp:dr} is either equivalent to, or weaker than, the standard assumptions in the matching estimation literature. Second, Assumption \ref{asp:dr-p} accommodates regression model misspecification, and its validity may be verified by leveraging the fundamental projection principles underlying regression techniques (see, also, the discussions in Section \ref{sec:regression}). Third, Assumption \ref{asp:dr-p'} constitutes a mild condition, notably satisfied by distribution families such as the Gaussian (copula) and Cauchy (copula). Lastly, Assumptions \ref{asp:dr-o} through \ref{asp:se2} merit more discussion: due to the shift from using $X_i$'s as inputs to $\hat U_i$'s in the regression function, direct verification using standard results from the nonparametric smoothing estimation literature is no longer possible. We return to this technical issue in Section \ref{sec:regression}, where we consider explicitly least squares series estimation \citep{newey1997convergence,huang03-local-asympt-for-polyn-spline-regres,cattaneo2013optimal,belloni2015some,cattaneo2020large} to illustrate the verification of Assumption \ref{asp:dr-p} and Assumptions \ref{asp:dr-o}--\ref{asp:se2}.

We are now ready to present our main theorem for Rosenbaum's rank-based matching estimator.

\begin{theorem}[Main Theorem]  \phantomsection $ $ \label{thm:bc} 
	
\begin{enumerate}[itemsep=-.5ex,label=(\roman*)]
	
  \item\label{thm:dr} (Double robustness of $\hat\tau$) If either Assumptions~\ref{asp:dr}, \ref{asp:dr-p}, \ref{asp:dr-p'}, $M\log n/n \to 0$, and $M \to \infty$ as $n \to \infty$ hold, or Assumptions~\ref{asp:dr} and \ref{asp:dr-o} hold, then 
  \[
  \hat{\tau} - \tau \stackrel{\sf p}{\longrightarrow} 0.
  \]
  
  \item\label{thm:mbc} (Semiparametric efficiency of $\hat\tau$) Assume Assumptions~\ref{asp:dr}, \ref{asp:dr-p'}, \ref{asp:se1}, \ref{asp:se2} hold. Define
  \[
  \gamma  =  \max \Big\{ \Big[1- \frac{1}{2} \frac{d}{\max\{\lfloor d/2 \rfloor,1\} + 1}\Big] , \min_{\ell \in \zahl{\max\{\lfloor d/2 \rfloor,1\}} } \Big[1-\Big(\frac{1}{2} - \gamma_\ell \Big)\frac{d}{\ell}\Big] \Big\},
  \]
  where the $(\gamma_\ell:\ell=1,2,\ldots,\max\{\lfloor d/2 \rfloor,1\})$ are introduced in Assumption \ref{asp:se2}. If $M \to \infty$ and $M/n^\gamma \to 0$ as $n \to \infty$, then
  \[
  \sqrt{n} (\hat\tau - \tau) \stackrel{\sf d}{\longrightarrow} N(0,\sigma^2),
  \]
  where 
  \[
  \sigma^2:=\E\Big[\mu_1(U)-\mu_0(U)+\frac{D(Y-\mu_1(U))}{e(U)}-\frac{(1-D)(Y-\mu_0(U))}{1-e(U)}-\tau \Big]^2.
  \]

  \item\label{thm:asyvar} (Variance Estimation) If the conditions in part \ref{thm:mbc} and Assumption~\ref{asp:dr-o} hold, then
  \[\hat\sigma^2 - \sigma^2 \stackrel{\sf p}{\longrightarrow} 0\]
  where
  \[\hat\sigma^2:=\frac{1}{n}\sum_{i=1}^n\Big[\hat\mu_1(\hat U_i)-\hat\mu_0(\hat U_i)+(2D_i-1)\Big(1+\frac{K_M(i)}{M}\Big)\hat R_i-\hat\tau\Big]^2.
  \]
\end{enumerate}
\end{theorem}  

This theorem establishes three main results. Part \ref{thm:dr} shows that the generic bias-corrected Rosenbaum's rank-based matching estimator is doubly robust for a fairly large class of regression estimators based on estimated ranks of the covariates. The result in part \ref{thm:mbc} gives general regularity conditions guaranteeing asymptotic normality of the estimator. It follows directly from the second result that the estimator is semiparametrically efficient for estimating the ATE \citep{hahn1998role}. Finally, part \ref{thm:asyvar} establishes consistency of the plug-in variance estimator under nearly the same conditions as required for consistency and asymptotic normality.

\section{Regression Adjustment using Series Least Squares}\label{sec:regression}

The only remaining issue concerning Theorem \ref{thm:bc} revolves around the use of pairs $(\hat U_i,Y_i)$ for bias correction, as opposed to $(X_i,Y_i)$, or the idealized oracle $(U_i,Y_i)$ pairs. The dependence among the estimated rank-adjusted $\hat U_1,\dots,\hat U_n$ poses a challenge, making it hard to apply existing results in nonparametric statistics for the direct verification of Assumption \ref{asp:dr-p} or Assumptions \ref{asp:dr-o}--\ref{asp:se2}. This section illustrates how these assumptions can be verified when using series least squares regression estimation, covering both canonical approximating functions (e.g., power series, fourier series, splines, wavelets, and piecewise polynomials) as well as general covariate transformations (e.g., high-dimensional least squares regression with structured regressors).

The main result in this section concerns general estimated transformations of the independent variables, based on the underlying regressors only, and allowing for possible misspecification in both fixed-dimension and increasing-dimension least squares regression settings. Thus, we consider a more general setup where we either observe or have approximate information about $n$ i.i.d. pairs $(Y_1,W_1),\dots,(Y_n,W_n)$ of $(Y,W)$, and the goal is to estimate the conditional expectation
\[\psi(w):=\E[Y\given W=w],\]
using only the outcome variables $\{Y_i\}_{i=1}^n$ and the {\it generated covariates} $\{\hat{W}_i\}_{i=1}^n$ that are ``approximately close'' to $\{W_i\}_{i=1}^n$, where $\hat{W}_i$'s are measurable with respect to some sigma field $\cF_n$, $n\geq1$. In practice, as it is the case in the rank-based transformation we consider in this paper, a useful choice is $\cF_n=\cS(W_1,W_2,\dots,W_n)$, where $\cS(Z)$ denotes the sigma field generated by the random variable $Z$.

Let $p_K(w)=(p_{1K}(w),\ldots,p_{KK}(w))^\top$ be a $K$-dimensional vector of basis functions so that their linear combination may approximate $\psi(\cdot)$ well when $K$ is sufficiently large, at least under some specific assumptions. However, a good approximation is not strictly required, as we also consider misspecified regression adjustments. This is important in the context of this paper because Theorem \ref{thm:bc} established a double-robust property of the regression adjusted Rosenbaum's rank-based matching estimator, thereby allowing for misspecified or inconsistent estimators $(\hat \mu_{0},\hat \mu_{1})$ of the regression functions $(\mu_{0},\mu_{1})$. Furthermore, we also allow for the possibility of $W$ having a Lebesgue density that may not be bounded away from zero, as it occurs when its support is unbounded. In our application, due to the rank transformation, the support of the rank-based covariates is bounded but their density may not be bounded away from zero in, e.g., the Gaussian copula case.

The following assumption summarizes our setup. Let $\lambda_{\min}(\cdot)$ be the smallest eigenvalue of the input matrix.

\begin{assumption} \label{asp:series}
	\phantomsection $ $
	\begin{enumerate}[itemsep=-.5ex,label=(\roman*)]
		\item $(Y_1,W_1),\dots,(Y_n,W_n)$ are i.i.d. draws from $(Y,W)\in\mathcal{Y}\times\mathcal{W}\subseteq\mathbb{R}\times\mathbb{R}^{d}$.
		\item $\E[Y^2]$ is bounded, and $\lambda_{\min}\big(\E[p_K(W)p_K(W)^\top]\big)>0$ for all $K$ and $n$.
		\item There exists a sequence of sigma fields $\{\cF_n\}$ such that $(\widehat{W}_1,\widehat{W}_2,\dots,\widehat{W}_n)$ is measurable with respect to $\cF_n$ for all $n$, $\sup_{i\leq n}\E[(Y_i-\psi(W_i))^2|\cF_n]$ is bounded uniformly in $n$, and $\E[(Y_i-\psi(W_i))(Y_j-\psi(W_j))|\cF_n]=0$ for all $i\ne j$ and $n$.
	\end{enumerate}
\end{assumption}

The series estimator with generated covariates is
\[
\hat{\psi}_K(w) = p_K(w)^\top\hat{\beta}_K, \qquad \hat{\beta}_K \in \argmin_{b \in \bR^K} \frac{1}{n}\sum_{i=1}^n (Y_i - p_K(\hat{W}_i)^\top b)^2,
\]
which gives $\hat{\beta}_K := (\mP_n^\top \mP_n)^{-}\mP_n^\top \mY$, where $\mP_n = (p_K(\hat{W}_1),\ldots,p_K(\hat{W}_n))^\top\in\mathbb{R}^{n\times k}$, $\mY := (Y_1,\ldots,Y_n)^\top$, and $\mA^-$ denotes a generalized inverse of the matrix $\mA$. Let $\mP := (p_K(W_1),\ldots,p_K(W_n))^\top$ be what $\mP_n$ shall approximate, and 
\[
\beta_K := \argmin_{b \in \bR^K} \E[(Y_1-p_K(W_1)^\top b)^2] = \argmin_{b \in \bR^K} \E[(\psi(W_1)-p_K(W_1)^\top b)^2] 
\]
be what $\hat{\beta}_K$ shall approximate. It follows that $\psi_K(w) := p_K(w)^\top \beta_K$ is the best $L^2$ approximation of $\psi(w)$ based on $p_K(w)$, where $\beta_K = \mQ^{-}\E[p_K(W_1)\psi(W_1)]$ with $\mQ :=\E[p_K(W_1) p_K(W_1)^\top]\in\mathbb{R}^{K\times K}$.

Our results rely on the following quantities characterizing different aspects of the series estimator and the approximation errors:
\begin{align*}
	&\lambda_K        :=\lambda_{\min}(\mQ),
	&&\zeta_{q,K}     := \max_{t \in \Lambda_q} \sup_{w\in\mathcal{W}} \Vert \partial^t p_K(w) \Vert,\\
	&\xi^2_{K}       := \E[(\psi(W_1)-\psi_K(W_1))^2],
	&&\vartheta_{q,K} := \max_{t \in \Lambda_q}  \Vert \partial^t \psi - \partial^t \psi_K \Vert_{\infty},
\end{align*}
and
\begin{align*}
	R_n := \Vert \mPsi - \mPsi_n \Vert^2/n,\qquad
	B_n := \Vert (\mP - \mP_n)\mQ^{-1/2} \Vert_2^2/n,
\end{align*}
where $\mPsi := (\psi(W_1),\ldots,\psi(W_n))^\top\in\mathbb{R}^n$, $\mPsi_n := (\psi(\hat{W}_1),\ldots,\psi(\hat{W}_n))^\top\in\mathbb{R}^n$, and $\Vert \cdot \Vert_2$ denotes the matrix spectral norm.

Let $\Vert g \Vert^2_{L^2} := \int |g(w)|^2 \d F_W(w)$, and consider first the $L^2$ rate of approximation of the series-based least squares estimator:
\begin{align*}
	\Vert \hat{\psi}_K-\psi \Vert^2_{L^2}
	\leq 2 \Vert \hat{\psi}_K - \psi_K \Vert^2_{L^2} + 2 \Vert \psi_K - \psi\Vert^2_{L^2},
\end{align*}
where it follows immediately that $\Vert \psi_K - \psi\Vert^2_{L^2} = \xi^2_{K} \leq \vartheta^2_{0,K}$, implying that the best mean square approximation $\psi_K(w) = p_K(w)^\top \beta_K$ will approximate $\psi(w)$ well if $\vartheta_{0,K}\to0$, or at least $\xi_{K}\to0$, which in turn requires $K\to\infty$ in general. However, in many applications the series estimator may be misspecified or inconsistent in the sense that $\xi_{K}\not\to0$. In those cases, it is natural to take $\psi_K(w)$ as the target ``parameter''. The following theorem establishes two distinct $L^2$ convergence rates for the series estimator relative to the latter quantity.

\begin{theorem}[$L_2$ Convergence]\label{thm:series}
	
	Let Assumption \ref{asp:series} hold, $\lambda^{-1}_K \zeta^2_{0,K} \log(K)/n \to 0$, and $B_n \to 0$. Then,
	\begin{align*}
		\Vert \hat{\psi}_K - \psi_K \Vert^2_{L^2}
		 = O_\P\Big(\frac{K}{n} + \mathsf{A}_n \Big),
	\end{align*}
	where the approximation error term can be taken to be either
	\begin{align*}
		\mathsf{A}_n = \min\Big\{B_n + \xi^2_{K}  \,,\, R_n + \vartheta_{0, K}^2\Big\},
	\end{align*}
	or
	\begin{align*}
		\mathsf{A}_n = B_n + B_n \min\Big\{B_n + \xi^2_{K}  \,,\, R_n + \vartheta_{0, K}^2\Big\} + \min\Big\{\lambda^{-1}_K \zeta_{0,K}^2 \xi^2_{K}/n \,,\, K \vartheta_{0, K}^2/n\Big\}.
	\end{align*}
\end{theorem}

This theorem provides new results relative to previously known mean square convergence rates for series estimation. More specifically, it allows for generated regressors based on covariates with a possibly vanishing minimum eigenvalue of the expected scaled Gram matrix ($\lambda_K$), as it may occur when the Lebesgue density of $W$ is positive but not bounded away from zero on $\mathcal{W}$. Furthermore, the second rate estimate allows for a non-vanishing $L^2$ approximation error ($\vartheta_{0,K}\geq\xi_K\not\to0$), thereby offering $L^2$ consistency results for general least squares approximations.

It is easy to deduce (suboptimal) uniform rates of approximation using Theorem \ref{thm:series} because
\begin{align*}
	\max_{t \in \Lambda_q} \Vert \partial^t \hat{\psi}_K - \partial^t \psi \Vert_\infty
	&\le \max_{t \in \Lambda_q} \Vert \partial^t p_K^\top (\hat{\beta}_{K}-\beta_K) \Vert_\infty + \max_{t \in \Lambda_q} \Vert \partial^t \psi_K - \partial^t \psi \Vert_\infty,\\
	&\leq \zeta_{q,K} \lambda_K^{-1/2} \Vert \hat{\psi}_K - \psi_K \Vert_{L^2} + \vartheta_{q, K},
\end{align*}
where, as noted before, the first term characterizes the error in approximation when perhaps $\vartheta_{q, K}\not\to0$, in which case the target ``parameter'' can be taken to be $\partial^t \psi_K(w) = \partial^tp_K(w)^\top\beta_K$, regardless of whether $K\to\infty$ or not.

Underlying the assumptions imposed in Theorem \ref{thm:series}, there are several parameters that need further discussion. From the standard series estimation literature, $\zeta_{q,K}=O(K^{1+q})$ for power series and $\zeta_{q,K} = O(K^{1/2+q})$ for fourier series, splines, compact supported wavelets, and piecewise polynomial regression. Lower bounds for $\lambda_K$ need to be established on a case-by-case basis when the density of $W$ is not assumed to be bounded away from zero, so we illustrate one such verification further below for the case of rank transformations and a Gaussian copula. As already mentioned, the parameter $\xi_K \leq \vartheta_{0,K}$ captures the degree of approximation (or misspecification) of the series regression estimator, and needs not to vanish, in which case the second rate result in Theorem \ref{thm:series}, and the implied uniform convergence rate, must be used. If $\psi$ is $s$-times differentiable and other regularity conditions hold, then $\xi_{K} = O(K^{-s/d})$ for all the usual approximation basis functions if appropriately specified. In general, however, the difference between the $L^2$ and $L^\infty$ approximaton errors, $\xi_K$ and $\vartheta_{0,K}$, depends on the basis functions employed and the data generating features. In particular, for instance, when employing locally supported basis functions, it can be verified that $\xi_K \asymp \vartheta_{0,K}$, in which case $\vartheta_{q, K} = O(K^{-(s-q)/d})$ under regularity conditions. See \cite{newey1997convergence}, \cite{huang03-local-asympt-for-polyn-spline-regres}, \cite{cattaneo2013optimal}, \cite{belloni2015some}, \cite{cattaneo2020large}, and references therein, for more details.
			
An important feature of Theorem \ref{thm:series} is that it allows for generated regressors based on the covariates, which introduces two additional quantities characterizing the approximation rate: $R_n$ and $B_n$. For example, if $\psi$ satisfies the Lipschitz condition $|\psi(a)-\psi(b)|\leq L \lVert a-b \rVert$ for some constant $L$, then
\[R_n = \frac{1}{n} \sum_{i=1}^n (\psi(\hat{W}_i) - \psi(W_i))^2
      \leq L^2 \max_{i =1,2,\dots,n} \lVert \hat{W}_i-W_i\rVert^2 = O_\P(r_n),
\]
and therefore the convergence rate of $R_n$ is determined by the uniform convergence rate of the transformation of the covariates. Recall that in our application, we consider the empirical rank transformation $(W_i,\hat{W}_i)=(U_i,\hat{U}_i)$, and therefore $r_n=1/n$. A similar calculation can be done to bound $B_n$ when $p_K(\cdot)$ is smooth, because
\begin{align*}
	B_n \leq \lambda_K^{-1} \frac{1}{n} \sum_{i=1}^n \Vert p_K(\hat{W}_i)-p_K(W_i) \Vert^2
	    \leq \lambda_K^{-1} L_K^{2} \max_{i =1,2,\dots,n} \lVert \hat{W}_i-W_i \rVert^2 = O_\P(b_n),
\end{align*}
provided that the basis functions are Lipschitz with constant $L_K$, and where $L_K^{2} = O(\zeta^2_{1,K})$. Thus, in the particular case of the empirical rank transformation, $b_n = \lambda_K^{-1} \zeta^2_{1,K}/n$.

It remains to illustrate how to lower bound the minimum eigenvalue $\lambda_K$. It is well-known that if $W$ admits a Lebesgue density bounded and bounded away from zero over the support $\mathcal{W}$, then $\lambda^{-1}_{K}$ is uniformly bounded in $K$, after possibly rotating the basis functions. The following lemma considers the more interesting case when the density is not bounded away from zero, as it occurs for example when the support of $W$ is unbounded.

\begin{lemma}[Lower Bound on $\lambda_K$]\label{lemma:eigen}
  Suppose that $W$ admits a Lebesgue density $f_W$, and $p_{1K},\ldots,p_{KK}$ are orthonormal with respect to the Lebesgue measure over the support of $W$. If there exists a universal constant $C>0$ such that, for all sufficiently small $t>0$, $\mathsf{Leb}(\{w:0<f_W(w)< t\}) \le C t^{\rho}$ for some $\rho>0$, then
  $\lambda^{-1}_K = O(\zeta_{0,K}^{2/\rho})$.

\end{lemma}

As expected, the additional conditions in the previous lemma restrict the tail of $f_W$. It remains to illustrate how to verify the result for the case when $W_i$ and $\hat W_i$ are taken to be $U_i$ and $\hat U_i$ as in Section \ref{sec:setup}. This is done in the following proposition for the case of a Gaussian copula.

\begin{proposition}[Sufficient Conditions for Gaussian Copula] \label{crl:series}
	\phantomsection $ $
	Suppose $U_1,U_2,\dots,U_n$ follow a Gaussian copula distribution with an invertible parameter correlation matrix $\Sigma$. Then, the conditions of Lemma~\ref{lemma:eigen} are satisfied with $\rho = \lambda^{-1}_{\max}\big(\Sigma^{-1}-\mI_d\big)/d$, where $\mI_d$ is the d-dimensional identity matrix. 
\end{proposition}

Using Theorem \ref{thm:series} in general, or Lemma \ref{lemma:eigen} and specific conditions on the underlying copula of the population rank-based transformed covariates, it is possible to verify the conditions of Theorem \ref{thm:bc}. To be more specific, Theorem \ref{thm:series} gives versatile (albeit suboptimal) uniform consistency rates for series estimators $(\hat \mu_{0},\hat \mu_{1})$ of either some approximate (misspecified) functions or the population functions $(\mu_{0},\mu_{1})$. It is worth noting that our verification aimed for generality in terms of high-level conditions, but for the case of partition-based (locally supported) series estimators it is possible to obtain better (in fact optimal in some cases) uniform convergence rates \citep{cattaneo2013optimal,belloni2015some,cattaneo2020large}. Specifically, \citet{cattaneo2023binscatter} studies the case of rank-based transformations for B-splines when $d=1$, and establishes optimal uniform convergence rates on compact support. Their results could be extended to obtain sharper uniform convergence rates with generated regressors based on the covariates as required by Theorem \ref{thm:bc}.

\section{Generalized Framework}\label{sec:extend}

This section extends Rosenbaum's rank standardization idea to a more general setting, and establishes general theory that will cover Theorem \ref{thm:bc} as a special case. More specifically, for $\omega =0,1$, consider the following general mappings
\[
\phi_\omega: \cX \to \cX_\phi \subset \bR^m,
\]
with $\cX$ representing the support of $X$ and $m$ not necessarily equal to $d$. Note that here we allow $\phi_0$ and $\phi_1$ to be different. Consider the setting when $\phi_\omega$ is possibly unknown, and we will approximate it based on the sample $\{(X_i, D_i, Y_i)\}_{i=1}^n$, leading to a generic estimator $\hat{\phi}_\omega$ that may differ with different $\omega$. We then define 
\[
U_{\phi,\omega}:= \phi_\omega(X)\qquad \text{and}\qquad \hat{U}_{\phi,\omega,i} := \hat{\phi}_\omega(X_i)\qquad \text{for}\quad i \in \zahl{n}.
\] 
Note that, when setting $\phi_0=\phi_1=\mF$ and $\hat\phi_0=\hat\phi_1=\hat\mF_n$, we recover the $U$ and $\hat U_i$'s introduced in Section \ref{sec:setup}.

Similar to Section \ref{sec:setup}, let $\cJ_{\phi}(i)$ represent the index set of $M$-NN matches of $\hat{U}_{\phi,1-D_i,i}$ in $\{\hat{U}_{\phi,1-D_i,j}:D_j=1-D_i\}_{j=1}^n$ with ties broken in an arbitrary way. In other words, for determining the nearest neighbors, we are going to measure the similarity based on the Euclidean distance between transformed data points with the transformation function probably also having to be learned from the same data. Additionally, let $\hat{\mu}_{\phi,\omega}(u)$ be a mapping from $\cX_\phi$ to $\bR$ that estimates the conditional means of the outcomes
\[
\mu_{\phi,\omega}(u):= \E [Y \given U_{\phi,\omega}=u,D=\omega].
\]
The general $\phi$-transformation based bias-corrected matching estimator is then defined to be
\begin{align*}
  \hat{\tau}_\phi := \frac{1}{n} \sum_{i=1}^n \big[\hat{Y}_{\phi,i}(1) -\hat{Y}_{\phi,i}(0)\big],
\end{align*}
where
\begin{align*}
    \hat{Y}_{\phi,i}(\omega) := \begin{cases}
    Y_i, & \mbox{ if } D_i=\omega,\\
      \frac{1}{M} \sum_{j \in \cJ_{\phi}(i)} (Y_j + \hat{\mu}_{\phi,\omega}(\hat{U}_{\phi,\omega,i}) - \hat{\mu}_{\phi,\omega}(\hat{U}_{\phi,\omega,j})) 
      & \mbox{ if } D_i=1-\omega
    \end{cases}
    .
\end{align*}
It follows that $\hat\tau_\phi$ generalizes $\hat\tau$ in \eqref{eq:tau}.

\subsection{General Theory}

In order to analyze $\hat\tau_\phi$, we introduce some additional notation and assumptions that are in parallel to those made in Section \ref{sec:theory}. Let the residuals from fitting the outcome models be
\[
\hat{R}_{\phi,i} := Y_i - \hat{\mu}_{\phi,D_i}(\hat{U}_{\phi,D_i,i}),\qquad i\in\zahl{n},
\]
and the estimator based on the outcome models be
\[
\hat{\tau}_{\phi}^{\rm reg}:= n^{-1} \sum_{i=1}^n \big[\hat{\mu}_{\phi,1}(\hat{U}_{\phi,1,i}) - \hat{\mu}_{\phi,0}(\hat{U}_{\phi,0,i})\big].
\]
Finally, let $K_{\phi}(i)$ be the number of matched times for the unit $i$ according to the distances between $\hat U_{\phi,D_i,i}$'s, i.e.,
\[
  K_{\phi}(i) := \sum_{j=1, D_j = 1-D_i}^n \ind\big(i \in \cJ_{\phi}(j)\big).
\]

The first lemma corresponds to a generalization of the AIPW representation of the bias-corrected rank-based estimator given in \eqref{eq:mbc} in Section \ref{sec:theory}.

\begin{lemma}\label{lemma:aipw} It holds true that
  \begin{align*}
    \hat{\tau}_\phi = \hat{\tau}_{\phi}^{\rm reg} + \frac{1}{n} \sum_{i=1}^n (2D_i-1)\Big(1 + \frac{K_{\phi}(i)}{M}\Big) \hat{R}_{\phi,i}.
  \end{align*}
\end{lemma}

The first two assumptions in this section parallel Assumptions \ref{asp:dr} and \ref{asp:dr-p}.

\begin{assumption} \label{asp:dr-g}
	\phantomsection $ $
  \begin{enumerate}[itemsep=-.5ex,label=(\roman*)]
    \item For almost all $x \in X$, $D$ is independent of $(Y(0),Y(1))$ conditional on $X=x$, and there exists some constant $c > 0$ such that $c < \P(D=1 \given X=x) < 1-c$.
    \item $\{(X_i,D_i,Y_i)\}_{i=1}^n$ are i.i.d. following the joint distribution of $(X,D,Y)$.
    \item $\E [(Y(\omega) - \mu_{\phi,\omega}(U_{\phi,\omega}))^2 \given U_{\phi,\omega} = u]$ is uniformly bounded for almost all $u \in \cX_\phi$ and $\omega =0,1$.
    \item $\E [\mu_{\phi,\omega}^2(U_{\phi,\omega})]$ is bounded for $\omega =0,1$.
  \end{enumerate}
\end{assumption}

\begin{assumption} \label{asp:dr-p-g}
	\phantomsection $ $
  \begin{enumerate}[itemsep=-.5ex,label=(\roman*)]
    \item There exists a deterministic, possibly changing with $n$, function $\bar{\mu}_{\phi,\omega}(\cdot):\bR^m \to \bR$  such that $\E [\bar{\mu}_{\phi,\omega}^2(U_{\phi,\omega})]$ is uniformly bounded and the estimator $\hat{\mu}_{\phi,\omega}(x)$ satisfies
    $
    \Vert \hat{\mu}_{\phi,\omega} - \bar{\mu}_{\phi,\omega} \Vert_\infty = o_\P(1)$ for $\omega =0,1$.
    \item $\max_{i \in \zahl{n}} \vert \bar{\mu}_{\phi,\omega}(\hat{U}_{\phi,\omega,i}) - \bar{\mu}_{\phi,\omega}(U_{\phi,\omega,i}) \vert = o_\P(1)$ for $\omega =0,1$.
  \end{enumerate}
\end{assumption}

The next assumption regulates the transformation $\phi_\omega$; cf. \citet[Section B]{lin2021estimation} and \citet[Assumption 3.3(ii)]{lin2022regression}. From a high level perspective, it roughly states that $M^{-1}K_{\phi}(i)$ should be a consistent density ratio estimator. A detailed discussion of this assumption is given in Section \ref{sec:key-assumption} ahead.

\begin{assumption} \label{asp:dr-p'-g}
	\phantomsection $ $
The number of matched times satisfies
\begin{align*}
\lim_{n \to \infty }\E \Big[ 
  \frac{K_{\phi}(1)}{M} - \Big(D_1 \frac{1-e(X_1)}{e(X_1)} + (1-D_1) \frac{e(X_1)}{1-e(X_1)} \Big)
\Big]^2 = 0,
\end{align*}
where, for any $x\in\cX$, $e(x):=\P(D=1\given X=x)$ is the propensity score.
\end{assumption}

The next three assumptions correspond to Assumptions \ref{asp:dr-o} through \ref{asp:se2} in Section \ref{sec:theory}.

\begin{assumption} \label{asp:dr-o-g}
	\phantomsection $ $
  \begin{enumerate}[itemsep=-.5ex,label=(\roman*)]
    \item The estimator $\hat{\mu}_{\phi,\omega}(x)$ satisfies
    $
    \Vert \hat{\mu}_{\phi,\omega} - \mu_{\phi,\omega} \Vert_\infty = o_\P(1)$ for $\omega =0,1$.
    \item $\max_{i \in \zahl{n}} \vert \mu_{\phi,\omega}(\hat{U}_{\phi,\omega,i}) - \mu_{\phi,\omega}(U_{\phi,\omega,i}) \vert = o_\P(1)$ for $\omega =0,1$.
    \item $\E[Y(\omega)\given X=x] = \mu_{\phi,\omega}(\phi_\omega(x))$ for almost all $x \in \cX$ and $\omega =0,1.$ 
  \end{enumerate}
\end{assumption}

\begin{assumption} \label{asp:se1-g} 
	\phantomsection $ $
\begin{enumerate}[itemsep=-.5ex,label=(\roman*)]
  \item $\E [(Y(\omega) - \mu_{\phi,\omega}(U_{\phi,\omega}))^2 \given U_{\phi,\omega} = u] $ is uniformly bounded away from zero for almost all $u \in \cX_\phi$ and $\omega =0,1$.
  \item There exists a constant $c>0$ such that $\E [\vert Y(\omega) - \mu_{\phi,\omega}(U_{\phi,\omega}) \vert ^{2+c} \given U_{\phi,\omega} = u]$ is uniformly bounded for almost all $u \in \cX_\phi$ and $\omega =0,1$.
  \item\label{asp:se1-3} $\max_{t \in \Lambda_{\max\{\lfloor m/2 \rfloor,1\} + 1}} \Vert \partial^t \mu_{\phi,\omega} \Vert_\infty$ is bounded. 
  \item $\E[Y(\omega)\given X=x] = \mu_{\phi,\omega}(\phi_\omega(x))$ for almost all $x \in \cX$ and $\omega = 0,1$. 
  \item The density of $\phi_\omega(X)$ is continuous over its support for $\omega =0,1$.
\end{enumerate}
\end{assumption}

\begin{assumption} \label{asp:se2-g} 
	\phantomsection $ $
	For $\omega =0,1$, the estimator $\hat{\mu}_\omega(x)$ satisfies
	\[
	  \max_{t \in \Lambda_{\max\{\lfloor m/2 \rfloor,1\} + 1}} \Vert \partial^t \hat{\mu}_{\phi,\omega} \Vert_\infty = O_\P(1)~~~{\rm and}~~~ 
	  \max_{t \in \Lambda_\ell} \Vert \partial^t \hat{\mu}_{\phi,\omega} - \partial^t \mu_{\phi,\omega} \Vert_\infty = O_\P(n^{-\gamma_\ell}) ~~\mbox{\rm for all}~~ \ell \in \zahl{\max\{\lfloor m/2 \rfloor,1\}},
	\]
	with some constants $\gamma_\ell > \max\{\frac{1}{2} - \frac{\ell}{m},0\} $ for $\ell=1,2,\ldots,\max\{\lfloor m/2 \rfloor,1\}$.
\end{assumption}

The next assumption poses a Donsker-type condition on the approximation accuracy of the estimated transformation $\hat\phi_\omega$ towards $\phi_\omega$. This assumption is usually needed when one wishes to avoid using sample splitting. 

\begin{assumption}  \phantomsection \label{asp:se3-g} 
For $\omega =0,1$, the estimator $\hat{\phi}_\omega$ satisfies
\[
  \lim_{\delta \to 0}\limsup_{n \to \infty}\P\Big(\sqrt{n} \sup_{x,y \in \cX, \Vert \phi_\omega(x)-\phi_\omega(y) \Vert \le \delta} \Vert (\hat{\phi}_\omega - \phi_\omega)(x) - (\hat{\phi}_\omega - \phi_\omega)(y)\Vert \ge \epsilon\Big) = 0.
\]
\end{assumption}

We are now ready to present the following theorem, which is a generalization to Theorem \ref{thm:bc}.

\begin{theorem}[Generalized Main Theorem] \label{thm:bcphi}
	\phantomsection $ $
\begin{enumerate}[itemsep=-.5ex,label=(\roman*)]
  \item\label{thm:drphi} (Double robustness of $\hat{\tau}_\phi$) If either Assumptions~\ref{asp:dr-g}, \ref{asp:dr-p-g}, \ref{asp:dr-p'-g} hold, or Assumptions~\ref{asp:dr-g} and \ref{asp:dr-o-g} hold, then 
  \[
  \hat{\tau}_\phi - \tau \stackrel{\sf p}{\longrightarrow} 0.
  \]
  \item\label{thm:sephi} (Semiparametric efficiency of $\hat{\tau}_\phi$) Assume the distribution of $(X,D,Y)$ satisfies Assumptions~\ref{asp:dr-g}, \ref{asp:dr-p'-g}, \ref{asp:se1-g}, \ref{asp:se2-g}, \ref{asp:se3-g}. Define
  \[
  \gamma  =  \max \Big\{ \Big[1- \frac{1}{2} \frac{m}{\max\{\lfloor m/2 \rfloor,1\} + 1}\Big] , \min_{\ell \in \zahl{\max\{\lfloor m/2 \rfloor,1\}} } \Big[1-\Big(\frac{1}{2} - \gamma_\ell \Big)\frac{m}{\ell}\Big] \Big\}, 
  \]
  recalling that $\gamma_\ell$'s were introduced in Assumption \ref{asp:se2-g}. Then, if $M \to \infty$ and $M/n^\gamma \to 0$ as $n \to \infty$, we have 
  \[
  \sqrt{n} (\hat{\tau}_\phi - \tau) \stackrel{\sf d}{\longrightarrow} N(0,\sigma^2).
  \]
\item If in addition Assumption~\ref{asp:dr-o-g} holds, then $\hat{\sigma}_\phi^2\stackrel{\sf p}{\longrightarrow}\sigma^2$, where
\[
\hat\sigma_\phi^2:=\frac{1}{n}\sum_{i=1}^n\Big[\hat\mu_{\phi,1}(\hat U_{\phi, 1,i})-\hat\mu_{\phi,0}(\hat U_{\phi,0,i})+(2D_i-1)\Big(1+\frac{K_\phi(i)}{M}\Big)\hat R_{\phi,i}-\hat\tau_\phi\Big]^2.
\]
\end{enumerate}
\end{theorem}

\subsection{Discussion on High-Level Assumption}\label{sec:key-assumption}


It remains to decipher the high-level condition in Assumption~\ref{asp:dr-p'-g}. To this end, we first give additional regularizations about the population-transformed data.

\begin{assumption} \label{asp:risk}
	\phantomsection $ $
  \begin{enumerate}[itemsep=-.5ex,label=(\roman*)]
    \item The diameter of $\cX_\phi$ and the surface area of $\cX_\phi$ are bounded.
    \item The density of $\phi_\omega(X)$ is continuous over its support for $\omega =0,1$.
    \item $\P(D=1\given \phi_\omega(X)=\phi_\omega(x)) = \P(D=1\given X=x)$ for almost all $x \in \cX$ and $\omega =0,1$.
  \end{enumerate}
\end{assumption}

Next, we give two different type of conditions for the estimator $\hat\phi_{\omega}$ to approximate $\phi_{\omega}$ so that Assumption~\ref{asp:dr-p'-g} can hold.

\begin{assumption} \label{asp:risk1}
	\phantomsection $ $
  For $\omega =0,1$,
  \[
    \lim_{n \to \infty} \E\Big[\Big(\frac{n}{M}\Big)^2 \sup_{x_1,x_2 \in \cX} \Vert \hat{\phi}_\omega(\cdot;x_1,x_2) - \phi_\omega\Vert_\infty^{2d}\Big] = 0,
  \]
  where $\hat{\phi}_\omega(\cdot;x_1,x_2)$ is the estimator constructed by inserting two more new points, $x_1$ and $x_2$, into the group with $D=1-\omega$ for some $x_1,x_2 \in \cX$.
\end{assumption}

\begin{assumption} \label{asp:risk2}
	\phantomsection $ $
  For $\omega =0,1$, we assume that for any fixed $\epsilon>0$, there exists a function $T_\epsilon(u)$ such that, for any $\delta>0$,
  \[
    \P\Big(\sup_{\delta \ge u}\delta^{-1}\sup_{\Vert \phi_\omega(s) - \phi_\omega(t)\Vert \le \delta}\sup_{x_1,x_2 \in \cX} \Vert (\hat{\phi}_\omega(\cdot;x_1,x_2) - \phi_\omega)(s) - (\hat{\phi}_\omega(\cdot;x_1,x_2) - \phi_\omega)(t) \Vert >  \epsilon \Big) \le T_\epsilon(u),
  \]
  and, for any $k \in \zahl{2}$,
  \[
    \lim_{n \to \infty} \Big(\frac{n}{M}\Big)^k \int_0^\infty u^{k-1}T_\epsilon(u^{1/m}) \d u =0
  \]
  and
  \begin{align*}
    \lim_{n \to \infty} \Big(\frac{n}{M}\Big)^2 \P\Big(\Vert \hat\phi_\omega-\phi_\omega \Vert_\infty > \epsilon\Big)= 0.
  \end{align*}
\end{assumption}

\begin{theorem}  \label{thm:suff}
	\phantomsection $ $
	Assume that Assumption~\ref{asp:risk} holds, $M\log n/n \to 0$ and $M \to \infty$ as $n \to \infty$. If either Assumption~\ref{asp:risk1} or Assumption~\ref{asp:risk2} holds, then Assumption~\ref{asp:dr-p'-g} holds.
\end{theorem}

\section{Conclusion}

This paper studied the large sample properties of Rosenbaum's rank-based matching estimator with regression adjustment, and established its consistency, doubly robutness, asymptotic normality, and semiparametric efficiency. Consistency of a plug-in variance estimator was also established. These results were obtained as a consequence of a more general theorem that allows for a class of transformations of the covariates, an leading special case being the empirical rank transformation proposed by \citet{rosenbaum2005exact,rosenbaum2010design}. To provide primitive conditions for regression adjustment, novel convergence rates for series estimators with generated regressors and possibly covariate Lebesgue density not bounded away from zero were derived, which may be of independent interest.

\section*{Acknowledgments}

We thank Peng Ding, Yingjie Feng, and Boris Shigida for insightful comments. 
Cattaneo gratefully acknowledge financial support from the National Science Foundation through grant SES-2241575. Han gratefully acknowledge financial support from the National Science Foundation through grants SES-2019363 and DMS-2210019.

\appendix

\section{Proofs of Main Results}\label{sec:main-proof}

\subsection{Proof of Theorem~\ref{thm:bc}}

We take $\phi_0 = \phi_1 = \mF$ and $\hat{\phi}_0 = \hat{\phi}_1 = \hat{\mF}_n$. Note that $\mF$ is a bijective function. Then Assumption~\ref{asp:dr} implies Assumption~\ref{asp:dr-g}. To show that Assumption~\ref{asp:dr-p} and Assumption~\ref{asp:dr-o} imply Assumption~\ref{asp:dr-p-g} and Assumption~\ref{asp:dr-o-g}, respectively, it remains to show 
\[
\max_{i \in \zahl{n}} \vert \bar{\mu}_{\omega}(\hat U_i) - \bar{\mu}_{\omega}(U_i) \vert = o_\P(1)~~ {\rm or}~~ \max_{i \in \zahl{n}} \vert \mu_{\omega}(\hat U_i) - \mu_{\omega}(U_i) \vert = o_\P(1). 
\]
Note that $[0,1]^d$ is compact. Then the continuity of $\mu_{\omega}$ implies uniform continuity. Then $\max_{i \in \zahl{n}} \vert \mu_{\omega}(\hat U_i) - \mu_{\omega}(U_i) \vert = o_\P(1)$ is directly from $\max_{i \in \zahl{n}} \Vert \hat U_i - U_i \Vert = o_\P(1)$. The same holds for $\bar{\mu}_{\omega}$.

To verify Assumption~\ref{asp:dr-p'-g}, we use Theorem~\ref{thm:suff}. Assumption~\ref{asp:risk} holds by Assumption~\ref{asp:dr-p'} and the bijection of $\mF$. We verify Assumption~\ref{asp:risk1} when $d \ge 2$ and Assumption~\ref{asp:risk2} when $d = 1$.

{\bf Part I.} $d \ge 2$.

Note that
\begin{align*}
  \sup_{x_1,x_2 \in \cX} \Vert \hat{\mF}_n(\cdot;x_1,x_2) - \mF\Vert_\infty^{2d} \lesssim \sup_{x_1,x_2 \in \cX} \Vert \hat{\mF}_n(\cdot;x_1,x_2) - \hat{\mF}_n \Vert_\infty^{2d} + \Vert \hat{\mF}_n - \mF\Vert_\infty^{2d},
\end{align*}
where $\lesssim$ means ``asymptotically small than''.

By the definition of $\hat{\mF}_n(\cdot;x_1,x_2)$ and $\hat{\mF}_n$, for any $x_1,x_2 \in \cX$,
\begin{align*}
  &\Vert \hat{\mF}_n(\cdot;x_1,x_2) - \hat{\mF}_n \Vert_\infty \\
  &\le \sqrt{d} \max_{k \in \zahl{d}} \max_{x \in \bR} \Big\vert \frac{1}{n+2}\Big( \sum_{i=1}^n \ind(X_{i,k} \le x) + \ind(x_1 \le x) + \ind(x_2 \le x)\Big) - \frac{1}{n}\sum_{i=1}^n \ind(X_{i,k} \le x) \Big\vert\\
  &= \sqrt{d} \max_{k \in \zahl{d}} \max_{x \in \bR} \Big\vert \frac{1}{n+2}\Big(\ind(x_1 \le x) + \ind(x_2 \le x)\Big) - \frac{2}{n(n+2)}\sum_{i=1}^n \ind(X_{i,k} \le x) \Big\vert\\
  &\le \frac{4\sqrt{d}}{n+2}.
\end{align*}
We then have
\begin{align*}
  \lim_{n \to \infty} \E\Big[\Big(\frac{n}{M}\Big)^2 \sup_{x_1,x_2 \in \cX} \Vert \hat{\mF}_n(\cdot;x_1,x_2) - \hat{\mF}_n \Vert_\infty^{2d}\Big] = 0.
\end{align*}
Note that
\begin{align*}
  \Vert \hat{\mF}_n - \mF\Vert_\infty^{2d}
  &\lesssim \max_{k \in \zahl{d}} \max_{x \in \bR} \Big\vert \frac{1}{n}\sum_{i=1}^n \ind(X_{i,k} \le x) -  \P(X_k \le x) \Big\vert^{2d} \\
  &\le \sum_{k=1}^d \max_{x \in \bR} \Big\vert \frac{1}{n}\sum_{i=1}^n \ind(X_{i,k} \le x) -  \P(X_k \le x) \Big\vert^{2d}.
\end{align*}
By the Dvoretzky–Kiefer–Wolfowitz inequality, we have
\begin{align*}
  \E[\Vert \hat{\mF}_n - \mF\Vert_\infty^{2d}] \lesssim n^{-d}.
\end{align*}
By $d \ge 2$ and $M \to \infty$, we have
\begin{align*}
  \lim_{n \to \infty} \E\Big[\Big(\frac{n}{M}\Big)^2 \Vert \hat{\mF}_n - \mF\Vert_\infty^{2d}\Big] = 0.
\end{align*}
The proof is now complete.

{\bf Part II.} $d = 1$.

Note that for any $s,t \in \bR$, by Part I,
\begin{align*}
  \sup_{x_1,x_2 \in \cX} \Big\vert (\hat{\mF}_n(\cdot;x_1,x_2) - \mF)(s) - (\hat{\mF}_n(\cdot;x_1,x_2) - \mF)(t) \Big\vert \le \frac{8\sqrt{d}}{n+2} + \Big\vert (\hat{\mF}_n - \mF)(s) - (\hat{\mF}_n - \mF)(t) \Big\vert.
\end{align*}
For any $\epsilon>0$, we can take $n$ sufficiently large such that $\frac{8\sqrt{d}}{n+2} < \epsilon$, and then
\begin{align*}
  &\P\Big(\sup_{\delta \ge u}\delta^{-1}\sup_{\vert \mF(s) - \mF(t)\vert \le \delta}\sup_{x_1,x_2 \in \cX} \Big\vert (\hat{\mF}_n(\cdot;x_1,x_2) - \mF)(s) - (\hat{\mF}_n(\cdot;x_1,x_2) - \mF)(t) \Big\vert >  2\epsilon \Big)\\
  &\le \P\Big(\sup_{\delta \ge u}\delta^{-1}\sup_{\vert \mF(s) - \mF(t)\vert \le \delta}\Big\vert (\hat{\mF}_n - \mF)(s) - (\hat{\mF}_n - \mF)(t) \Big\vert >  \epsilon \Big)\\
  &= \P\Big(\sup_{\delta \ge u}\sup_{s \le t:\P(s < X \le t) \le \delta} \delta^{-1}\Big\vert \frac{1}{n} \sum_{i=1}^n \ind(s < X_i \le t) - \P(s < X \le t) \Big\vert >  \epsilon \Big).
\end{align*}

Fix $u\ge0$ and $\epsilon$. Let $T_i = \delta^{-1} u (\ind(s < X_i \le t) - \P(s < X \le t))$ for $i \in \zahl{n}$. Then for any $i \in \zahl{n}$, we have $\E[T_i]=0$, and that $\vert T_i \vert \le 1$ for $\delta \ge u$.

Note that
\begin{align*}
  \sup_{\delta \ge u}\sup_{s \le t:\P(s < X \le t) \le \delta} \sum_{i=1}^n \E[T_i^2] \le \sup_{\delta \ge u}\sup_{s \le t:\P(s < X \le t) \le \delta} \sum_{i=1}^n \delta^{-2} u^2 \P(s < X_i \le t) \le u n,
\end{align*}
and
\begin{align*}
  \E[\sup_{\delta \ge u}\sup_{s \le t:\P(s < X \le t) \le \delta} \sum_{i=1}^n T_i^2] \le \E[\sup_{\delta \ge u}\sup_{s \le t:\P(s < X \le t) \le \delta} \sum_{i=1}^n \delta^{-2} u^2 \ind(s < X_i \le t)] \lesssim u n,
\end{align*}
by standard empirical process theory. Then by the concentration inequality for bounded processes \citep[Theorem 12.2]{boucheron2013concentration}, we have for $n$ sufficiently large,
\begin{align*}
  & \P\Big(\sup_{\delta \ge u}\sup_{s \le t:\P(s < X \le t) \le \delta} \delta^{-1}\Big\vert \frac{1}{n} \sum_{i=1}^n \ind(s < X_i \le t) - \P(s < X \le t) \Big\vert >  \epsilon \Big)\\
  &= \P\Big(\sup_{\delta \ge u}\sup_{s \le t:\P(s < X \le t) \le \delta} \Big\vert \sum_{i=1}^n T_i \Big\vert > u n  \epsilon \Big)\\
  &\le \exp\Big(-\frac{u^2n^2\epsilon^2}{Cun+u n  \epsilon}\Big) = \exp\Big(-\frac{\epsilon^2}{C+\epsilon} u n\Big).
\end{align*}

The proof for this part is now complete by taking integral using $m=1$ and $M \to \infty$.

Lastly, it is easy to see Assumption~\ref{asp:se1} and Assumption~\ref{asp:se2} imply Assumption~\ref{asp:se1-g} and Assumption~\ref{asp:se2-g}, respectively. Assumption~\ref{asp:se3-g} is followed by the Donsker's theorem applied to empirical distribution function.

\subsection{Proof of Theorem~\ref{thm:series}}

Let $\hat{\mQ} := \mP^\top \mP/n$. Using $\Vert \mQ^{-1/2} \Vert_2 = \lambda^{-1/2}_K$,
\begin{align*}
	\Vert \mQ^{-1/2} p_K(W_i)p_K(W_i)^\top \mQ^{-1/2}\Vert_2 = \Vert \mQ^{-1/2} p_K(W_i)\Vert^2 \leq \lambda^{-1}_K \zeta^2_{0,K},
\end{align*}
\begin{align*}
	&\Vert \E[\mQ^{-1/2} p_K(W_i)p_K(W_i)^\top \mQ^{-1} p_K(W_i)p_K(W_i)^\top\mQ^{-1/2}] \Vert_2^2\\
	&\leq \lambda^{-1}_K \zeta^2_{0,K} \Vert \E[\mQ^{-1/2} p_K(W_i)p_K(W_i)^\top\mQ^{-1/2}] \Vert_2^2 = \lambda^{-1}_K \zeta^2_{0,K},
\end{align*}
and a standard exponential concentration inequality for random matrices \citep[Section 6]{tropp2012user},
\begin{align*}
	\Vert \mQ^{-1/2}\hat{\mQ}\mQ^{-1/2} - \mI_K \Vert_2 = O_\P\big(\lambda^{-1/2}_K\zeta_{0,K} \sqrt{\log(K)/n} + \lambda^{-1}_K\zeta^2_{0,K} \log(K)/n \big) = o_\P(1)
\end{align*}
because $\lambda^{-1}_K\zeta^2_{0,K} \log(K)/n\to 0$ by assumption.

Let $\hat{\mQ}_n := \mP_n^\top \mP_n/n$. Then,
\begin{align*}
	\Vert \mQ^{-1/2}\mP^\top\Vert^2_2
	&= \Vert \mQ^{-1/2}\mP^\top\mP\mQ^{-1/2}\Vert_2\\
	&\leq n\Vert \mQ^{-1/2}(\hat{\mQ}-\mQ)\mQ^{-1/2}\Vert_2
	+ n\Vert \mQ^{-1/2}\mQ\mQ^{-1/2}\Vert_2\\
	&=O_\P(n\lambda^{-1/2}_K\zeta_{0,K} \sqrt{\log(K)/n} + n),
\end{align*}
and therefore
\begin{align*}
	&\Vert \mQ^{-1/2}(\hat{\mQ}_n - \hat{\mQ})\mQ^{-1/2} \Vert_2\\
	&= \Vert \mQ^{-1/2}(\mP_n^\top \mP_n -  \mP^\top \mP)\mQ^{-1/2} \Vert_2 /n\\
	&\le \Vert \mQ^{-1/2}(\mP_n-\mP)^\top(\mP_n-\mP)\mQ^{-1/2} \Vert_2 /n
	 + 2 \Vert \mQ^{-1/2}\mP^\top (\mP_n-\mP)\mQ^{-1/2} \Vert_2 /n\\
	&\le \Vert (\mP_n-\mP)\mQ^{-1/2} \Vert^2_2 /n
	 + 2 \Vert \mQ^{-1/2}\mP^\top\Vert_2 \Vert(\mP_n-\mP)\mQ^{-1/2} \Vert_2 /n\\
	&= O_\P(B_n + \lambda^{-1/4}_K\zeta_{0,K} (\log(K)/n)^{1/4} B^{1/2}_n + B^{1/2}_n)
	 =o_\P(1),
\end{align*}
because $\lambda^{-1}_K\zeta^2_{0,K} \log(K)/n\to 0$ and $B_n\to0$ by assumption.

Putting the two results together,
\begin{align*}
  \Vert \mQ^{-1/2} \hat{\mQ}_n \mQ^{-1/2} - \mI_K \Vert_2
  &\leq \Vert \mQ^{-1/2}(\hat{\mQ}_n - \hat{\mQ})\mQ^{-1/2} \Vert_2 + \Vert \mQ^{-1/2}\hat{\mQ}\mQ^{-1/2} - \mI_K \Vert_2\\
  &= O_\P( B^{1/2}_n + \lambda^{-1/2}_K\zeta_{0,K} \sqrt{\log(K)/n} ) = o_\P(1),
\end{align*}
given the rate restrictions imposed in the theorem.

Let $\ind_n = \ind(\lambda_{\min} (\mQ^{-1/2} \hat{\mQ}_n \mQ^{-1/2})>1/2)$. Then, $\lim_{n \to \infty}\P(\ind_n=1) = 1$. Letting $\varepsilon := \mY - \mPsi$,
\begin{align*}
	\ind_n\Vert \hat{\psi}_K - \psi_K \Vert^2_{L^2}
	&= \ind_n \int [p_K(w)^\top\hat{\beta}_{K} - p_K(w)^\top\beta_K]^2 dF_W(w)\\
	&= \ind_n (\hat{\beta}_{K} - \beta_K)^\top \E[p_K(W) p_K(W)^\top](\hat{\beta}_{K} - \beta_K)
	 = \ind_n \Vert \mQ^{1/2} (\hat{\beta}_{K} - \beta_K) \Vert^2\\
	&\leq 2\ind_n \Vert \mQ^{1/2}\hat{\mQ}_n^{-1}\mP_n^\top \varepsilon/n \Vert^2
	 + 2\ind_n \Vert \mQ^{1/2}\hat{\mQ}_n^{-1}\mP_n^\top (\mPsi - \mP_n\beta_K)/n \Vert^2.
\end{align*}
For the first term, we have
\begin{align*}
	\ind_n \Vert \mQ^{1/2}\hat{\mQ}_n^{-1}\mP_n^\top \varepsilon/n \Vert^2 = O_\P(K/n)
\end{align*}
because
\begin{align*}
	\ind_n \Vert \mQ^{1/2}\hat{\mQ}_n^{-1} \mP_n^\top \varepsilon/n \Vert^2
	\leq \ind_n \Vert \mQ^{1/2}\hat{\mQ}_n^{-1/2} \Vert^2_2 \Vert \hat{\mQ}_n^{-1/2}\mP_n^\top \varepsilon/n \Vert^2
	= O_\P(1) \ind_n \Vert\hat{\mQ}_n^{-1/2}\mP_n^\top \varepsilon/n \Vert^2,
\end{align*}
using $\ind_n \Vert \mQ^{1/2}\hat{\mQ}_n^{-1/2} \Vert^2_2=\ind_n \Vert \mQ^{1/2}\hat{\mQ}_n^{-1}\mQ^{1/2} \Vert_2=O_\P(1)$, and
\begin{align*}
\E[\Vert\hat{\mQ}_n^{-1/2}\mP_n^\top \varepsilon/n \Vert^2 | \cF_n] = {\rm tr}( \hat{\mQ}_n^{-1/2}\mP_n^\top\E(\varepsilon\varepsilon^\top | \cF_n) \mP_n\hat{\mQ}_n^{-1/2})/n^2=O_\P(K/n).
\end{align*}

We can bound the second term in different ways, depending on the approximation errors considered. The first two bounds rely on vanishing approximation errors ($\xi_{K}\to0$ or $\vartheta_{0,K}\to0$), and thus (implicitly) require $K\to\infty$ in general:
\begin{align*}
	&\ind_n \Vert \mQ^{1/2}\hat{\mQ}_n^{-1}\mP_n^\top (\mPsi - \mP_n\beta_K)/n \Vert^2\\
	&\leq \ind_n \Vert \mQ^{1/2}\hat{\mQ}_n^{-1}\mP_n^\top/\sqrt{n} \Vert^2_2 \Vert (\mPsi - \mP_n\beta_K)/\sqrt{n}\Vert^2 \\
	&\leq O_\P(1) \Vert \mPsi - \mP_n\beta_K \Vert^2/n\\
	&= O_\P(\min\{B_n + \xi^2_{K}  \,,\, R_n + \vartheta_{0, K}^2\}),
\end{align*}
because $\ind_n \Vert \mQ^{1/2}\hat{\mQ}_n^{-1}\mP_n^\top/\sqrt{n} \Vert^2_2 = \ind_n \Vert \mQ^{1/2}\hat{\mQ}_n^{-1}\mQ^{1/2} \Vert_2= O_\P(1)$, and because the term $\Vert \mPsi - \mP_n\beta_K \Vert^2/n$ can be bounded in two different ways:
\begin{align*}
	\Vert \mPsi - \mP_n\beta_K \Vert^2/n \leq 2\Vert \mPsi - \mPsi_n\Vert^2/n + 2\Vert \mPsi_n - \mP_n\beta_K \Vert^2/n = O_\P(R_n + \vartheta_{0, K}^2),
\end{align*}
or
\begin{align*}
	\Vert \mPsi - \mP_n\beta_K \Vert^2/n
	&\leq 2\Vert \mPsi - \mP\beta_K \Vert^2/n + 2\Vert (\mP-\mP_n)\beta_K\Vert^2/n\\
	&\leq O_\P(\xi^2_{K}) + 2\Vert (\mP-\mP_n)\mQ^{-1/2}\Vert^2_2 \Vert \E[\mQ^{-1/2}p_K(W_1)\psi(W_1)] \Vert^2/n\\
	&= O_\P(B_n + \xi^2_{K}),
\end{align*}
because $\Vert \E[\mQ^{-1/2}p_K(W_i)\psi(W_i)] = \beta_K^\top \mQ \beta_K = \E[(p_K(W_i)^\top \beta_K)^2] = \E[\psi_K(W_i)^2]\leq \E[\psi(W_i)^2]=O(1)$. Therefore, $\Vert (\mP-\mP_n)\beta_K\Vert^2/n = O_\P(B_n)$.

Next, for other possible bounds that do not require vanishing approximation errors ($\vartheta_{0,K}\geq \xi_{K} \not\to0$), even when possibly $K\to\infty$, notice that
\begin{align*}
	&\ind_n \Vert \mQ^{1/2}\hat{\mQ}_n^{-1}\mP_n^\top (\mPsi - \mP_n\beta_K)/n \Vert^2\\
	&\leq \ind_n \Vert \mQ^{1/2}\hat{\mQ}_n^{-1}\mQ^{1/2}\Vert^2_2 \Vert \mQ^{-1/2}\mP_n^\top (\mPsi - \mP_n\beta_K) \Vert^2 /n^2\\
	&\leq O_\P(1) \Vert \mQ^{-1/2}(\mP_n-\mP)^\top (\mPsi - \mP_n\beta_K) \Vert^2 /n^2 \\
	& \quad + O_\P(1) \Vert \mQ^{-1/2}\mP^\top (\mP_n-\mP)\beta_K \Vert^2 /n^2\\
	& \quad + O_\P(1) \Vert \mQ^{-1/2}\mP^\top (\mPsi - \mP\beta_K) \Vert^2 /n^2,
\end{align*}
where each of the three terms are bounded as follows. For the first term,
\begin{align*}
	&\Vert \mQ^{-1/2}(\mP_n-\mP)^\top (\mPsi - \mP_n\beta_K) \Vert^2 /n^2\\
	&\leq \Vert \mQ^{-1/2}(\mP_n-\mP)^\top\Vert^2_2 \Vert\mPsi - \mP_n\beta_K \Vert^2 /n^2 \\
	&= O_\P(B_n) O_\P(\min\{B_n + \xi^2_{K}  \,,\, R_n + \vartheta_{0, K}^2\}),
\end{align*}
using the calcultions above. For the second term,
\begin{align*}
	&\Vert \mQ^{-1/2}\mP^\top (\mP_n-\mP)\beta_K \Vert^2 /n^2\\
	&\leq \Vert \mQ^{-1/2}\mP^\top \Vert^2_2 \Vert (\mP_n-\mP)\beta_K \Vert^2 /n^2 \\
	&= O_\P(\lambda^{-1/2}_K\zeta_{0,K} \sqrt{\log(K)/n} + 1) O_\P(B_n),
\end{align*}
also using the calculations above. Finally, for the third and the last term,
\begin{align*}
	\Vert \mQ^{-1/2}\mP^\top (\mPsi - \mP\beta_K) \Vert^2 /n^2
	= O_\P(\min\{\lambda^{-1}_K \zeta_{0,K}^2 \xi^2_{K}, K \vartheta_{0, K}^2\} /n)
\end{align*}
because, by the orthogonality of the $L^2$ projection,
\begin{align*}
	&\E[\Vert \mQ^{-1/2}\mP^\top (\mPsi - \mP\beta_K) \Vert^2]/n^2\\
	&= \frac{1}{n^2} \E\Big[ \Big( \sum_{i=1}^n \mQ^{-1/2}p_{K}(W_i)(\psi(W_i) - \psi_K(W_i)) \Big)^\top \Big( \sum_{i=1}^n \mQ^{-1/2}p_{K}(W_i)(\psi(W_i) - \psi_K(W_i)) \Big)\Big] \\
	&= \frac{1}{n} \E\Big[ p_K(W_i)^\top \mQ^{-1}p_K(W_i)(\psi(W_i) - \psi_K(W_i))^2\Big].
\end{align*}

The final result in the theorem follows because $\lim_{n \to \infty}\P(\ind_n=1) = 1$.

\subsection{Proof of Lemma~\ref{lemma:eigen}}

For any $u \in \bR^K$ with $\Vert u \Vert =1$, by the orthonormality of the basis functions, we have
\begin{align*}
  1 = \Vert u \Vert^2 = \int \Big(\sum_{k=1}^K u_k p_{kK}(w) \Big)^2 \d w.
\end{align*}

Note that
\begin{align*}
  \lambda_{\min}(\E[p_K(W_1)p_K(W_1)^\top])
  &= \min_{u\in\bR^K:\Vert u \Vert=1} u^\top \E[p_K(W_1)p_K(W_1)^\top] u
   = \min_{u\in\bR^K:\Vert u \Vert=1} \E\Big[\Big(\sum_{k=1}^K u_k p_{kK}(W_1)\Big)^2\Big]\\
  &= \min_{u\in\bR^K:\Vert u \Vert=1} \int \Big(\sum_{k=1}^K u_k p_{kK}(w)\Big)^2 f_W(w) \d w,
\end{align*}
and since $\E[p_K(W_1)p_K(W_1)^\top]$ is positive semidefinite,
\begin{align*}
  \Vert \E[p_K(W_1)p_K(W_1)^\top]\Vert_2 = \max_{u\in\bR^K:\Vert u \Vert=1} \int \Big(\sum_{k=1}^K u_k p_{kK}(w)\Big)^2 f_W(w) \d w.
\end{align*}

If $f_W$ is bounded away from zero over the support of $W$, then for any $u \in \bR^K$ with $\Vert u \Vert =1$,
\begin{align*}
  \int \Big(\sum_{k=1}^K u_k p_{kK}(w)\Big)^2 f_W(w) \d w \ge c \int \Big(\sum_{k=1}^K u_k p_{kK}(w) \Big)^2 \d w = c,
\end{align*}
for some constants $c>0$. If $f_W$ is bounded over the support of $W$, then
\begin{align*}
  \int \Big(\sum_{k=1}^K u_k p_{kK}(w)\Big)^2 f_W(w) \d w \le C \int \Big(\sum_{k=1}^K u_k p_{kK}(w) \Big)^2 \d w = C,
\end{align*}
for some constants $C>0$.

More generally, for any $t>0$,
\begin{align*}
  &\int \Big(\sum_{k=1}^K u_k p_{kK}(w)\Big)^2 f_W(w) \d w \ge \int \Big(\sum_{k=1}^K u_k p_{kK}(w)\Big)^2 f_W(w) \ind(f_W(w) \ge t ) \d w\\
  &\ge t \int \Big(\sum_{k=1}^K u_k p_{kK}(w)\Big)^2 \ind(f_W(w) \ge t ) \d w\\
  &= t \Big[1-\int \Big(\sum_{k=1}^K u_k p_{kK}(w)\Big)^2 \ind(0 < f_W(w) < t ) \d w\Big].
\end{align*}

By the Cauchy-Schwarz inequality, for all sufficiently small $t>0$,
\begin{align*}
  &\int \Big(\sum_{k=1}^K u_k p_{kK}(w)\Big)^2 \ind(0 < f_W(w) < t ) \d w
   \le \int \Big(\sum_{k=1}^K p_{kK}^2(w)\Big) \ind(0 < f_W(w) < t ) \d w \\
  &\le \zeta^2_{0,K} \int \ind(0 < f_W(w) < t ) \d w \le C \zeta^2_{0,K} t^\rho.
\end{align*}

Take $t=c' \zeta^{-2/\rho}_{0,K}$ for some sufficiently small $c'>0$ such that $C \zeta^2_{0,K} t^\rho < 1/2$. We then obtain $\lambda_{\min}(\E[p_K(W_1)p_K(W_1)^\top]) \ge c' \zeta^{-2/\rho}_{0,K}/2$, as desired.

\subsection{Proof of Proposition~\ref{crl:series}}

For Part (i), $W$ follows the distribution of the Gaussian copula from the multivariate normal distribution with correlation matrix $\Sigma$, and thus the Lebesgue density of $W$ is
\begin{align*}
  f_W(w)= \frac{1}{\sqrt{\det \Sigma}} \exp\Big(-\frac{1}{2} (\Phi^{-1}(w_1),\ldots,\Phi^{-1}(w_d))(\Sigma^{-1}-\mI_d)(\Phi^{-1}(w_1),\ldots,\Phi^{-1}(w_d))^\top\Big),
\end{align*}
where $\Phi^{-1}(\cdot)$ is the inverse cumulative distribution function of a standard normal. Then,
\begin{align*}
  & \Big\{w:0<f_W(w)< t\Big\} \\
  &= \Big\{w:(\Phi^{-1}(w_1),\ldots,\Phi^{-1}(w_d))(\Sigma^{-1}-\mI_d)(\Phi^{-1}(w_1),\ldots,\Phi^{-1}(w_d))^\top > 2 \log \Big(\frac{1}{t\sqrt{\det \Sigma}}\Big)\Big\}\\
  &\subset \Big\{w:\Vert (\Phi^{-1}(w_1),\ldots,\Phi^{-1}(w_d)) \Vert^2 > \frac{2}{\lambda_{\max}(\Sigma^{-1}-\mI_d)} \log \Big(\frac{1}{t\sqrt{\det \Sigma}}\Big)\Big\}\\
  &\subset \bigcup_{k=1}^d \Big\{w:\Phi^{-1}(w_k)^2 > \frac{2}{d \lambda_{\max}(\Sigma^{-1}-\mI_d)} \log \Big(\frac{1}{t\sqrt{\det \Sigma}}\Big)\Big\}\\
  &= \bigcup_{k=1}^d \Big\{w:w_k > \Phi\Big(\sqrt{\frac{2}{d \lambda_{\max}(\Sigma^{-1}-\mI_d)} \log \Big(\frac{1}{t\sqrt{\det \Sigma}}\Big)}\Big)\Big\}.
\end{align*}

For any $k \in \zahl{d}$, by the Chernoff bound,
\begin{align*}
  & \mathsf{Leb}\Big(\Big\{w:w_k > \Phi\Big(\sqrt{\frac{2}{d \lambda_{\max}(\Sigma^{-1}-\mI_d)} \log \Big(\frac{1}{t\sqrt{\det \Sigma}}\Big)}\Big)\Big\}\Big)\\
  & = 1-\Phi\Big(\sqrt{\frac{2}{d \lambda_{\max}(\Sigma^{-1}-\mI_d)} \log \Big(\frac{1}{t\sqrt{\det \Sigma}}\Big)}\Big)\\
  &\le \exp\Big(-\frac{1}{d \lambda_{\max}(\Sigma^{-1}-\mI_d)} \log \Big(\frac{1}{t\sqrt{\det \Sigma}}\Big)\Big)\\
  &= \Big(t\sqrt{\det \Sigma}\Big)^{\frac{1}{d \lambda_{\max}(\Sigma^{-1}-\mI_d)}}.
\end{align*}

Then, we have
\begin{align*}
  \mathsf{Leb}(\{w:0<f_W(w)< t\}) \le d \Big(t\sqrt{\det \Sigma}\Big)^{\frac{1}{d \lambda_{\max}(\Sigma^{-1}-\mI_d)}},
\end{align*}
as desired.

\subsection{Proof of Lemma~\ref{lemma:aipw}}

By simple algebra, we have
\begin{align*}
  &\hat{\tau}_\phi = \frac{1}{n} \sum_{i=1}^n \Big[\hat{Y}_{\phi,i}(1) - \hat{Y}_{\phi,i}(0)\Big] \\
  &= \frac{1}{n} \sum_{i=1}^n D_i \Big[Y_i - \frac{1}{M} \sum_{j \in \cJ_{\phi}(i)} (Y_j + \hat{\mu}_{\phi,0}(\hat{U}_{\phi,0,i}) - \hat{\mu}_{\phi,0}(\hat{U}_{\phi,0,j})) \Big] \\
  & + \frac{1}{n} \sum_{i=1}^n (1-D_i)\Big[\frac{1}{M} \sum_{j \in \cJ_{\phi}(i)} (Y_j + \hat{\mu}_{\phi,1}(\hat{U}_{\phi,1,i}) - \hat{\mu}_{\phi,1}(\hat{U}_{\phi,1,j})) - Y_i \Big]\\
  &= \frac{1}{n} \sum_{i=1,D_i = 1}^n \Big[\hat{R}_{\phi,i} + \hat{\mu}_{\phi,1}(\hat{U}_{\phi,1,i}) - \hat{\mu}_{\phi,0}(\hat{U}_{\phi,0,i}) - \frac{1}{M} \sum_{j \in \cJ_{\phi}(i)} \hat{R}_{\phi,j} \Big] \\
  & + \frac{1}{n} \sum_{i=1,D_i = 0}^n \Big[\frac{1}{M} \sum_{j \in \cJ_{\phi}(i)} \hat{R}_{\phi,j} - \hat{R}_{\phi,i} + \hat{\mu}_{\phi,1}(\hat{U}_{\phi,1,i}) - \hat{\mu}_{\phi,0}(\hat{U}_{\phi,0,i}) \Big]\\
  &= \frac{1}{n} \sum_{i=1}^n \Big[\hat{\mu}_{\phi,1}(\hat{U}_{\phi,1,i}) - \hat{\mu}_{\phi,0}(\hat{U}_{\phi,0,i})\Big] + \frac{1}{n} \Big[ \sum_{i=1,D_i = 1}^n \Big(1 + \frac{K_{\phi}(i)}{M}\Big) \hat{R}_{\phi,i} - \sum_{i=1,D_i = 0}^n \Big(1 + \frac{K_{\phi}(i)}{M}\Big) \hat{R}_{\phi,i} \Big].
\end{align*}
This completes the proof.

\subsection{Proof of Theorem~\ref{thm:bcphi}\ref{thm:drphi}}

{\bf Part I.} Suppose the propensity score model is correct, i.e., Assumption~\ref{asp:dr-p-g} and \ref{asp:dr-p'-g} hold. For any $i \in \zahl{n}$, let $\bar{R}_{\phi,i} := Y_i - \bar{\mu}_{\phi,D_i}(U_{\phi,D_i,i})$. By Lemma~\ref{lemma:aipw},
\begin{align*}
    & \hat{\tau}_\phi = \hat{\tau}_{\phi}^{\rm reg} + \frac{1}{n} \sum_{i=1}^n (2D_i-1)\Big(1 + \frac{K_{\phi}(i)}{M}\Big) \hat{R}_{\phi,i}\\
    &= \frac{1}{n} \sum_{i=1}^n \Big[\hat{\mu}_{\phi,1}(\hat{U}_{\phi,1,i}) - \bar{\mu}_{\phi,1}(U_{\phi,1,i})\Big] - \frac{1}{n} \sum_{i=1}^n \Big[\hat{\mu}_{\phi,0}(\hat{U}_{\phi,0,i}) - \bar{\mu}_{\phi,0}(U_{\phi,0,i})\Big]\\
    &+  \frac{1}{n} \Big[ \sum_{i=1}^n D_i \Big(1 + \frac{K_{\phi}(i)}{M} \Big) \Big(\bar{\mu}_{\phi,1}(U_{\phi,1,i}) - \hat{\mu}_{\phi,1}(\hat{U}_{\phi,1,i}) \Big) - \sum_{i=1}^n (1-D_i)\Big(1 + \frac{K_{\phi}(i)}{M} \Big) \Big(\bar{\mu}_{\phi,0}(U_{\phi,0,i}) - \hat{\mu}_{\phi,0}(\hat{U}_{\phi,0,i}) \Big) \Big]\\
    &+  \frac{1}{n} \Big[ \sum_{i=1}^n D_i \Big(1 + \frac{K_{\phi}(i)}{M} - \frac{1}{e(X_i)}\Big) \bar{R}_{\phi,i} - \sum_{i=1}^n (1-D_i)\Big(1 + \frac{K_{\phi}(i)}{M} - \frac{1}{1-e(X_i)}\Big) \bar{R}_{\phi,i} \Big]\\
    & + \frac{1}{n} \Big[ \sum_{i=1}^n \Big(1 - \frac{D_i}{e(X_i)}\Big) \bar{\mu}_{\phi,1}(U_{\phi,1,i}) - \sum_{i=1}^n \Big(1 - \frac{1-D_i}{1-e(X_i)}\Big) \bar{\mu}_{\phi,0}(U_{\phi,0,i}) \Big]\\
    & + \frac{1}{n} \Big[ \sum_{i=1}^n \frac{D_i}{e(X_i)} Y_i - \sum_{i=1}^n \frac{1-D_i}{1-e(X_i)} Y_i \Big].
    \yestag\label{eq:dml1}
\end{align*}

For each pair of terms, we only establish the first half part under treatment, and the second half under control can be established in the same way.

For the first term in \eqref{eq:dml1}, by Assumption~\ref{asp:dr-p-g},
\begin{align*}
  &\Big\lvert \frac{1}{n} \sum_{i=1}^n \Big[\hat{\mu}_{\phi,1}(\hat{U}_{\phi,1,i}) - \bar{\mu}_{\phi,1}(U_{\phi,1,i})\Big] \Big\rvert \le \Big\lvert \frac{1}{n} \sum_{i=1}^n \Big[\hat{\mu}_{\phi,1}(\hat{U}_{\phi,1,i}) - \bar{\mu}_{\phi,1}(\hat{U}_{\phi,1,i})\Big] \Big\rvert + \Big\lvert \frac{1}{n} \sum_{i=1}^n \Big[\bar{\mu}_{\phi,1}(\hat{U}_{\phi,1,i}) - \bar{\mu}_{\phi,1}(U_{\phi,1,i})\Big] \Big\rvert \\
  &\le \lVert \hat{\mu}_1 - \bar{\mu}_1 \rVert_\infty + \max_{i \in \zahl{n}} \Big\lvert \bar{\mu}_{\phi,1}(\hat{U}_{\phi,1,i}) - \bar{\mu}_{\phi,1}(U_{\phi,1,i}) \Big\rvert   = o_\P(1).
\end{align*}
Then
\begin{align}\label{eq:dml11}
  \frac{1}{n} \sum_{i=1}^n \Big[\hat{\mu}_{\phi,1}(\hat{U}_{\phi,1,i}) - \bar{\mu}_{\phi,1}(U_{\phi,1,i})\Big] = o_\P(1).
\end{align}

For the second term in \eqref{eq:dml1}, by Assumption~\ref{asp:dr-p-g},
\begin{align*}
    & \Big\lvert \frac{1}{n} \sum_{i=1}^n D_i \Big(1 + \frac{K_{\phi}(i)}{M} \Big) \Big(\bar{\mu}_{\phi,1}(U_{\phi,1,i}) - \hat{\mu}_{\phi,1}(\hat{U}_{\phi,1,i}) \Big) \Big\rvert \\
    &\le \max_{i \in \zahl{n}} \Big\lvert \bar{\mu}_{\phi,1}(U_{\phi,1,i}) - \hat{\mu}_{\phi,1}(\hat{U}_{\phi,1,i}) \Big\rvert \cdot \frac{1}{n} \sum_{i=1}^n D_i \Big(1 + \frac{K_{\phi}(i)}{M} \Big) = \max_{i \in \zahl{n}} \Big\lvert \bar{\mu}_{\phi,1}(U_{\phi,1,i}) - \hat{\mu}_{\phi,1}(\hat{U}_{\phi,1,i}) \Big\rvert\\
    &\le \lVert \hat{\mu}_1 - \bar{\mu}_1 \rVert_\infty + \max_{i \in \zahl{n}} \Big\lvert \bar{\mu}_{\phi,1}(\hat{U}_{\phi,1,i}) - \bar{\mu}_{\phi,1}(U_{\phi,1,i}) \Big\rvert = o_\P(1).
\end{align*}
We then have
\begin{align}\label{eq:dml12}
  \frac{1}{n} \sum_{i=1}^n D_i \Big(1 + \frac{K_{\phi}(i)}{M} \Big) \Big(\bar{\mu}_{\phi,1}(U_{\phi,1,i}) - \hat{\mu}_{\phi,1}(\hat{U}_{\phi,1,i}) \Big) = o_\P(1).
\end{align}

For the third term in \eqref{eq:dml1}, by the Cauchy-Schwarz inequality,
\begin{align*}
    & \Big\lvert\frac{1}{n} \sum_{i=1}^n D_i \Big(1 + \frac{K_{\phi}(i)}{M} - \frac{1}{e(X_i)}\Big) \bar{R}_{\phi,i} \Big\rvert\\
    &\le \Big\{\frac{1}{n} \sum_{i=1}^n D_i \Big(1 + \frac{K_{\phi}(i)}{M} - \frac{1}{e(X_i)}\Big)^2 \Big\}^{1/2} \Big\{\frac{1}{n} \sum_{i=1}^n D_i \bar{R}_{\phi,i}^2\Big\}^{1/2}.
\end{align*}
Note that by Assumptions~\ref{asp:dr-g} and \ref{asp:dr-p-g},
\begin{align*}
  & \E \Big[\frac{1}{n} \sum_{i=1}^n D_i \bar{R}_{\phi,i}^2\Big] = \E\Big[D_1 \bar{R}_{\phi,1}^2\Big] = \E\Big[D_1 \Big(Y_1(1) - \bar{\mu}_{\phi,1}(U_{\phi,1,1})\Big)^2\Big]\\
  &\le 2\E\Big[D_1 \Big(\sigma_1^2(U_{\phi,1,1})+[\mu_{\phi,1}(U_{\phi,1,1})-\bar{\mu}_{\phi,1}(U_{\phi,1,1})]^2\Big)\Big] < \infty,
\end{align*}
where $\sigma_1^2(u) := \E [(Y(1)-\mu_{\phi,1}(u))^2 \given U_{\phi,1}=u]$ for $u \in \cX_\phi$. We then obtain by Assumption~\ref{asp:dr-p'-g} and the Markov inequality that
\begin{align}\label{eq:dml13}
  \frac{1}{n} \sum_{i=1}^n D_i \Big(1 + \frac{K_{\phi}(i)}{M} - \frac{1}{e(X_i)}\Big) \bar{R}_{\phi,i} = o_\P(1).
\end{align}

For the fourth term in \eqref{eq:dml1}, notice that $\bar{\mu}_{\phi,1}(U_{\phi,1,i})$ is a function of $X_i$. Then by the definition of the propensity score and Assumption~\ref{asp:dr-g},
\[
  \E\Big[\Big(1 - \frac{D_i}{e(X_i)}\Big) \bar{\mu}_{\phi,1}(U_{\phi,1,i})\Big] = 0, ~~ \E\Big[\Big\lvert\Big(1 - \frac{D_i}{e(X_i)}\Big) \bar{\mu}_{\phi,1}(U_{\phi,1,i})\Big\rvert\Big] < \infty.
\]
By the i.i.d of $[(X_i,D_i)]_{i=1}^n$ and the weak law of large numbers, we have
\begin{align}\label{eq:dml15}
  \frac{1}{n} \sum_{i=1}^n \Big(1 - \frac{D_i}{e(X_i)}\Big) \bar{\mu}_{\phi,1}(U_{\phi,1,i}) = o_\P(1).
\end{align}

For the fifth term in \eqref{eq:dml1}, notice that $\E [\lvert Y \rvert]$ is bounded from Assumption~\ref{asp:dr} and $[(X_i,D_i,Y_i)]_{i=1}^n$ are i.i.d.. Using the weak law of large numbers yields
\begin{align}\label{eq:dml16}
    \frac{1}{n} \Big[ \sum_{i=1}^n \frac{D_i}{e(X_i)} Y_i - \sum_{i=1}^n \frac{1-D_i}{1-e(X_i)} Y_i \Big] \stackrel{\sf p}{\longrightarrow} \E\Big[Y_i(1)-Y_i(0)\Big] = \tau.
\end{align}

Plugging \eqref{eq:dml11}, \eqref{eq:dml12}, \eqref{eq:dml13}, \eqref{eq:dml15} into \eqref{eq:dml1} completes the proof.

{\bf Part II.} Suppose the outcome model is correct, i.e., Assumption~\ref{asp:dr-o-g} holds. Using the representation \eqref{eq:mbc},
\begin{align*}
  & \hat{\tau}_\phi = \hat{\tau}_{\phi}^{\rm reg} + \frac{1}{n} \sum_{i=1}^n (2D_i-1)\Big(1 + \frac{K_{\phi}(i)}{M}\Big) \hat{R}_{\phi,i}\\
  &= \frac{1}{n} \sum_{i=1}^n \Big[\hat{\mu}_{\phi,1}(\hat{U}_{\phi,1,i}) - \mu_{\phi,1}(U_{\phi,1,i})\Big] - \frac{1}{n} \sum_{i=1}^n \Big[\hat{\mu}_{\phi,0}(\hat{U}_{\phi,0,i}) - \mu_{\phi,0}(U_{\phi,0,i})\Big]\\
    &+  \frac{1}{n} \Big[ \sum_{i=1}^n D_i \Big(1 + \frac{K_{\phi}(i)}{M} \Big) \Big(\mu_{\phi,1}(U_{\phi,1,i}) - \hat{\mu}_{\phi,1}(\hat{U}_{\phi,1,i}) \Big) - \sum_{i=1}^n (1-D_i)\Big(1 + \frac{K_{\phi}(i)}{M} \Big) \Big(\mu_{\phi,0}(U_{\phi,0,i}) - \hat{\mu}_{\phi,0}(\hat{U}_{\phi,0,i}) \Big) \Big]\\
    &+  \frac{1}{n} \Big[ \sum_{i=1}^n D_i \Big(1 + \frac{K_{\phi}(i)}{M} \Big) \Big(Y_i - \mu_{\phi,1}(U_{\phi,1,i}) \Big) - \sum_{i=1}^n (1-D_i) \Big(1 + \frac{K_{\phi}(i)}{M} \Big) \Big(Y_i - \mu_{\phi,0}(U_{\phi,0,i}) \Big) \Big]\\
    &+  \frac{1}{n} \sum_{i=1}^n \Big[\mu_{\phi,1}(U_{\phi,1,i}) - \mu_{\phi,0}(U_{\phi,0,i})\Big].
    \yestag\label{eq:dml2}
\end{align*}

For the first term in \eqref{eq:dml2}, in the same way as \eqref{eq:dml11},
\begin{align}\label{eq:dml21}
  \frac{1}{n} \sum_{i=1}^n \Big[\hat{\mu}_{\phi,1}(\hat{U}_{\phi,1,i}) - \mu_{\phi,1}(U_{\phi,1,i})\Big] = o_\P(1).
\end{align}

For the second term in \eqref{eq:dml2}, in the same way as \eqref{eq:dml12},
\begin{align}\label{eq:dml22}
  \frac{1}{n} \sum_{i=1}^n D_i \Big(1 + \frac{K_{\phi}(i)}{M} \Big) \Big(\mu_{\phi,1}(U_{\phi,1,i}) - \hat{\mu}_{\phi,1}(\hat{U}_{\phi,1,i}) \Big) = o_\P(1).
\end{align}

For the third term in \eqref{eq:dml2}, noticing that $[K_{\phi}(i)]_{i=1}^n$ is a function of $\{(X_i,D_i)\}_{i=1}^n$, by Assumption~\ref{asp:dr-g} and Assumption~\ref{asp:dr-o-g}, we can obtain for any $i \in \zahl{n}$,
\begin{align*}
  &\E \Big[ D_i \Big(1 + \frac{K_{\phi}(i)}{M} \Big) \Big(Y_i - \mu_{\phi,1}(U_{\phi,1,i}) \Big) \Biggiven \{(X_i,D_i)\}_{i=1}^n \Big] \\
  &= D_i \Big(1 + \frac{K_{\phi}(i)}{M} \Big) \Big(\E[Y_i\given X_i,D_i=1] - \mu_{\phi,1}(U_{\phi,1,i}) \Big)\\
  &= D_i \Big(1 + \frac{K_{\phi}(i)}{M} \Big) \Big(\E[Y_i(1)\given X_i] - \mu_{\phi,1}(U_{\phi,1,i}) \Big) = 0,
\end{align*}
and
\begin{align*}
    &\E \Big[ \Big\lvert \frac{1}{n} \sum_{i=1}^n D_i \Big(1 + \frac{K_{\phi}(i)}{M} \Big) \Big(Y_i - \mu_{\phi,1}(U_{\phi,1,i}) \Big) \Big\rvert\Big] \le \E \Big[ \Big\lvert \frac{1}{n} \sum_{i=1}^n D_i \Big(1 + \frac{K_{\phi}(i)}{M} \Big) \Big\rvert\Big] \lVert \sigma_1 \rVert_\infty\\
    &\lesssim \lVert \sigma_1 \rVert_\infty = O(1).
\end{align*}
Accordingly, by the martingale convergence theorem in the same way as \cite{abadie2012martingale}, we obtain
\begin{align}\label{eq:dml23}
  \frac{1}{n} \sum_{i=1}^n D_i \Big(1 + \frac{K_{\phi}(i)}{M} \Big) \Big(Y_i - \mu_{\phi,1}(U_{\phi,1,i}) \Big) = o_\P(1).
\end{align}

For the fourth term in \eqref{eq:dml2}, notice that $\E [\mu_{\phi,\omega}^2(U_{\phi,\omega})]$ is bounded for $\omega =0,1$. Using the weak law of large number, we obtain
\begin{align}\label{eq:dml24}
  \frac{1}{n} \sum_{i=1}^n \Big[\mu_{\phi,1}(U_{\phi,1,i}) - \mu_{\phi,0}(U_{\phi,0,i})\Big] \stackrel{\sf p}{\longrightarrow} \E \Big[\mu_1(X_1) - \mu_0(X_1)\Big] = \tau.
\end{align}

Plugging \eqref{eq:dml21}, \eqref{eq:dml22}, \eqref{eq:dml23}, \eqref{eq:dml24} into \eqref{eq:dml2} completes the proof.

\subsection{Proof of Theorem~\ref{thm:bcphi}\ref{thm:sephi}}

We decompose $\hat{\tau}_\phi$ as
\begin{align*}
  & \hat{\tau}_\phi = \hat{\tau}_{\phi}^{\rm reg} + \frac{1}{n} \sum_{i=1}^n (2D_i-1)\Big(1 + \frac{K_{\phi}(i)}{M}\Big) \hat{R}_{\phi,i}\\
  &= \frac{1}{n} \sum_{i=1}^n \Big[\mu_{\phi,1}(U_{\phi,1,i}) - \mu_{\phi,0}(U_{\phi,0,i})\Big] + \frac{1}{n} \sum_{i=1}^n (2D_i-1)\Big(1 + \frac{K_{\phi}(i)}{M}\Big) \Big(Y_i - \mu_{\phi,D_i}(U_{\phi,D_i,i})\Big)\\
  & + \frac{1}{n}\sum_{i=1}^n (2D_i-1) \Big[\mu_{\phi,1-D_i}(U_{\phi,1-D_i,i}) - \frac{1}{M} \sum_{j \in \cJ_{\phi}(i)} \mu_{\phi,1-D_i}(U_{\phi,1-D_i,j}) \Big]\\
  & - \frac{1}{n}\sum_{i=1}^n (2D_i-1) \Big[\hat{\mu}_{\phi,1-D_i}(\hat{U}_{\phi,1-D_i,i}) - \frac{1}{M} \sum_{j \in \cJ_{\phi}(i)} \hat{\mu}_{\phi,1-D_i}(\hat{U}_{\phi,1-D_i,j}) \Big]\\
  =:& \bar{\tau}_{\phi} + E_n + B_n - \hat{B}_n.
\end{align*}

In the same way as Lemmas A.1 and A.2 in \cite{lin2022regression}, we have the following central limit theorem on $\bar{\tau}_{\phi} + E_n$.

\begin{lemma}\label{lemma:mbc,clt}
    Under Assumptions \ref{asp:dr-g}, \ref{asp:dr-p'-g}, \ref{asp:se1-g},
    \begin{align*}
        \sqrt{n} \sigma^{-1} \Big(\bar{\tau}_{\phi} + E_n - \tau \Big)\stackrel{\sf d}{\longrightarrow} N\Big(0,1\Big).
    \end{align*}
\end{lemma}

For the bias term $B_M-\hat B_M$, in light of the smoothness conditions on $\mu_\omega$ and approximation conditions on $\hat{\mu}_\omega$ for $\omega =0,1$, one can establish the following lemma.

\begin{lemma}\label{lemma:mbc,bias}
    Under Assumptions \ref{asp:dr-g}, \ref{asp:dr-p'-g}, \ref{asp:se1-g}, \ref{asp:se2-g}, \ref{asp:se3-g},
    \begin{align*}
        \sqrt{n} \Big(B_n - \hat{B}_n \Big)\stackrel{\sf p}{\longrightarrow} 0.
    \end{align*}
\end{lemma}

Combining Lemma~\ref{lemma:mbc,clt} and Lemma~\ref{lemma:mbc,bias} completes the proof.

The consistency of the variance estimator can be established in a similar way as the proof of Theorem 4.1 in \cite{lin2021estimation}.

\subsection{Proof of Theorem~\ref{thm:suff}}

Similar to \cite{lin2021estimation}, we first consider a two-sample density ratio estimation problem.

With an abuse of notation, restricted to this section let's consider two general random vectors $X,Z$ in $\bR^d$ that are defined on the same probability space. Let $\cX, \cZ \subset \bR^d$ be the supports of $X$ and $Z$ respectively with $\cZ \subset \cX$. Consider a general function $\phi:\cX \to \bR^m$. Let $\nu_0$ and $\nu_1$ represent the probability measures of $\phi(X)$ and $\phi(Z)$, respectively. Assume $\nu_0$ and $\nu_1$ are absolutely continuous with respect to the Lebesgue measure $\lambda$ on $\bR^m$ equipped with the Euclidean norm $\lVert \cdot \rVert$; denote the corresponding densities (Radon-Nikodym derivatives) by $f_0$ and $f_1$. Assume further that $\nu_1$ is absolutely continuous with respect to $\nu_0$ and write the corresponding density ratio, $f_1/f_0$, as $r$; set $0/0=0$ by default.

Assume $X_1,\ldots,X_{N_0}$ are $N_0$ independent copies of $X$, $Z_1,\ldots,Z_{N_1}$ are $N_1$ independent copies of $Z$, and $[X_i]_{i=1}^{N_0}$ and $[Z_j]_{j=1}^{N_1}$ are mutually independent. We aim to estimate the density ratio $r(\phi(x))$ for any $x \in \cX$ based on $\{X_1,\ldots,X_{N_0},Z_1,\ldots,Z_{N_1}\}$.

For any $x \in \cX$, we consider a general estimator $\hat\phi$ estimating $\phi$, which may depend on 
\[
\{X_1,\ldots,X_{N_0},Z_1,\ldots,Z_{N_1}\} ~{\rm and}~ x. 
\]
Define the catchment area of $x$:
\begin{align}
  A_\phi(x) = A_\phi\big(x,\{X_i\}_{i=1}^{N_0},\hat\phi\big) := \Big\{z \in \cZ: \lVert \hat\phi(x) - \hat\phi(z) \rVert \le \hat\Phi_M(z)\Big\},
\end{align}
where $\hat\Phi_M(z)$ is the $M$-th order statistics of $\{\lVert \hat\phi(X_i) - \hat\phi(z) \rVert\}_{i=1}^{N_0}$, and the number of matched times of $x$:
\begin{align}
  K_\phi(x) = K_\phi\big(x,\{X_i\}_{i=1}^{N_0},\{Z_j\}_{j=1}^{N_1}\big) := \sum_{j=1}^{N_1} \ind\Big(Z_j \in A_\phi(x)\Big).
\end{align}

Then the density ratio estimator is defined as:
\begin{align}
  \hat{r}_\phi(x) = \hat{r}_\phi\big(x,\{X_i\}_{i=1}^{N_0},\{Z_j\}_{j=1}^{N_1}\big):= \frac{N_0}{N_1} \frac{K_\phi(x)}{M}.
\end{align}

For any positive integer $p$, let $\hat\phi_{(Z_1,\ldots,Z_p) \to z}$ be the estimator replacing $(Z_1,\ldots,Z_p)$ by $z$ for $z \in \cZ^p$. We consider the following two assumptions, which are analogies of Assumption~\ref{asp:risk1} and Assumption \ref{asp:risk2} in the two-sample problem.

\begin{assumption} \phantomsection \label{asp:proof1}
  \[
    \lim_{N_0 \to \infty} \E\Big[\Big(\frac{N_0}{M}\Big)^p \sup_{z \in \cZ^p} \lVert \hat\phi_{(Z_1,\ldots,Z_p) \to z}-\phi \rVert_\infty^{pd}\Big] = 0.
  \] 
\end{assumption}

\begin{assumption} \phantomsection \label{asp:proof2}
 For any $\epsilon>0$ and $\delta>0$,
  \[
    \P\Big(\sup_{\delta \in \bR: \delta \ge u}\delta^{-1}\sup_{s,t \in \cX:\lVert \phi(s) - \phi(t)\rVert \le \delta}\sup_{z \in \cZ^p}\lVert (\hat\phi_{(Z_1,\ldots,Z_p) \to z}-\phi)(s) - (\hat\phi_{(Z_1,\ldots,Z_p) \to z}-\phi)(t) \rVert >  \epsilon \Big) \le T_\epsilon(u),
  \]
  for $T_\epsilon(u)$ satisfying for any $k \in \zahl{p}$,
  \[
    \lim_{N_0 \to \infty} \Big(\frac{N_0}{M}\Big)^k \int_0^\infty u^{k-1} T_\epsilon(u^{1/m}) \d u =0,
  \]
  and
  \begin{align*}
    \lim_{N_0 \to \infty} \Big(\frac{N_0}{M}\Big)^p \P\Big(\lVert \hat\phi-\phi \rVert_\infty > \epsilon\Big)= 0.
  \end{align*}
\end{assumption}

The following theorem considers the asymptotic $L^p$ moments of $\hat{r}_\phi$.

\begin{theorem}[Asymptotic $L^p$ moments of $\hat{r}_\phi$]  \label{lemma:moment} Let $p$ be any positive integer. Assume Assumption~\ref{asp:proof1} or \ref{asp:proof2} holds for $p$. Assume $M\log N_0/N_0 \to 0$, $MN_1/N_0 \to \infty$ and $M \to \infty$ as $N_0 \to \infty$. We then have
  \[
    \lim_{N_0\to\infty} \E\big[(\hat{r}_\phi(x))^p\big] = \big[r(\phi(x))\big]^p
  \]
  holds for all $x \in \cX$ such that $f_0(\phi(x))>0$ and $f_0,f_1$ are continuous at $\phi(x)$.
\end{theorem}

The proof of Theorem~\ref{lemma:moment} will use the following lemma.

\begin{lemma}\label{lemma:moment,p} Under the same conditions of Theorem~\ref{lemma:moment}, we have
  \begin{align*}
    \lim_{N_0 \to \infty} \Big(\frac{N_0}{M}\Big)^p \P\Big(Z_1,\ldots,Z_p \in A_\phi(x)\Big) = \big[r(\phi(x))\big]^p.
  \end{align*}
  holds for all $x \in \cX$ such that $f_0(\phi(x))>0$ and $f_0,f_1$ are continuous at $\phi(x)$.
\end{lemma}

As a direct result of Theorem~\ref{lemma:moment}, we can establish the pointwise consistency of the estimator $\hat{r}_\phi$.

\begin{corollary}[Pointwise consistency]\label{thm:cons,lp} Under the same conditions as Theorem~\ref{lemma:moment}, if $p$ is even, we have
    \[
      \lim_{N_0\to\infty} \E \big[ \lvert \hat{r}_\phi(x) - r(\phi(x))\rvert^p\big] = 0
    \]
    holds for all $x \in \cX$ such that $f_0(\phi(x))>0$ and $f_0,f_1$ are continuous at $\phi(x)$.
\end{corollary}

The pointwise consistency of $\hat{r}_\phi$ can then be generalized to global consistency under the following assumptions on $\cX$.

\begin{assumption} \phantomsection \label{asp:proof3}
  \begin{enumerate}[itemsep=-.5ex,label=(\roman*)]
    \item $\cX$ is compact and the surface areas of $\cX$ and $\cZ$ are bounded.
    \item $r$ is bounded over $\cX$.
    \item $f_0$ is continuous over $\cX$ and $f_1$ is continuous over $\cZ$.
  \end{enumerate}
\end{assumption}

\begin{theorem}[Global consistency] \label{thm:cons} Under the same conditions of Theorem~\ref{lemma:moment} and Assumption~\ref{asp:proof3}, if $p$ is even, we have
  \[
    \lim_{N_0\to\infty} \E \Big[ \Big\lvert \hat{r}_\phi(X) - r(\phi(X))\Big\rvert^p\Big] = 0.
  \]
\end{theorem}

Now back to the causal setting, the density ratio $r(\phi(x))$ is not necessarily equal to the density ratio we are interested in. We consider a lemma showing the equivalence of the two density ratios under additional assumptions.
\begin{lemma}\label{lemma:drequal} 
Let $f_{X \given D=1}$ and $f_{X \given D=0}$ be the density of $X \given D=1$ and $X \given D=0$, respectively. Let $f_{\phi,X \given D=1}$ and $f_{\phi,X \given D=0}$ be the density of $\phi(X) \given D=1$ and $\phi(X) \given D=0$. Then for any $x \in \cX$ such that $\P(D=1\given \phi(X)=\phi(x)) = \P(D=1\given X=x)$, we have
  \[
    \frac{f_{\phi,X \given D=1}(\phi(x))}{f_{\phi,X \given D=0}(\phi(x))} = \frac{f_{X \given D=1}(x)}{f_{X \given D=0}(x)}.
  \]
\end{lemma}

Note that
\begin{align*}
  &\E \Big[\frac{K_{\phi}(1)}{M} - \Big(D_1 \frac{1-e(X_1)}{e(X_1)} + (1-D_1) \frac{e(X_1)}{1-e(X_1)} \Big)\Big]^2\\
  &= \E\Big[\E \Big[\Big[\frac{K_{\phi}(1)}{M} - \Big(D_1 \frac{1-e(X_1)}{e(X_1)} + (1-D_1) \frac{e(X_1)}{1-e(X_1)} \Big)\Big]^2 \Biggiven \{D_i\}_{i=1}^n \Big]\Big]\\
  &= \E\Big[\E \Big[\Big(\frac{K_{\phi}(1)}{M} - \frac{1-e(X_1)}{e(X_1)} \Big)^2 \Biggiven \{D_i\}_{i=1}^n, D_1=1 \Big] \ind\Big(D_1=1\Big)\Big] \\
  &+ \E\Big[\E \Big[\Big(\frac{K_{\phi}(1)}{M} - \frac{e(X_1)}{1-e(X_1)} \Big)^2 \Biggiven \{D_i\}_{i=1}^n, D_1=0 \Big] \ind\Big(D_1=0\Big)\Big].
\end{align*}

We consider the second term for example. Conditional on $\{D_i\}_{i=1}^n$, $[X_i]_{i:D_i=0}$ and $[X_i]_{i:D_i=1}$ are two samples from $X \given D=0$ and $X \given D=1$, respectively. Note that
\begin{align*}
  \E \Big[\Big(\frac{K_{\phi}(1)}{M} - \frac{e(X_1)}{1-e(X_1)} \Big)^2 \Biggiven \{D_i\}_{i=1}^n, D_1=0 \Big] = \Big(\frac{N_1}{N_0}\Big)^2\E \Big[\Big(\frac{N_0}{N_1}\frac{K_{\phi}(1)}{M} - \frac{N_0}{N_1} \frac{e(X_1)}{1-e(X_1)} \Big)^2 \Biggiven \{D_i\}_{i=1}^n, D_1=0 \Big].
\end{align*}

By the strong law of large number, we have $(N_0/N_1) \to \P(D=0)/\P(D=1)$ with probability one. Note that
\begin{align*}
  \frac{\P(D=0)}{\P(D=1)} \frac{e(X_1)}{1-e(X_1)} = \frac{f_{X \given D=1}(x)}{f_{X \given D=0}(x)}.
\end{align*}

To apply Theorem~\ref{thm:cons} and Lemma~\ref{lemma:drequal}, the last thing is to compare the definition of $K_{\phi}(1)$ with $K_\phi(x)$. Note that if we define
\begin{align*}
  A_\phi'(x) := \Big\{z \in \cZ: \lVert \hat\phi(x) - \hat\phi(z) \rVert < \hat\Phi_M(z)\Big\},~~ K_\phi'(x) := \sum_{j=1}^{N_1} \ind\Big(Z_j \in A_\phi'(x)\Big),
\end{align*}
as long as the ties are broken in arbitrary way, we can check that $K_\phi'(X_1) \le K_{\phi}(1) \le K_\phi(X_1)$. Note that all the previous results for $K_\phi(x)$ are also hold for $K_\phi'(x)$. Then the proof is complete.

\section{Proofs of Auxiliary Results}

\subsection{Proof of Lemma~\ref{lemma:mbc,bias}}

Note that for any $i \in \zahl{n}$ and $\omega =0,1$,
\begin{align*}
  &\Big\lvert\mu_{\phi,\omega}(U_{\phi,\omega,i}) -  \mu_{\phi,\omega}(U_{\phi,\omega,j}) - \hat{\mu}_{\phi,\omega}(\hat{U}_{\phi,\omega,i}) + \hat{\mu}_{\phi,\omega}(\hat{U}_{\phi,\omega,j})\Big\rvert\\
  &\le \Big\lvert\mu_{\phi,\omega}(U_{\phi,\omega,i}) -  \mu_{\phi,\omega}(U_{\phi,\omega,j}) - \hat{\mu}_{\phi,\omega}(U_{\phi,\omega,i}) + \hat{\mu}_{\phi,\omega}(U_{\phi,\omega,j})\Big\rvert\\
  &+ \Big\lvert\hat{\mu}_{\phi,\omega}(U_{\phi,\omega,i}) -  \hat{\mu}_{\phi,\omega}(U_{\phi,\omega,j}) - \hat{\mu}_{\phi,\omega}(\hat{U}_{\phi,\omega,i}) + \hat{\mu}_{\phi,\omega}(\hat{U}_{\phi,\omega,j})\Big\rvert.
  \yestag\label{eq:mbc,bias6}
\end{align*}

We can also decompose it in another way:
\begin{align*}
  &\Big\lvert\mu_{\phi,\omega}(U_{\phi,\omega,i}) -  \mu_{\phi,\omega}(U_{\phi,\omega,j}) - \hat{\mu}_{\phi,\omega}(\hat{U}_{\phi,\omega,i}) + \hat{\mu}_{\phi,\omega}(\hat{U}_{\phi,\omega,j})\Big\rvert\\
  &\le \Big\lvert\mu_{\phi,\omega}(U_{\phi,\omega,i}) -  \mu_{\phi,\omega}(U_{\phi,\omega,j}) - \mu_{\phi,\omega}(\hat{U}_{\phi,\omega,i}) + \mu_{\phi,\omega}(\hat{U}_{\phi,\omega,j})\Big\rvert\\
  &+ \Big\lvert\mu_{\phi,\omega}(\hat{U}_{\phi,\omega,i}) - \mu_{\phi,\omega}(\hat{U}_{\phi,\omega,j}) - \hat{\mu}_{\phi,\omega}(\hat{U}_{\phi,\omega,i}) + \hat{\mu}_{\phi,\omega}(\hat{U}_{\phi,\omega,j})\Big\rvert.
\end{align*}

We consider the proof under \eqref{eq:mbc,bias6}, and the proof under the second decomposition is similar.

For the first term in \eqref{eq:mbc,bias6}, by Taylor expansion to $k$-th order with $k=\max\{\lfloor m/2 \rfloor,1\} + 1$,
\begin{align*}
    &\Big\lvert \mu_{\phi,\omega}(U_{\phi,\omega,j}) - \mu_{\phi,\omega}(U_{\phi,\omega,i}) - \sum_{\ell = 1}^{k-1} \sum_{t \in \Lambda_\ell} \frac{1}{t!}\partial^t \mu_{\phi,\omega}(U_{\phi,\omega,i}) (U_{\phi,\omega,j}-U_{\phi,\omega,i})^t \Big\rvert\notag \\
    &\le \max_{t \in \Lambda_k} \lVert \partial^t \mu_{\phi,\omega} \rVert_\infty \sum_{t \in \Lambda_k} \frac{1}{t!} \lVert U_{\phi,\omega,j}-U_{\phi,\omega,i} \rVert^k.
\end{align*}
In the same way,
\begin{align*}
  &\Big\lvert \hat{\mu}_{\phi,\omega}(U_{\phi,\omega,j}) - \hat{\mu}_{\phi,\omega}(U_{\phi,\omega,i}) - \sum_{\ell = 1}^{k-1} \sum_{t \in \Lambda_\ell} \frac{1}{t!} \partial^t \hat{\mu}_{\phi,\omega}(U_{\phi,\omega,i}) (U_{\phi,\omega,j}-U_{\phi,\omega,i})^t \Big\rvert\notag \\
  &\le \max_{t \in \Lambda_k} \lVert \partial^t \hat{\mu}_{\phi,\omega} \rVert_\infty \sum_{t \in \Lambda_k} \frac{1}{t!} \lVert U_{\phi,\omega,j}-U_{\phi,\omega,i} \rVert^k.
\end{align*}
We also have
\begin{align*}
  & \Big\lvert \sum_{\ell = 1}^{k-1} \sum_{t \in \Lambda_\ell} \frac{1}{t!}  (\partial^t \hat{\mu}_{\phi,\omega}(U_{\phi,\omega,i}) - \partial^t \mu_{\phi,\omega}(U_{\phi,\omega,i})) (U_{\phi,\omega,j}-U_{\phi,\omega,i})^t \Big\rvert \notag \\ 
  &\le \sum_{\ell = 1}^{k-1} \max_{t \in \Lambda_\ell} \lVert \partial^t \hat{\mu}_{\phi,\omega}(U_{\phi,\omega,i}) - \partial^t \mu_{\phi,\omega}(U_{\phi,\omega,i}) \rVert \sum_{t \in \Lambda_\ell} \frac{1}{t!}  \lVert U_{\phi,\omega,j}-U_{\phi,\omega,i} \rVert^\ell.
\end{align*}

For the second term in \eqref{eq:mbc,bias6}, by Taylor expansion,
\begin{align*}
  &\Big\lvert\hat{\mu}_{\phi,\omega}(U_{\phi,\omega,i}) -  \hat{\mu}_{\phi,\omega}(U_{\phi,\omega,j}) - \hat{\mu}_{\phi,\omega}(\hat{U}_{\phi,\omega,i}) + \hat{\mu}_{\phi,\omega}(\hat{U}_{\phi,\omega,j})\Big\rvert\\
  &= \Big\lvert \partial \hat{\mu}_{\phi,\omega}(\bar{u}_j)^\top (\hat{U}_{\phi,\omega,j} - U_{\phi,\omega,j}) - \partial \hat{\mu}_{\phi,\omega}(\bar{u}_i)^\top (\hat{U}_{\phi,\omega,i} - U_{\phi,\omega,i})\Big\rvert\\
  &\le \Big\lvert [\partial \hat{\mu}_{\phi,\omega}(\bar{u}_j) - \partial \hat{\mu}_{\phi,\omega}(\bar{u}_i)]^\top (\hat{U}_{\phi,\omega,i} - U_{\phi,\omega,i})\Big\rvert + \Big\lvert \partial \hat{\mu}_{\phi,\omega}(\bar{u}_j)^\top (\hat{U}_{\phi,\omega,j} - \hat{U}_{\phi,\omega,i} - U_{\phi,\omega,j} + U_{\phi,\omega,i}) \Big\rvert\\
  &\lesssim \max_{t \in \Lambda_2} \lVert \partial^t \hat{\mu}_{\phi,\omega} \rVert_\infty \lVert \bar{u}_j - \bar{u}_i \rVert \lVert \hat{U}_{\phi,\omega,i} - U_{\phi,\omega,i} \rVert + \lVert \partial \hat{\mu}_{\phi,\omega} \rVert_\infty \lVert \hat{U}_{\phi,\omega,j} - \hat{U}_{\phi,\omega,i} - U_{\phi,\omega,j} + U_{\phi,\omega,i} \rVert,
\end{align*}
where $\bar{u}_i$ is between $U_{\phi,\omega,i}$ and $\hat{U}_{\phi,\omega,i}$, and $\bar{u}_j$ is between $U_{\phi,\omega,j}$ and $\hat{U}_{\phi,\omega,j}$. Since $\lVert \hat{U}_{\phi,\omega,i} - U_{\phi,\omega,i} \rVert \le \lVert \hat{\phi}_\omega - \phi_\omega \rVert_\infty$ and $\lVert \bar{u}_j - \bar{u}_i \rVert \le \lVert U_{\phi,\omega,j} - U_{\phi,\omega,i}\rVert + \lVert \hat{\phi}_\omega - \phi_\omega \rVert_\infty$, we have
\begin{align*}
  &\Big\lvert\hat{\mu}_{\phi,\omega}(U_{\phi,\omega,i}) -  \hat{\mu}_{\phi,\omega}(U_{\phi,\omega,j}) - \hat{\mu}_{\phi,\omega}(\hat{U}_{\phi,\omega,i}) + \hat{\mu}_{\phi,\omega}(\hat{U}_{\phi,\omega,j})\Big\rvert\\
  &\lesssim \max_{t \in \Lambda_2} \lVert \partial^t \hat{\mu}_{\phi,\omega} \rVert_\infty[\lVert \hat{\phi}_\omega - \phi_\omega \rVert_\infty^2 + \lVert \hat{\phi}_\omega - \phi_\omega \rVert_\infty \lVert U_{\phi,\omega,j} - U_{\phi,\omega,i}\rVert] + \lVert \partial \hat{\mu}_{\phi,\omega} \rVert_\infty \lVert \hat{U}_{\phi,\omega,j} - \hat{U}_{\phi,\omega,i} - U_{\phi,\omega,j} + U_{\phi,\omega,i} \rVert.
\end{align*}

Since we have $\lvert \cJ_{\phi}(i) \rvert = M$ for any $i \in \zahl{n}$, then
\begin{align*}
  & \lvert B_n - \hat{B}_n \rvert \\
  &\le \frac{1}{n}\sum_{i=1}^n \frac{1}{M} \sum_{j \in \cJ_{\phi}(i)} \Big\lvert\mu_{\phi,1-D_i}(U_{\phi,1-D_i,i}) -  \mu_{\phi,1-D_i}(U_{\phi,1-D_i,j}) - \hat{\mu}_{\phi,1-D_i}(\hat{U}_{\phi,1-D_i,i}) + \hat{\mu}_{\phi,1-D_i}(\hat{U}_{\phi,1-D_i,j})\Big\rvert\\
  &\le \frac{1}{n}\sum_{i=1}^n \max_{j \in \cJ_{\phi}(i)} \Big\lvert\mu_{\phi,1-D_i}(U_{\phi,1-D_i,i}) -  \mu_{\phi,1-D_i}(U_{\phi,1-D_i,j}) - \hat{\mu}_{\phi,1-D_i}(\hat{U}_{\phi,1-D_i,i}) + \hat{\mu}_{\phi,1-D_i}(\hat{U}_{\phi,1-D_i,j})\Big\rvert\\
  &\lesssim \max_{\omega \in \{0,1\}} \Big(\max_{t \in \Lambda_k} \lVert \partial^t \mu_{\phi,\omega} \rVert_\infty + \max_{t \in \Lambda_k} \lVert \partial^t \hat{\mu}_{\phi,\omega} \rVert_\infty \Big) \Big(\frac{1}{n}\sum_{i=1}^n \max_{j \in \cJ_{\phi}(i)} \lVert U_{\phi,1-D_i,j}-U_{\phi,1-D_i,i} \rVert^k\Big)\\
  &+ \sum_{\ell = 1}^{k-1} \Big(\frac{1}{n}\sum_{i=1}^n \max_{t \in \Lambda_\ell} \lVert \partial^t \hat{\mu}_{\phi,1-D_i}(U_{\phi,1-D_i,i}) - \partial^t \mu_{\phi,1-D_i}(U_{\phi,1-D_i,i}) \rVert \max_{j \in \cJ_{\phi}(i)} \lVert U_{\phi,1-D_i,j}-U_{\phi,1-D_i,i} \rVert^\ell\Big)\\
  &+ \max_{\omega \in \{0,1\}} \max_{t \in \Lambda_2} \lVert \partial^t \hat{\mu}_{\phi,\omega} \rVert_\infty \Big[\lVert \hat{\phi}_\omega - \phi_\omega \rVert_\infty^2 + \lVert \hat{\phi}_\omega - \phi_\omega \rVert_\infty \Big(\frac{1}{n}\sum_{i=1}^n \max_{j \in \cJ_{\phi}(i)} \lVert U_{\phi,1-D_i,j} - U_{\phi,1-D_i,i}\rVert\Big)\Big] \\
  &+ \max_{\omega \in \{0,1\}}\lVert \partial \hat{\mu}_{\phi,\omega} \rVert_\infty \Big(\frac{1}{n}\sum_{i=1}^n \max_{j \in \cJ_{\phi}(i)} \lVert \hat{U}_{\phi,1-D_i,j} - \hat{U}_{\phi,1-D_i,i} - U_{\phi,1-D_i,j} + U_{\phi,1-D_i,i} \rVert\Big).
  \yestag\label{eq:mbc,bias1}
\end{align*}

For any $i \in \zahl{n}$, let $\tilde{\cJ}_{\phi}(i)$ be the index set of $M$-NNs of $U_{\phi,1-D_i,i}$ in $\{U_{\phi,1-D_i,j}:D_j=1-D_i\}_{j=1}^n$ with ties broken in arbitrary way. Then
\begin{align*}
  \max_{j \in \cJ_{\phi}(i)} \lVert U_{\phi,1-D_i,j} - U_{\phi,1-D_i,i}\rVert &\le\max_{j \in \cJ_{\phi}(i)} \lVert \hat{U}_{\phi,1-D_i,j} - \hat{U}_{\phi,1-D_i,i}\rVert + 2\lVert \hat{\phi}_{1-D_i} - \phi_{1-D_i} \rVert_\infty\\
  &\le \max_{j \in \tilde{\cJ}_{\phi}(i)} \lVert \hat{U}_{\phi,1-D_i,j} - \hat{U}_{\phi,1-D_i,i}\rVert + 2\lVert \hat{\phi}_{1-D_i} - \phi_{1-D_i} \rVert_\infty\\
  &\le \max_{j \in \tilde{\cJ}_{\phi}(i)} \lVert U_{\phi,1-D_i,j} - U_{\phi,1-D_i,i}\rVert + 4\lVert \hat{\phi}_{1-D_i} - \phi_{1-D_i} \rVert_\infty.
\end{align*}

By \citet[Lemma 14.1]{li2023nonparametric}, as long as the density of $U_{\phi,\omega}$ is continuous for $\omega =0,1$, we have for any positive integer $p$,
\begin{align*}
  \E\Big[\frac{1}{n}\sum_{i=1}^n \max_{j \in \tilde{\cJ}_{\phi}(i)} \lVert U_{\phi,1-D_i,j}-U_{\phi,1-D_i,i} \rVert^p\Big] \lesssim \Big(\frac{M}{n}\Big)^{p/m}.
\end{align*}

Then for any positive integer $p$, by Assumption~\ref{asp:se3-g}, we have
\begin{align*}
  &\frac{1}{n}\sum_{i=1}^n \max_{j \in \cJ_{\phi}(i)} \lVert U_{\phi,1-D_i,j}-U_{\phi,1-D_i,i} \rVert^p\\
   &\lesssim  \frac{1}{n}\sum_{i=1}^n \Big(\max_{j \in \tilde{\cJ}_{\phi}(i)} \lVert U_{\phi,1-D_i,j}-U_{\phi,1-D_i,i} \rVert^p + \lVert \hat{\phi}_{1-D_i} - \phi_{1-D_i} \rVert_\infty^p\Big)\\
  &= O_\P((M/n)^{p/m} + n^{-p/2}).
  \yestag\label{eq:mbc,bias2}
\end{align*}

For any positive integer $\ell \in \zahl{k-1}$, by Assumption~\ref{asp:se2-g}, we have
\begin{align*}
  &\frac{1}{n}\sum_{i=1}^n \max_{t \in \Lambda_\ell} \lVert \partial^t \hat{\mu}_{\phi,1-D_i}(U_{\phi,1-D_i,i}) - \partial^t \mu_{\phi,1-D_i}(U_{\phi,1-D_i,i}) \rVert \max_{j \in \cJ_{\phi}(i)} \lVert U_{\phi,1-D_i,j}-U_{\phi,1-D_i,i} \rVert^\ell\\
  &\le \max_{\omega \in \{0,1\}}\max_{t \in \Lambda_\ell} \lVert \partial^t \hat{\mu}_{\phi,\omega} - \partial^t \mu_{\phi,\omega} \rVert_\infty\Big(\frac{1}{n}\sum_{i=1}^n \max_{j \in \cJ_{\phi}(i)} \lVert U_{\phi,1-D_i,j}-U_{\phi,1-D_i,i} \rVert^{\ell}\Big)\\
  &= O_\P(n^{-\gamma_\ell}((M/n)^{\ell/m}+ n^{-\ell/2})).
  \yestag\label{eq:mbc,bias3}
\end{align*}

For any $\epsilon>0$ and $\omega =0,1$, we have
\begin{align*}
  &\P\Big(\frac{1}{n}\sum_{i=1}^n \max_{j \in \cJ_{\phi}(i)} \lVert \hat{U}_{\phi,\omega,j} - \hat{U}_{\phi,\omega,i} - U_{\phi,\omega,j} + U_{\phi,\omega,i} \rVert \ge n^{-\frac{1}{2}} \epsilon\Big)\\
  &\le \P\Big(\sqrt{n} \sup_{x,y \in \cX, \lVert \phi_\omega(x)-\phi_\omega(y) \rVert \le \delta} \lVert (\hat{\phi}_\omega - \phi_\omega)(x) - (\hat{\phi}_\omega - \phi_\omega)(y)\rVert \ge \epsilon\Big) + \P\Big(\max_{i \in \zahl{n}} \max_{j \in \cJ_{\phi}(i)} \lVert U_{\phi,\omega,j} - U_{\phi,\omega,i} \rVert \ge \delta \Big),
\end{align*}
which holds for any $\delta>0$.

Taking $n \to \infty$ and then $\delta \to 0$, by Assumption~\ref{asp:se3-g} and $M/n \to 0$, we have
\begin{align}\label{eq:mbc,bias4}
  \frac{1}{n}\sum_{i=1}^n \max_{j \in \cJ_{\phi}(i)} \lVert \hat{U}_{\phi,\omega,j} - \hat{U}_{\phi,\omega,i} - U_{\phi,\omega,j} + U_{\phi,\omega,i} \rVert = o_\P(n^{-1/2}). 
\end{align}

Plugging \eqref{eq:mbc,bias2}, \eqref{eq:mbc,bias3}, \eqref{eq:mbc,bias4} into \eqref{eq:mbc,bias1} and using Assumption~\ref{asp:se2-g} yields
\begin{align*}
  & \lvert B_n - \hat{B}_n \rvert\\
  &\lesssim  O_\P((M/n)^{k/m} + n^{-k/2}) + \sum_{\ell = 1}^{k-1} O_\P(n^{-\gamma_\ell}((M/n)^{\ell/m}+ n^{-\ell/2})) + O_\P((M/n)^{1/m}n^{-1/2}) + o_\P(n^{-1/2}).
\end{align*}

This completes the proof by the selection of $M$.

\subsection{Proof of Theorem~\ref{lemma:moment}}

Note that by the multinomial theorem,
\begin{align*}
  & \E\big[(\hat{r}_\phi(x))^p\big] = \E\Big[\Big(\frac{N_0}{N_1} \frac{K_\phi(x)}{M}\Big)^p\Big] = \Big(\frac{N_0}{N_1 M}\Big)^p \E\Big[\Big(\sum_{j=1}^{N_1} \ind\Big(Z_j \in A_\phi(x)\Big)\Big)^p\Big]\\
  &= \Big(\frac{N_0}{N_1 M}\Big)^p \sum_{p_1+\cdots+p_{N_1}=p;~p_1,\ldots,p_{N_1}\ge0} \binom{p}{p_1,\ldots,p_{N_1}}\E\Big[\prod_{j=1}^{N_1} \ind\Big(Z_j \in A_\phi(x)\Big)^{p_j}\Big]\\
  &= \Big(\frac{N_0}{N_1 M}\Big)^p \sum_{p_1+\cdots+p_{N_1}=p;~p_1,\ldots,p_{N_1}\ge0} \binom{p}{p_1,\ldots,p_{N_1}}\P\Big(Z_j \in A_\phi(x): p_j > 0\Big).
\end{align*}
Then by Lemma~\ref{lemma:moment,p}, we have
\begin{align*}
  \lim_{N_0 \to \infty} \Big(\frac{N_0}{M}\Big)^{\sum_{j=1}^{N_1} \ind(p_j > 0)} \P\Big(Z_j \in A_\phi(x): p_j > 0\Big) = \big[r(\phi(x))\big]^{\sum_{j=1}^{N_1} \ind(p_j > 0)}.
\end{align*}

Note that for $p_1,\ldots,p_{N_1}\ge0$ with $p_1+\cdots+p_{N_1}=p$, the number of terms such that $\sum_{j=1}^{N_1} \ind(p_j > 0) = k$ is of order $N_1^k$ for any $k \in \zahl{p}$. Also note that $\binom{p}{p_1,\ldots,p_{N_1}}$ is bounded. Therefore if $MN_1/N_0 \to \infty$, we have
\begin{align*}
  \lim_{N_0\to\infty} \E\big[(\hat{r}_\phi(x))^p\big] = \lim_{N_0\to\infty} \frac{1}{N_1^p} \binom{N_1}{p} \binom{p}{1,\ldots,1} \big[r(\phi(x))\big]^p = \big[r(\phi(x))\big]^p.
\end{align*}
This completes the proof.

\subsection{Proof of Lemma~\ref{lemma:moment,p}}

We only consider those $x \in \cX$ such that $f_0(\phi(x))>0$ and $\phi(x)$ is a continuous point of $f_0$ and $f_1$. We seperate the proof into two cases depending on whether $f_1(\phi(x))$ is zero.

{\bf Part I.} We first consider the simple case where $p=1$ and Assumption~\ref{asp:proof1} holds for $p$.

{\bf Case I.} $f_1(\phi(x))>0$.
Since $\phi(x)$ is a continuous point of $f_0$ and $f_1$, for any $\epsilon \in (0,1)$, there exists some $\delta = \delta_x>0$ such that for any $z \in \cX$ with $\lVert \phi(z) - \phi(x) \rVert \le 3\delta$, we have $\lvert f_0(\phi(z)) - f_0(\phi(x)) \rvert \le \epsilon f_0(\phi(x))$ and $\lvert f_1(\phi(z)) - f_1(\phi(x)) \rvert \le \epsilon f_1(\phi(x))$. Denote the closed ball in $\bR^m$ centered at $x$ with radius $\delta$ by $B_{x,\delta}$, and the Lebesgue measure by $\lambda$. Then for any $z \in \cX$ with $\lVert \phi(z) - \phi(x) \rVert \le \delta$, we have
\begin{align*}
  & \Big\lvert \frac{\nu_0(B_{\phi(x),\lVert \phi(z) - \phi(x) \rVert})}{\lambda(B_{\phi(x),\lVert \phi(z) - \phi(x) \rVert})} - f_0(\phi(x)) \Big\rvert \le \epsilon f_0(\phi(x)),~~ \Big\lvert \frac{\nu_0(B_{\phi(z),\lVert \phi(z) - \phi(x) \rVert})}{\lambda(B_{\phi(z),\lVert \phi(z) - \phi(x) \rVert})} - f_0(\phi(x)) \Big\rvert \le \epsilon f_0(\phi(x)),\\
  & \Big\lvert \frac{\nu_1(B_{\phi(x),\lVert \phi(z) - \phi(x) \rVert})}{\lambda(B_{\phi(x),\lVert \phi(z) - \phi(x) \rVert})} - f_1(\phi(x)) \Big\rvert \le \epsilon f_1(\phi(x)),~~ \Big\lvert \frac{\nu_1(B_{\phi(z),\lVert \phi(z) - \phi(x) \rVert})}{\lambda(B_{\phi(z),\lVert \phi(z) - \phi(x) \rVert})} - f_1(\phi(x)) \Big\rvert \le \epsilon f_1(\phi(x)).
\end{align*}
Accordingly, if $\lVert \phi(z) - \phi(x) \rVert \le \delta$, we have
\[
  \frac{1-\epsilon}{1+\epsilon} \frac{f_0(\phi(x))}{f_1(\phi(x))} \le \frac{\nu_0(B_{\phi(z),\lVert \phi(z) - \phi(x) \rVert})}{\lambda(B_{\phi(z),\lVert \phi(z) - \phi(x) \rVert})} \frac{\lambda(B_{\phi(x),\lVert \phi(z) - \phi(x) \rVert})}{\nu_1(B_{\phi(x),\lVert \phi(z) - \phi(x) \rVert})}  \le \frac{1+\epsilon}{1-\epsilon} \frac{f_0(\phi(x))}{f_1(\phi(x))}.
\]
Since $\lambda(B_{\phi(z),\lVert \phi(z) - \phi(x) \rVert}) = \lambda(B_{\phi(x),\lVert \phi(z) - \phi(x) \rVert})$, we then have
\begin{align*}
  \frac{1-\epsilon}{1+\epsilon} \frac{f_0(\phi(x))}{f_1(\phi(x))} \le \frac{\nu_0(B_{\phi(z),\lVert \phi(z) - \phi(x) \rVert})}{\nu_1(B_{\phi(x),\lVert \phi(z) - \phi(x) \rVert})} \le \frac{1+\epsilon}{1-\epsilon} \frac{f_0(\phi(x))}{f_1(\phi(x))}.
\end{align*}

On the other hand, consider any $\epsilon' \in (0,1)$. For any $z \in \cX$ such that $\lVert \phi(z) - \phi(x) \rVert > \delta$, as long as $\epsilon'$ small enough such that $\epsilon' {\rm diam}(\cX) < \delta/2$, where ${\rm diam}(\cX)$ is the diameter of $\cX$, we have $B_{y,\delta/2} \subset B_{\phi(z),\lVert \phi(z) - \phi(x) \rVert - \delta/2} \subset B_{\phi(z),(1-\epsilon')\lVert \phi(z) - \phi(x) \rVert}$, where $y\in\bR^m$ is taken such that $y$ is the intersection point of the surface of $B_{\phi(x),\delta}$ and the line connecting $\phi(z)$ and $\phi(x)$. Then
\[
  \nu_0(B_{\phi(z),(1-\epsilon')\lVert \phi(z) - \phi(x) \rVert}) \ge \nu_0(B_{y,\delta/2}) \ge (1-\epsilon)f_0(\phi(x)) \lambda(B_{y,\delta/2}) = (1-\epsilon)f_0(\phi(x)) \lambda(B_{0,\delta/2}).
\]

Let $\eta_N = 4\log(N_0/M)$. Since $M \log N_0 /N_0 \to 0$, we can take $N_0$ large enough so that
\[
  \eta_N \frac{M}{N_0} = 4\frac{M}{N_0} \log\Big(\frac{N_0}{M}\Big) < (1-\epsilon)f_0(\phi(x)) \lambda(B_{0,\delta/2}).
\]

Then for any $z \in \cX$ such that $\nu_0(B_{\phi(z),(1-\epsilon')\lVert \phi(z) - \phi(x) \rVert}) \le \eta_N M/N_0$, we have $\lVert \phi(z) - \phi(x) \rVert < \delta$ since otherwise it would contradict the selection of $\eta_N$.

{\bf Upper bound.} Let $\Phi_M(z)$ be the $M$-th order statistics of $\{\lVert \phi(X_i) - \phi(z) \rVert\}_{i=1}^{N_0}$. By the definition of $A_\phi(x)$, we have for any $\epsilon' \in (0,1)$,
\begin{align*}
  &\P\Big(Z_1 \in A_\phi(x)\Big) = \P\Big(\lVert \hat\phi(x) - \hat\phi(Z_1) \rVert \le \hat\Phi_M(Z_1)\Big) \\
  &\le \P\Big(\lVert \phi(x) - \phi(Z_1) \rVert - 2\lVert \hat\phi-\phi \rVert_\infty \le \Phi_M(Z_1) + 2\lVert \hat\phi-\phi \rVert_\infty\Big)\\
  &= \P\Big(\lVert \phi(x) - \phi(Z_1) \rVert - 4\lVert \hat\phi-\phi \rVert_\infty \le \Phi_M(Z_1), 4\lVert \hat\phi-\phi \rVert_\infty \le \epsilon' \lVert \phi(x) - \phi(Z_1) \rVert\Big)\\
  &+ \P\Big(\lVert \phi(x) - \phi(Z_1) \rVert - 4\lVert \hat\phi-\phi \rVert_\infty \le \Phi_M(Z_1), 4\lVert \hat\phi-\phi \rVert_\infty > \epsilon' \lVert \phi(x) - \phi(Z_1) \rVert\Big).
  \yestag\label{eq:momentcatch1}
\end{align*}

For the first term in \eqref{eq:momentcatch1}, note that $[\phi(X_i)]_{i=1}^{N_0}$ are i.i.d from $\nu_0$, and then $\nu_0(B_{\phi(Z_1),\lVert \phi(X_i) - \phi(Z_1) \rVert})$ are i.i.d from $U(0,1)$ and are independent of $Z_1$ by the probability integral transform. Then
\begin{align*}
  &\P\Big(\lVert \phi(x) - \phi(Z_1) \rVert - 4\lVert \hat\phi-\phi \rVert_\infty \le \Phi_M(Z_1), 4\lVert \hat\phi-\phi \rVert_\infty \le \epsilon' \lVert \phi(x) - \phi(Z_1) \rVert\Big)\\
  &\le \P\Big((1-\epsilon')\lVert \phi(x) - \phi(Z_1) \rVert \le \Phi_M(Z_1)\Big)\\
  &= \P\Big(\nu_0(B_{\phi(Z_1),(1-\epsilon')\lVert \phi(x) - \phi(Z_1) \rVert }) \le \nu_0(B_{\phi(Z_1),\Phi_M(Z_1) })\Big)\\
  &\le \P\Big(\nu_0(B_{\phi(Z_1),(1-\epsilon')\lVert \phi(x) - \phi(Z_1) \rVert }) \le \nu_0(B_{\phi(Z_1),\Phi_M(Z_1) }) \le \eta_N \frac{M}{N_0} \Big) + \P\Big(U_{(M)} > \eta_N \frac{M}{N_0}\Big),
\end{align*}
where $U_{(M)}$ is the $M$-th order statistic of $N_0$ independent random variables from $U(0,1)$.

By the selection of $\eta_N$, and taking $\epsilon'$ small and $N_0$ large enough, we have
\begin{align*}
  &\P\Big(\nu_0(B_{\phi(Z_1),(1-\epsilon')\lVert \phi(x) - \phi(Z_1) \rVert }) \le \nu_0(B_{\phi(Z_1),\Phi_M(Z_1) }) \le \eta_N \frac{M}{N_0} \Big)\\
  &\le \P\Big(\nu_0(B_{\phi(Z_1),(1-\epsilon')\lVert \phi(x) - \phi(Z_1) \rVert }) \le \nu_0(B_{\phi(Z_1),\Phi_M(Z_1) }), \lVert \phi(x) - \phi(Z_1) \rVert \le \delta \Big).
\end{align*}

Under the event $\{\lVert \phi(x) - \phi(Z_1) \rVert \le \delta\}$, we have
\begin{align*}
  &\nu_0(B_{\phi(Z_1),\lVert \phi(x) - \phi(Z_1) \rVert }) - \nu_0(B_{\phi(Z_1),(1-\epsilon')\lVert \phi(x) - \phi(Z_1) \rVert})\\
  &= \int_{B_{\phi(Z_1),\lVert \phi(x) - \phi(Z_1) \rVert } \setminus B_{\phi(Z_1),(1-\epsilon')\lVert \phi(x) - \phi(Z_1) \rVert}} f_0(y) \d y\\
  &\le (1+\epsilon)f_0(\phi(x)) \lambda(B_{\phi(Z_1),\lVert \phi(x) - \phi(Z_1) \rVert } \setminus B_{\phi(Z_1),(1-\epsilon')\lVert \phi(x) - \phi(Z_1) \rVert})\\
  &= (1+\epsilon)f_0(\phi(x)) V_m [1-(1-\epsilon')^d] \lVert \phi(x) - \phi(Z_1) \rVert^d\\
  &\le (1+\epsilon) f_0(\phi(x)) V_m d \epsilon' \lVert \phi(x) - \phi(Z_1) \rVert^d \\
  &= (1+\epsilon) f_0(\phi(x)) d \epsilon' \lambda(B_{\phi(x),\lVert \phi(x) - \phi(Z_1) \rVert })\\
  &\le \frac{(1+\epsilon)f_0(\phi(x)) d \epsilon'}{(1-\epsilon)f_1(\phi(x))} \nu_1(B_{\phi(x),\lVert \phi(x) - \phi(Z_1) \rVert }),
\end{align*}
where $V_m$ is the Lebesgue measure of the $m-$dimensional unit ball, and
\begin{align*}
  \nu_0(B_{\phi(Z_1),\lVert \phi(x) - \phi(Z_1) \rVert }) \ge \frac{(1-\epsilon)f_0(\phi(x))}{(1+\epsilon)f_1(\phi(x))} \nu_1(B_{\phi(x),\lVert \phi(x) - \phi(Z_1) \rVert }).
\end{align*}

From the probability integral transform, we have $\nu_1(B_{\phi(x),\lVert \phi(x) - \phi(Z_1) \rVert })$ is from $U(0,1)$ and then for $U \sim U(0,1)$,
\begin{align*}
  &\P\Big(\nu_0(B_{\phi(Z_1),(1-\epsilon')\lVert \phi(x) - \phi(Z_1) \rVert }) \le \nu_0(B_{\phi(Z_1),\Phi_M(Z_1) }), \lVert \phi(x) - \phi(Z_1) \rVert \le \delta \Big)\\
  &\le \P\Big(\Big(\frac{1-\epsilon}{1+\epsilon} - \frac{1+\epsilon}{1-\epsilon} d \epsilon' \Big) \frac{f_0(\phi(x))}{f_1(\phi(x))}\nu_1(B_{\phi(x),\lVert \phi(x) - \phi(Z_1) \rVert }) \le \nu_0(B_{\phi(Z_1),\Phi_M(Z_1) })\Big)\\
  &= \P\Big(\Big(\frac{1-\epsilon}{1+\epsilon} - \frac{1+\epsilon}{1-\epsilon} d \epsilon' \Big) \frac{f_0(\phi(x))}{f_1(\phi(x))} U \le U_{(M)} \Big).
\end{align*}
We can check that
\begin{align*}
  \lim_{N_0 \to \infty} \frac{N_0}{M} \P\Big(\Big(\frac{1-\epsilon}{1+\epsilon} - \frac{1+\epsilon}{1-\epsilon} d \epsilon' \Big) \frac{f_0(\phi(x))}{f_1(\phi(x))} U \le U_{(M)} \Big) = \Big(\frac{1-\epsilon}{1+\epsilon} - \frac{1+\epsilon}{1-\epsilon} d \epsilon' \Big)^{-1} \frac{f_1(\phi(x))}{f_0(\phi(x))}.
\end{align*}
Note that $\eta_N \to \infty$ as $N_0 \to \infty$ since $M/N_0 \to 0$. Then from the Chernoff bound and for $N_0$ sufficiently large, we have
\begin{align*}
  & \frac{N_0}{M} \P\Big(U_{(M)} > \eta_N \frac{M}{N_0}\Big) = \frac{N_0}{M} \P\Big({\rm Bin}\Big(N_0,\eta_N \frac{M}{N_0}\Big) \le M \Big)\\
  &\le \frac{N_0}{M}\exp\Big((1+\log \eta_N - \eta_N) M\Big) \le \frac{N_0}{M} \exp\Big(-\frac{1}{2} \eta_N M \Big) = \Big(\frac{N_0}{M}\Big)^{1-2M}.
\end{align*}
Since $M /N_0 \to 0$ and $M \ge 1$, we then obtain
\[
  \lim_{N_0 \to \infty} \frac{N_0}{M} \P\Big(U_{(M)} > \eta_N \frac{M}{N_0}\Big) = 0.
\]
Then we obtain
\begin{align*}
  &\limsup_{N_0 \to \infty} \frac{N_0}{M} \P\Big(\lVert \phi(x) - \phi(Z_1) \rVert - 4\lVert \hat\phi-\phi \rVert_\infty \le \Phi_M(Z_1), 4\lVert \hat\phi-\phi \rVert_\infty \le \epsilon' \lVert \phi(x) - \phi(Z_1) \rVert\Big)\\
  &\le \Big(\frac{1-\epsilon}{1+\epsilon} - \frac{1+\epsilon}{1-\epsilon} d \epsilon' \Big)^{-1} \frac{f_1(\phi(x))}{f_0(\phi(x))}.
\end{align*}

For the second term in \eqref{eq:momentcatch1}, we have for any $\delta>0$,
\begin{align*}
  &\P\Big(\lVert \phi(x) - \phi(Z_1) \rVert - 4\lVert \hat\phi-\phi \rVert_\infty \le \Phi_M(Z_1), 4\lVert \hat\phi-\phi \rVert_\infty > \epsilon' \lVert \phi(x) - \phi(Z_1) \rVert\Big) \\
  &\le \P\Big( \lVert \phi(x) - \phi(Z_1) \rVert < 4\lVert \hat\phi-\phi \rVert_\infty/ \epsilon'\Big)\\
  &\le \P\Big( \lVert \phi(x) - \phi(Z_1) \rVert < 4\lVert \hat\phi-\phi \rVert_\infty/ \epsilon' \le \delta\Big) + \P\Big(4\lVert \hat\phi-\phi \rVert_\infty/ \epsilon' > \delta\Big).
\end{align*}
Note that $\hat\phi$ may depend on $Z_1$. Recall that $\hat\phi_{Z_1 \to z}$ is the estimator of $\phi$ replacing $Z_1$ by some $z \in \cZ$. Then $\sup_{z \in \cZ} \lVert \hat\phi_{Z_1 \to z}-\phi \rVert_\infty$ is independent of $Z_1$ and then we have
\begin{align*}
  & \P\Big( \lVert \phi(x) - \phi(Z_1) \rVert < 4\lVert \hat\phi-\phi \rVert_\infty/ \epsilon' \le \delta\Big)\\
  &\le \P\Big(\lambda(B_{\phi(x),\lVert \phi(x) - \phi(Z_1) \rVert }) \le V_m(4/\epsilon')^d \lVert \hat\phi-\phi \rVert_\infty^d, \lVert \phi(x) - \phi(Z_1) \rVert \le \delta \Big)\\
  &\le \P\Big(\nu_1(B_{\phi(x),\lVert \phi(x) - \phi(Z_1) \rVert }) \le (1+\epsilon) f_1(\phi(x)) V_m(4/\epsilon')^d \lVert \hat\phi-\phi \rVert_\infty^d \Big)\\
  &\le \P\Big(\nu_1(B_{\phi(x),\lVert \phi(x) - \phi(Z_1) \rVert }) \le (1+\epsilon) f_1(\phi(x)) V_m(4/\epsilon')^d \sup_{z \in \cZ} \lVert \hat\phi_{Z_1 \to z}-\phi \rVert_\infty^d \Big)\\
  &= \E\{[(1+\epsilon) f_1(\phi(x)) V_m(4/\epsilon')^d \sup_{z \in \cZ} \lVert \hat\phi_{Z_1 \to z}-\phi \rVert_\infty^d] \wedge 1\}.
\end{align*}
By the Markov inequality,
\begin{align*}
  \P\Big(4\lVert \hat\phi-\phi \rVert_\infty/ \epsilon' > \delta\Big) \le \delta^{-d} (4/\epsilon')^d \E[\lVert \hat\phi-\phi \rVert_\infty^d] \le \delta^{-d} (4/\epsilon')^d \E[\sup_{z \in \cZ} \lVert \hat\phi_{Z_1 \to z}-\phi \rVert_\infty^d].
\end{align*}
By Assumption~\ref{asp:proof1}, we have
\begin{align*}
  \lim_{N_0 \to \infty} \frac{N_0}{M} \P\Big(\lVert \phi(x) - \phi(Z_1) \rVert - 4\lVert \hat\phi-\phi \rVert_\infty \le \Phi_M(Z_1), 4\lVert \hat\phi-\phi \rVert_\infty > \epsilon' \lVert \phi(x) - \phi(Z_1) \rVert\Big) = 0.
\end{align*}

By \eqref{eq:momentcatch1} and $\epsilon,\epsilon'$ are arbitrary, we obtain
\begin{align*}
  \limsup_{N_0 \to \infty} \frac{N_0}{M} \P\Big(Z_1 \in A_\phi(x)\Big) \le \frac{f_1(\phi(x))}{f_0(\phi(x))}.
\end{align*}

{\bf Lower bound.} For any $\epsilon' \in (0,1)$, we have
\begin{align*}
  &\P\Big(Z_1 \in A_\phi(x)\Big) = \P\Big(\lVert \hat\phi(x) - \hat\phi(Z_1) \rVert \le \hat\Phi_M(Z_1)\Big) \\
  &\ge \P\Big(\lVert \phi(x) - \phi(Z_1) \rVert + 2\lVert \hat\phi-\phi \rVert_\infty \le \Phi_M(Z_1) - 2\lVert \hat\phi-\phi \rVert_\infty\Big)\\
  &\ge \P\Big(\lVert \phi(x) - \phi(Z_1) \rVert + 4\lVert \hat\phi-\phi \rVert_\infty \le \Phi_M(Z_1), 4\lVert \hat\phi-\phi \rVert_\infty \le \epsilon' \lVert \phi(x) - \phi(Z_1) \rVert\Big)\\
  &\ge \P\Big((1+\epsilon')\lVert \phi(x) - \phi(Z_1) \rVert \le \Phi_M(Z_1), 4\lVert \hat\phi-\phi \rVert_\infty \le \epsilon' \lVert \phi(x) - \phi(Z_1) \rVert\Big)\\
  &\ge \P\Big((1+\epsilon')\lVert \phi(x) - \phi(Z_1) \rVert \le \Phi_M(Z_1)\Big) - \P\Big(4\lVert \hat\phi-\phi \rVert_\infty > \epsilon' \lVert \phi(x) - \phi(Z_1) \rVert\Big)\\
  &= \P\Big(\nu_0(B_{\phi(Z_1),(1+\epsilon')\lVert \phi(x) - \phi(Z_1) \rVert }) \le \nu_0(B_{\phi(Z_1),\Phi_M(Z_1) })\Big) - \P\Big(4\lVert \hat\phi-\phi \rVert_\infty > \epsilon' \lVert \phi(x) - \phi(Z_1) \rVert\Big)\\
  &\ge \P\Big(\nu_0(B_{\phi(Z_1),(1+\epsilon')\lVert \phi(x) - \phi(Z_1) \rVert }) \le \nu_0(B_{\phi(Z_1),\Phi_M(Z_1) }) \le \eta_N \frac{M}{N_0} \Big) - \P\Big(4\lVert \hat\phi-\phi \rVert_\infty > \epsilon' \lVert \phi(x) - \phi(Z_1) \rVert\Big).
\end{align*}

Note that by the selection of $\eta_N$, and taking $\epsilon'$ small and $N_0$ large enough, we have
\begin{align*}
  & \P\Big(\nu_0(B_{\phi(Z_1),(1+\epsilon')\lVert \phi(x) - \phi(Z_1) \rVert }) \le \nu_0(B_{\phi(Z_1),\Phi_M(Z_1) }) \le \eta_N \frac{M}{N_0} \Big)\\
  &= \P\Big(\nu_0(B_{\phi(Z_1),(1+\epsilon')\lVert \phi(x) - \phi(Z_1) \rVert }) \le \nu_0(B_{\phi(Z_1),\Phi_M(Z_1) }) \le \eta_N \frac{M}{N_0}, \lVert \phi(x) - \phi(Z_1) \rVert \le \delta\Big).
\end{align*}
Under the event $\{\lVert \phi(x) - \phi(Z_1) \rVert \le \delta\}$, we have $B_{\phi(Z_1),(1+\epsilon')\lVert \phi(x) - \phi(Z_1) \rVert } \subset B_{\phi(x),3\delta}$, and then
\begin{align*}
  &\nu_0(B_{\phi(Z_1),(1+\epsilon')\lVert \phi(x) - \phi(Z_1) \rVert }) - \nu_0(B_{\phi(Z_1),\lVert \phi(x) - \phi(Z_1) \rVert })\\
  &= \int_{B_{\phi(Z_1),(1+\epsilon')\lVert \phi(x) - \phi(Z_1) \rVert }\setminus B_{\phi(Z_1),\lVert \phi(x) - \phi(Z_1) \rVert } } f_0(y) \d y\\
  &\le (1+\epsilon)f_0(\phi(x)) \lambda(B_{\phi(Z_1),(1+\epsilon')\lVert \phi(x) - \phi(Z_1) \rVert }\setminus B_{\phi(Z_1),\lVert \phi(x) - \phi(Z_1) \rVert })\\
  &= (1+\epsilon)f_0(\phi(x)) V_m [(1+\epsilon')^d-1] \lVert \phi(x) - \phi(Z_1) \rVert^d\\
  &\le (1+\epsilon) f_0(\phi(x)) V_m d \epsilon' (1+\epsilon')^{d-1} \lVert \phi(x) - \phi(Z_1) \rVert^d \\
  &= (1+\epsilon) f_0(\phi(x)) d \epsilon' (1+\epsilon')^{d-1} \lambda(B_{\phi(x),\lVert \phi(x) - \phi(Z_1) \rVert })\\
  &\le \frac{(1+\epsilon)f_0(\phi(x)) d \epsilon'(1+\epsilon')^{d-1}}{(1-\epsilon)f_1(\phi(x))} \nu_1(B_{\phi(x),\lVert \phi(x) - \phi(Z_1) \rVert }),
\end{align*}
and
\begin{align*}
  \nu_0(B_{\phi(Z_1),\lVert \phi(x) - \phi(Z_1) \rVert }) \le \frac{(1+\epsilon)f_0(\phi(x))}{(1-\epsilon)f_1(\phi(x))} \nu_1(B_{\phi(x),\lVert \phi(x) - \phi(Z_1) \rVert }).
\end{align*}
Then
\begin{align*}
  & \P\Big(\nu_0(B_{\phi(Z_1),(1+\epsilon')\lVert \phi(x) - \phi(Z_1) \rVert }) \le \nu_0(B_{\phi(Z_1),\Phi_M(Z_1) }) \le \eta_N \frac{M}{N_0}, \lVert \phi(x) - \phi(Z_1) \rVert \le \delta\Big)\\
  &\ge \P\Big(\frac{(1+\epsilon)f_0(\phi(x))}{(1-\epsilon)f_1(\phi(x))} \Big(1+d \epsilon'(1+\epsilon')^{d-1}\Big) \nu_1(B_{\phi(x),\lVert \phi(x) - \phi(Z_1) \rVert })\le \nu_0(B_{\phi(Z_1),\Phi_M(Z_1) }) \le \eta_N \frac{M}{N_0}, \lVert \phi(x) - \phi(Z_1) \rVert \le \delta\Big)\\
  &= \P\Big(\frac{(1+\epsilon)f_0(\phi(x))}{(1-\epsilon)f_1(\phi(x))} \Big(1+d \epsilon'(1+\epsilon')^{d-1}\Big) \nu_1(B_{\phi(x),\lVert \phi(x) - \phi(Z_1) \rVert })\le \nu_0(B_{\phi(Z_1),\Phi_M(Z_1) }) \le \eta_N \frac{M}{N_0}\Big)\\
  &\ge \P\Big(\frac{(1+\epsilon)f_0(\phi(x))}{(1-\epsilon)f_1(\phi(x))} \Big(1+d \epsilon'(1+\epsilon')^{d-1}\Big) \nu_1(B_{\phi(x),\lVert \phi(x) - \phi(Z_1) \rVert }) \le \nu_0(B_{\phi(Z_1),\Phi_M(Z_1) }) \Big) - \P\Big(U_{(M)} > \eta_N \frac{M}{N_0}\Big)\\
  &= \P\Big(\frac{(1+\epsilon)f_0(\phi(x))}{(1-\epsilon)f_1(\phi(x))} \Big(1+d \epsilon'(1+\epsilon')^{d-1}\Big) U \le U_{(M)} \Big) - \P\Big(U_{(M)} > \eta_N \frac{M}{N_0}\Big).
\end{align*}
The second last equality is from the fact that for $z \in \cX$ such that $\lVert \phi(z) - \phi(x) \rVert > \delta$,
\begin{align*}
  &\frac{(1+\epsilon)f_0(\phi(x))}{(1-\epsilon)f_1(\phi(x))} \Big(1+d \epsilon'(1+\epsilon')^{d-1}\Big) \nu_1(B_{\phi(x),\lVert \phi(x) - \phi(z) \rVert }) \ge \frac{(1+\epsilon)f_0(\phi(x))}{(1-\epsilon)f_1(\phi(x))} \nu_1(B_{\phi(x),\delta})\\
  &\ge \frac{(1+\epsilon)f_0(\phi(x))}{(1-\epsilon)f_1(\phi(x))} f_1(\phi(x)) (1-\epsilon) \lambda(B_{0,\delta}) > \eta_N \frac{M}{N_0}
\end{align*}
by the selection of $\eta_N$.

We can check that
\begin{align*}
  \lim_{N_0 \to \infty} \frac{N_0}{M} \P\Big(\frac{(1+\epsilon)f_0(\phi(x))}{(1-\epsilon)f_1(\phi(x))} \Big(1+d \epsilon'(1+\epsilon')^{d-1}\Big) U \le U_{(M)} \Big) = \frac{1-\epsilon}{1+\epsilon}\Big(1+d \epsilon'(1+\epsilon')^{d-1}\Big)^{-1} \frac{f_1(\phi(x))}{f_0(\phi(x))}.
\end{align*}
By $\epsilon,\epsilon'$ are arbitrary, we obtain
\begin{align*}
  \liminf_{N_0 \to \infty} \frac{N_0}{M} \P\Big(Z_1 \in A_\phi(x)\Big) \ge \frac{f_1(\phi(x))}{f_0(\phi(x))}.
\end{align*}

Combining the upper bound and the lower bound yields
\begin{align*}
  \lim_{N_0 \to \infty} \frac{N_0}{M} \P\Big(Z_1 \in A_\phi(x)\Big) = \frac{f_1(\phi(x))}{f_0(\phi(x))}.
\end{align*}

{\bf Case II.} $f_1(\phi(x))=0$.

For any $\epsilon \in (0,1)$, there exists some $\delta = \delta_x>0$ such that for any $z \in \cX$ with $\lVert \phi(z) - \phi(x) \rVert \le 3\delta$, we have $\lvert f_0(\phi(z)) - f_0(\phi(x)) \rvert \le \epsilon f_0(\phi(x))$ and $f_1(\phi(z)) \le \epsilon$. Then for any $z \in \cX$ with $\lVert \phi(z) - \phi(x) \rVert \le \delta$, we have
\begin{align*}
  \Big\lvert \frac{\nu_0(B_{\phi(z),\lVert \phi(z) - \phi(x) \rVert})}{\lambda(B_{\phi(z),\lVert \phi(z) - \phi(x) \rVert})} - f_0(\phi(x)) \Big\rvert \le \epsilon f_0(\phi(x)), ~~ \Big\lvert \frac{\nu_1(B_{\phi(x),\lVert \phi(z) - \phi(x) \rVert})}{\lambda(B_{\phi(x),\lVert \phi(z) - \phi(x) \rVert})}\Big\rvert \le \epsilon.
\end{align*}

We consider the same decomposition as \eqref{eq:momentcatch1}. For the first term in \eqref{eq:momentcatch1}, we still have
\begin{align*}
  & \P\Big(\lVert \phi(x) - \phi(Z_1) \rVert - 4\lVert \hat\phi-\phi \rVert_\infty \le \Phi_M(Z_1), 4\lVert \hat\phi-\phi \rVert_\infty \le \epsilon' \lVert \phi(x) - \phi(Z_1) \rVert\Big)\\
  &\le \P\Big(\nu_0(B_{\phi(Z_1),(1-\epsilon')\lVert \phi(x) - \phi(Z_1) \rVert }) \le \nu_0(B_{\phi(Z_1),\Phi_M(Z_1) }), \lVert \phi(x) - \phi(Z_1) \rVert \le \delta \Big) + \P\Big(U_{(M)} > \eta_N \frac{M}{N_0}\Big).
\end{align*}
Note that for any $z \in \cX$ with $\lVert \phi(x) - \phi(z) \rVert \le \delta$,
\begin{align*}
  \Big\lvert \frac{\nu_0(B_{\phi(z),(1-\epsilon')\lVert \phi(z) - \phi(x) \rVert})}{\lambda(B_{\phi(z),(1-\epsilon')\lVert \phi(z) - \phi(x) \rVert})} - f_0(\phi(x)) \Big\rvert \le \epsilon f_0(\phi(x)),
\end{align*}
and
\begin{align*}
  & \frac{\nu_1(B_{\phi(x),\lVert \phi(x) - \phi(z) \rVert })}{\lambda(B_{\phi(z),(1-\epsilon')\lVert \phi(z) - \phi(x) \rVert})} = \frac{\lambda(B_{\phi(x),\lVert \phi(z) - \phi(x) \rVert})}{\lambda(B_{\phi(z),(1-\epsilon')\lVert \phi(z) - \phi(x) \rVert})} \frac{\nu_1(B_{\phi(x),\lVert \phi(z) - \phi(x) \rVert})}{\lambda(B_{\phi(x),\lVert \phi(z) - \phi(x) \rVert})}\\
  &\le \frac{\lambda(B_{\phi(x),\lVert \phi(z) - \phi(x) \rVert})}{\lambda(B_{\phi(z),(1-\epsilon')\lVert \phi(z) - \phi(x) \rVert})} \epsilon = (1-\epsilon')^{-d} \epsilon.
\end{align*}
Then
\begin{align*}
  & \P\Big(\nu_0(B_{\phi(Z_1),(1-\epsilon')\lVert \phi(x) - \phi(Z_1) \rVert }) \le \nu_0(B_{\phi(Z_1),\Phi_M(Z_1) }), \lVert \phi(x) - \phi(Z_1) \rVert \le \delta \Big)\\
  &\le \P\Big((1-\epsilon')^{d} \epsilon^{-1}(1-\epsilon) f_0(\phi(x))\nu_1(B_{\phi(x),\lVert \phi(x) - \phi(Z_1) \rVert }) \le \nu_0(B_{\phi(Z_1),\Phi_M(Z_1) })\Big)\\
  &= \P\Big((1-\epsilon')^{d} \epsilon^{-1}(1-\epsilon) f_0(\phi(x)) U \le U_{(M)} \Big).
\end{align*}
We can check that
\begin{align*}
  \lim_{N_0 \to \infty} \frac{N_0}{M} \P\Big((1-\epsilon')^{d} \epsilon^{-1}(1-\epsilon) f_0(\phi(x)) U \le U_{(M)} \Big) = \epsilon (1-\epsilon')^{-d} (1-\epsilon)^{-1} \frac{1}{f_0(\phi(x))}.
\end{align*}
Then we obtain
\begin{align*}
  &\limsup_{N_0 \to \infty} \frac{N_0}{M} \P\Big(\lVert \phi(x) - \phi(Z_1) \rVert - 4\lVert \hat\phi-\phi \rVert_\infty \le \Phi_M(Z_1), 4\lVert \hat\phi-\phi \rVert_\infty \le \epsilon' \lVert \phi(x) - \phi(Z_1) \rVert\Big)\\
  &\le \epsilon (1-\epsilon')^{-d} (1-\epsilon)^{-1} \frac{1}{f_0(\phi(x))}.
\end{align*}

For the second term in \eqref{eq:momentcatch1}, we still have
\begin{align*}
  &\P\Big(\lVert \phi(x) - \phi(Z_1) \rVert - 4\lVert \hat\phi-\phi \rVert_\infty \le \Phi_M(Z_1), 4\lVert \hat\phi-\phi \rVert_\infty > \epsilon' \lVert \phi(x) - \phi(Z_1) \rVert\Big) \\
  &\le \P\Big( \lVert \phi(x) - \phi(Z_1) \rVert < 4\lVert \hat\phi-\phi \rVert_\infty/ \epsilon' \le \delta\Big) + \P\Big(4\lVert \hat\phi-\phi \rVert_\infty/ \epsilon' > \delta\Big).
\end{align*}
Note that
\begin{align*}
  & \P\Big( \lVert \phi(x) - \phi(Z_1) \rVert < 4\lVert \hat\phi-\phi \rVert_\infty/ \epsilon' \le \delta\Big)\\
  &\le \P\Big(\lambda(B_{\phi(x),\lVert \phi(x) - \phi(Z_1) \rVert }) \le V_m(4/\epsilon')^d \lVert \hat\phi-\phi \rVert_\infty^d, \lVert \phi(x) - \phi(Z_1) \rVert \le \delta \Big)\\
  &\le \P\Big(\nu_1(B_{\phi(x),\lVert \phi(x) - \phi(Z_1) \rVert }) \le \epsilon V_m(4/\epsilon')^d \lVert \hat\phi-\phi \rVert_\infty^d \Big)\\
  &\le \P\Big(\nu_1(B_{\phi(x),\lVert \phi(x) - \phi(Z_1) \rVert }) \le \epsilon V_m(4/\epsilon')^d \sup_{z \in \cZ} \lVert \hat\phi_{Z_1 \to z}-\phi \rVert_\infty^d \Big)\\
  &= \E\{[\epsilon  V_m(4/\epsilon')^d \sup_{z \in \cZ} \lVert \hat\phi_{Z_1 \to z}-\phi \rVert_\infty^d] \wedge 1\}.
\end{align*}
By $\epsilon,\epsilon'$ are arbitrary, we obtain
\begin{align*}
  \lim_{N_0 \to \infty} \frac{N_0}{M} \P\Big(Z_1 \in A_\phi(x)\Big) = 0 = \frac{f_1(\phi(x))}{f_0(\phi(x))}.
\end{align*}

{\bf Part II.} We then consider the general case where $p$ is a fixed positive integer and Assumption~\ref{asp:proof1} holds for $p$. We only consider the case where $f_1(\phi(x))>0$. The case with $f_1(\phi(x))=0$ can be established in a similar way.

Let $\eta_N = \eta_{N,p} = 4p\log(N_0/M)$. We also take $N_0$ sufficiently large so that
\[
  \eta_N \frac{M}{N_0} = 4p\frac{M}{N_0} \log\Big(\frac{N_0}{M}\Big) < (1-\epsilon)f_0(\phi(x)) \lambda(B_{0,\delta/2}).
\]

Then
\begin{align*}
  &\P\Big(Z_1,\ldots,Z_p \in A_\phi(x)\Big) = \P\Big(\lVert \hat\phi(x) - \hat\phi(Z_k) \rVert \le \hat\Phi_M(Z_k), \forall k \in \zahl{p} \Big) \\
  &\le \P\Big(\lVert \phi(x) - \phi(Z_k) \rVert - 2\lVert \hat\phi-\phi \rVert_\infty \le \Phi_M(Z_k) + 2\lVert \hat\phi-\phi \rVert_\infty, \forall k \in \zahl{p} \Big)\\
  &= \sum_{S \subset \zahl{p}}\P\Big(\lVert \phi(x) - \phi(Z_k) \rVert - 4\lVert \hat\phi-\phi \rVert_\infty \le \Phi_M(Z_k), 4\lVert \hat\phi-\phi \rVert_\infty \le \epsilon' \lVert \phi(x) - \phi(Z_k) \rVert {\rm ~for~} k \in S,\\ &4\lVert \hat\phi-\phi \rVert_\infty > \epsilon' \lVert \phi(x) - \phi(Z_k) \rVert {\rm ~for~} k \notin S\Big)\\
  &\le \sum_{S \subset \zahl{p}}\P\Big(\nu_0(B_{\phi(Z_k),(1-\epsilon')\lVert \phi(x) - \phi(Z_k) \rVert }) \le \nu_0(B_{\phi(Z_k),\Phi_M(Z_k) }){\rm ~for~} k \in S, 4\lVert \hat\phi-\phi \rVert_\infty > \epsilon' \max_{k \notin S}\lVert \phi(x) - \phi(Z_k) \rVert\Big).
  \yestag\label{eq:catch1}
\end{align*}

Let $W_k = \nu_0(B_{\phi(Z_k),(1-\epsilon')\lVert \phi(x) - \phi(Z_k) \rVert })$ and $V_k = \nu_0(B_{\phi(Z_k),\Phi_M(Z_k) })$ for any $k \in \zahl{p}$. Then $[W_k]_{k=1}^p$ are i.i.d. since $[Z_k]_{k=1}^p$ are i.i.d.. For any $k \in \zahl{p}$ and $Z_k \in \cX$ given, $V_k \given Z_k $ has the same distribution as $U_{(M)}$. Then for any $k \in \zahl{p}$, $V_k$ has the same distribution as $U_{(M)}$, and $V_k$ is independent of $Z_k$.

Fix $S \subset \zahl{p}$. Let $W_{\max} = \max_{k \in S} W_k$ and $V_{\max} = \max_{k \in S} V_k$. Then
\begin{align*}
  & \P\Big(\nu_0(B_{\phi(Z_k),(1-\epsilon')\lVert \phi(x) - \phi(Z_k) \rVert }) \le \nu_0(B_{\phi(Z_k),\Phi_M(Z_k) }){\rm ~for~} k \in S, 4\lVert \hat\phi-\phi \rVert_\infty > \epsilon' \max_{k \notin S}\lVert \phi(x) - \phi(Z_k) \rVert\Big)\\
  &\le \P\Big(W_{\max} < V_{\max}, 4\lVert \hat\phi-\phi \rVert_\infty > \epsilon' \max_{k \notin S}\lVert \phi(x) - \phi(Z_k) \rVert\Big)\\
  &\le \P\Big(W_{\max} < V_{\max} \le \eta_N \frac{M}{N_0}, \max_{k \notin S} \lVert \phi(x) - \phi(Z_k) \rVert < 4\lVert \hat\phi-\phi \rVert_\infty/ \epsilon' \le \delta \Big) \\
  &+ \P\Big(V_{\max} > \eta_N \frac{M}{N_0}\Big) + \P\Big(4\lVert \hat\phi-\phi \rVert_\infty/ \epsilon' > \delta\Big).
  \yestag\label{eq:catch2}
\end{align*}

For the first term in \eqref{eq:catch2}, by the selection of $\eta_N$, and taking $\epsilon'<1/2$ and $N_0$ large enough, we have
\begin{align*}
  &\P\Big(W_{\max} < V_{\max} \le \eta_N \frac{M}{N_0}, \max_{k \notin S} \lVert \phi(x) - \phi(Z_k) \rVert < 4\lVert \hat\phi-\phi \rVert_\infty/ \epsilon' \le \delta \Big)\\
  &\le \P\Big(W_{\max} < V_{\max} \le \eta_N \frac{M}{N_0}, \max_{k \notin S} \lVert \phi(x) - \phi(Z_k) \rVert < 4\lVert \hat\phi-\phi \rVert_\infty/ \epsilon' \le \delta, \max_{k \in S} \lVert \phi(x) - \phi(Z_k) \rVert \le \delta \Big).
\end{align*}
Let $W_k' = \nu_0(B_{\phi(Z_k),\lVert \phi(x) - \phi(Z_k) \rVert })$ and $W_{\max}' = \max_{k \in S} W_k'$. Under the event $\{\max_{k \in S} \lVert \phi(x) - \phi(Z_k) \rVert \le \delta\}$, we have
\begin{align*}
  & W_{\max}' - W_{\max} \le \max_{k \in S} [\nu_0(B_{\phi(Z_k),\lVert \phi(x) - \phi(Z_k) \rVert }) - \nu_0(B_{\phi(Z_k),(1-\epsilon')\lVert \phi(x) - \phi(Z_k) \rVert })]\\
  &\le \frac{(1+\epsilon)f_0(\phi(x)) d \epsilon'}{(1-\epsilon)f_1(\phi(x))} \max_{k \in S} \nu_1(B_{\phi(x),\lVert \phi(x) - \phi(Z_k) \rVert }).
\end{align*}
and
\begin{align*}
  W_{\max}' \ge \frac{(1-\epsilon)f_0(\phi(x))}{(1+\epsilon)f_1(\phi(x))} \max_{k \in S} \nu_1(B_{\phi(x),\lVert \phi(x) - \phi(Z_k) \rVert }).
\end{align*}
On the other hand, recall that $\hat\phi_{(Z_1,\ldots,Z_p) \to z}$ is the estimator replacing $(Z_1,\ldots,Z_p)$ by $z$ for $z \in \cZ^p$. Then $\sup_{z \in \cZ^p} \lVert \hat\phi_{(Z_1,\ldots,Z_p) \to z}-\phi \rVert_\infty$ is independent of $(Z_1,\ldots,Z_p)$. Note that $$\max_{k \notin S} \lVert \phi(x) - \phi(Z_k) \rVert < 4\lVert \hat\phi-\phi \rVert_\infty/ \epsilon' \le \delta$$ implies that $$\max_{k \notin S} \nu_1(B_{\phi(x),\lVert \phi(x) - \phi(Z_k) \rVert }) \le (1+\epsilon) f_1(\phi(x)) V_m(4/\epsilon')^d \sup_{z \in \cZ^p} \lVert \hat\phi_{(Z_1,\ldots,Z_p) \to z}-\phi \rVert_\infty^d.$$
Then
\begin{align*}
  &\P\Big(W_{\max} < V_{\max} \le \eta_N \frac{M}{N_0}, \max_{k \notin S} \lVert \phi(x) - \phi(Z_k) \rVert < 4\lVert \hat\phi-\phi \rVert_\infty/ \epsilon' \le \delta, \max_{k \in S} \lVert \phi(x) - \phi(Z_k) \rVert \le \delta \Big)\\
  &\le \P\Big(\Big(\frac{1-\epsilon}{1+\epsilon} - \frac{1+\epsilon}{1-\epsilon} d \epsilon' \Big) \frac{f_0(\phi(x))}{f_1(\phi(x))}\max_{k \in S}\nu_1(B_{\phi(x),\lVert \phi(x) - \phi(Z_k) \rVert }) < V_{\max}, \\
  & \max_{k \notin S} \nu_1(B_{\phi(x),\lVert \phi(x) - \phi(Z_k) \rVert }) \le (1+\epsilon) f_1(\phi(x)) V_m(4/\epsilon')^d \sup_{z \in \cZ^p} \lVert \hat\phi_{(Z_1,\ldots,Z_p) \to z}-\phi \rVert_\infty^d\Big)\\
  &= \E\Big[\ind\Big(\Big(\frac{1-\epsilon}{1+\epsilon} - \frac{1+\epsilon}{1-\epsilon} d \epsilon' \Big) \frac{f_0(\phi(x))}{f_1(\phi(x))}\max_{k \in S}\nu_1(B_{\phi(x),\lVert \phi(x) - \phi(Z_k) \rVert }) < V_{\max}\Big)\\
  &\Big([(1+\epsilon) f_1(\phi(x)) V_m(4/\epsilon')^d \sup_{z \in \cZ^p} \lVert \hat\phi_{(Z_1,\ldots,Z_p) \to z}-\phi \rVert_\infty^d]\wedge 1\Big)^{p-\lvert S \rvert}\Big],
\end{align*}
since $\sup_{z \in \cZ^p} \lVert \hat\phi_{(Z_1,\ldots,Z_p) \to z}-\phi \rVert_\infty$, $V_{\max}$ and $[Z_k]_{k \in S}$ are all independent with $[Z_k]_{k \notin S}$.

Note that
\begin{align*}
  & \E\Big[\ind\Big(\Big(\frac{1-\epsilon}{1+\epsilon} - \frac{1+\epsilon}{1-\epsilon} d \epsilon' \Big) \frac{f_0(\phi(x))}{f_1(\phi(x))}\max_{k \in S}\nu_1(B_{\phi(x),\lVert \phi(x) - \phi(Z_k) \rVert }) < V_{\max}\Big)\\
  &\Big([(1+\epsilon) f_1(\phi(x)) V_m(4/\epsilon')^d \sup_{z \in \cZ^p} \lVert \hat\phi_{(Z_1,\ldots,Z_p) \to z}-\phi \rVert_\infty^d]\wedge 1\Big)^{p-\lvert S \rvert}\Big]\\
  &=\int_0^1 \lvert S \rvert u^{\lvert S \rvert-1} \E\Big[\ind\Big(\Big(\frac{1-\epsilon}{1+\epsilon} - \frac{1+\epsilon}{1-\epsilon} d \epsilon' \Big) \frac{f_0(\phi(x))}{f_1(\phi(x))}u < V_{\max}\Big)\\
  &\Big([(1+\epsilon) f_1(\phi(x)) V_m(4/\epsilon')^d \sup_{z \in \cZ^p} \lVert \hat\phi_{(Z_1,\ldots,Z_p) \to z}-\phi \rVert_\infty^d]\wedge 1\Big)^{p-\lvert S \rvert} \Biggiven \max_{k \in S}\nu_1(B_{\phi(x),\lVert \phi(x) - \phi(Z_k) \rVert }) = u\Big] \d u\\
  &= \lvert S \rvert \Big[\Big(\frac{1-\epsilon}{1+\epsilon} - \frac{1+\epsilon}{1-\epsilon} d \epsilon' \Big)^{-1} \frac{f_1(\phi(x))}{f_0(\phi(x))} \frac{M}{N_0}\Big]^{\lvert S \rvert} \int_0^{(\frac{1-\epsilon}{1+\epsilon} - \frac{1+\epsilon}{1-\epsilon} d \epsilon') \frac{f_0(\phi(x))}{f_1(\phi(x))} \frac{N_0}{M}} u^{\lvert S \rvert-1} \E\Big[\ind\Big(V_{\max} > \frac{M}{N_0} u \Big)\\
  &\Big([(1+\epsilon) f_1(\phi(x)) V_m(4/\epsilon')^d \sup_{z \in \cZ^p} \lVert \hat\phi_{(Z_1,\ldots,Z_p) \to z}-\phi \rVert_\infty^d]\wedge 1\Big)^{p-\lvert S \rvert} \Biggiven\\ 
  & \max_{k \in S}\nu_1(B_{\phi(x),\lVert \phi(x) - \phi(Z_k) \rVert }) = \Big(\frac{1-\epsilon}{1+\epsilon} - \frac{1+\epsilon}{1-\epsilon} d \epsilon' \Big)^{-1} \frac{f_1(\phi(x))}{f_0(\phi(x))} \frac{M}{N_0} u\Big] \d u.
  \yestag\label{eq:catch3}
\end{align*}

We split the above integral into two parts using 1. For the first part, note that $\sup_{z \in \cZ^p} \lVert \hat\phi_{(Z_1,\ldots,Z_p) \to z}-\phi \rVert_\infty$ is independent with $\max_{k \in S}\nu_1(B_{\phi(x),\lVert \phi(x) - \phi(Z_k) \rVert })$. Then
\begin{align*}
  &\Big(\frac{N_0}{M}\Big)^{p-\lvert S \rvert}\int_0^1 u^{\lvert S \rvert-1} \E\Big[\ind\Big(V_{\max} > \frac{M}{N_0} u \Big) \Big([(1+\epsilon) f_1(\phi(x)) V_m(4/\epsilon')^d \sup_{z \in \cZ^p} \lVert \hat\phi_{(Z_1,\ldots,Z_p) \to z}-\phi \rVert_\infty^d]\wedge 1\Big)^{p-\lvert S \rvert} \Biggiven\\ 
  & \max_{k \in S}\nu_1(B_{\phi(x),\lVert \phi(x) - \phi(Z_k) \rVert }) = \Big(\frac{1-\epsilon}{1+\epsilon} - \frac{1+\epsilon}{1-\epsilon} d \epsilon' \Big)^{-1} \frac{f_1(\phi(x))}{f_0(\phi(x))} \frac{M}{N_0} u\Big] \d u \\
  &\le \Big(\frac{N_0}{M}\Big)^{p-\lvert S \rvert}\int_0^1 u^{\lvert S \rvert-1} \E\Big[ \Big([(1+\epsilon) f_1(\phi(x)) V_m(4/\epsilon')^d \sup_{z \in \cZ^p} \lVert \hat\phi_{(Z_1,\ldots,Z_p) \to z}-\phi \rVert_\infty^d]\wedge 1\Big)^{p-\lvert S \rvert}\Big] \d u. 
\end{align*}

If $\lvert S \rvert = p$, we have
\begin{align*}
  &\Big(\frac{N_0}{M}\Big)^{p-\lvert S \rvert}\int_0^1 u^{\lvert S \rvert-1} \E\Big[ \Big([(1+\epsilon) f_1(\phi(x)) V_m(4/\epsilon')^d \sup_{z \in \cZ^p} \lVert \hat\phi_{(Z_1,\ldots,Z_p) \to z}-\phi \rVert_\infty^d]\wedge 1\Big)^{p-\lvert S \rvert}\Big] \d u \\
  &= \int_0^1 u^{p-1} \d u = \frac{1}{p}.
  \yestag\label{eq:catch4}
\end{align*}

If $\lvert S \rvert < p$, by Assumption~\ref{asp:proof1}, we have
\begin{align*}
  &\limsup_{N_0 \to \infty}\Big(\frac{N_0}{M}\Big)^{p-\lvert S \rvert}\int_0^1 u^{\lvert S \rvert-1} \E\Big[ \Big([(1+\epsilon) f_1(\phi(x)) V_m(4/\epsilon')^d \sup_{z \in \cZ^p} \lVert \hat\phi_{(Z_1,\ldots,Z_p) \to z}-\phi \rVert_\infty^d]\wedge 1\Big)^{p-\lvert S \rvert}\Big] \d u \\
  &\lesssim \limsup_{N_0 \to \infty} \E\Big[ \Big(\frac{N_0}{M}\sup_{z \in \cZ^p} \lVert \hat\phi_{(Z_1,\ldots,Z_p) \to z}-\phi \rVert_\infty^d\Big)^{p-\lvert S \rvert}\Big] = 0.
  \yestag\label{eq:catch5}
\end{align*}

For the second part, we have
\begin{align*}
  & \int_1^{(\frac{1-\epsilon}{1+\epsilon} - \frac{1+\epsilon}{1-\epsilon} d \epsilon') \frac{f_0(\phi(x))}{f_1(\phi(x))} \frac{N_0}{M}} u^{\lvert S \rvert-1} \E\Big[\ind\Big(V_{\max} > \frac{M}{N_0} u \Big) \Big([(1+\epsilon) f_1(\phi(x)) V_m(4/\epsilon')^d \sup_{z \in \cZ^p} \lVert \hat\phi_{(Z_1,\ldots,Z_p) \to z}-\phi \rVert_\infty^d]\wedge 1\Big)^{p-\lvert S \rvert} \Biggiven\\ 
  & \max_{k \in S}\nu_1(B_{\phi(x),\lVert \phi(x) - \phi(Z_k) \rVert }) = \Big(\frac{1-\epsilon}{1+\epsilon} - \frac{1+\epsilon}{1-\epsilon} d \epsilon' \Big)^{-1} \frac{f_1(\phi(x))}{f_0(\phi(x))} \frac{M}{N_0} u\Big] \d u\\
  &\le \int_1^\infty u^{\lvert S \rvert-1} \E\Big[\ind\Big(V_{\max} > \frac{M}{N_0} u \Big) \Big([(1+\epsilon) f_1(\phi(x)) V_m(4/\epsilon')^d \sup_{z \in \cZ^p} \lVert \hat\phi_{(Z_1,\ldots,Z_p) \to z}-\phi \rVert_\infty^d]\wedge 1\Big)^{p-\lvert S \rvert} \Biggiven\\ 
  & \max_{k \in S}\nu_1(B_{\phi(x),\lVert \phi(x) - \phi(Z_k) \rVert }) = \Big(\frac{1-\epsilon}{1+\epsilon} - \frac{1+\epsilon}{1-\epsilon} d \epsilon' \Big)^{-1} \frac{f_1(\phi(x))}{f_0(\phi(x))} \frac{M}{N_0} u\Big] \d u\\
  &\le \sum_{k \in S} \int_1^\infty u^{\lvert S \rvert-1} \E\Big[\ind\Big(V_k > \frac{M}{N_0} u \Big) \Big([(1+\epsilon) f_1(\phi(x)) V_m(4/\epsilon')^d \sup_{z \in \cZ^p} \lVert \hat\phi_{(Z_1,\ldots,Z_p) \to z}-\phi \rVert_\infty^d]\wedge 1\Big)^{p-\lvert S \rvert} \Big] \d u,
  \yestag\label{eq:catch6}
\end{align*}
where the last step is from the fact that $V_k$ and $\sup_{z \in \cZ^p} \lVert \hat\phi_{(Z_1,\ldots,Z_p) \to z}-\phi \rVert_\infty$ are independent of $[Z_k]_{k \in S}$ for any $k \in S$.

For any $k \in S$, by the Hölder inequality,
\begin{align*}
  & \Big(\frac{N_0}{M}\Big)^{p-\lvert S \rvert} \int_1^\infty u^{\lvert S \rvert-1} \E\Big[\ind\Big(V_k > \frac{M}{N_0} u \Big) \Big([(1+\epsilon) f_1(\phi(x)) V_m(4/\epsilon')^d \sup_{z \in \cZ^p} \lVert \hat\phi_{(Z_1,\ldots,Z_p) \to z}-\phi \rVert_\infty^d]\wedge 1\Big)^{p-\lvert S \rvert} \Big] \d u\\
  &\lesssim \int_1^\infty u^{\lvert S \rvert-1} \E\Big[\ind\Big(V_k > \frac{M}{N_0} u \Big) \Big(\frac{N_0}{M} \sup_{z \in \cZ^p} \lVert \hat\phi_{(Z_1,\ldots,Z_p) \to z}-\phi \rVert_\infty^d\Big)^{p-\lvert S \rvert} \Big] \d u\\
  &\le \int_1^\infty u^{\lvert S \rvert-1} \Big\{\P\Big(V_k > \frac{M}{N_0} u \Big)\Big\}^{\frac{\lvert S \rvert}{p}}\Big\{\E\Big[\Big(\frac{N_0}{M} \sup_{z \in \cZ^p} \lVert \hat\phi_{(Z_1,\ldots,Z_p) \to z}-\phi \rVert_\infty^d\Big)^p\Big]\Big\}^{\frac{p-\lvert S \rvert}{p}} \d u\\
  &= \Big\{\E\Big[\Big(\frac{N_0}{M} \sup_{z \in \cZ^p} \lVert \hat\phi_{(Z_1,\ldots,Z_p) \to z}-\phi \rVert_\infty^d\Big)^p\Big]\Big\}^{\frac{p-\lvert S \rvert}{p}} \int_1^\infty u^{\lvert S \rvert-1} \Big\{\P\Big(V_k > \frac{M}{N_0} u \Big)\Big\}^{\frac{\lvert S \rvert}{p}} \d u.
\end{align*}

Using the Chernoff bound,
\begin{align*}
  & \int_1^\infty u^{\lvert S \rvert-1} \Big\{\P\Big(V_k > \frac{M}{N_0} u \Big)\Big\}^{\frac{\lvert S \rvert}{p}} \d u = \int_0^\infty (1+u)^{\lvert S \rvert-1} \Big\{\P\Big(U_{(M)} > \frac{M}{N_0} (1+u) \Big)\Big\}^{\frac{\lvert S \rvert}{p}} \d u\\
  &\le \int_0^\infty (1+u)^{\lvert S \rvert-1} (1+u)^{M\lvert S \rvert/p} \exp(-uM\lvert S \rvert/p) \d u \\
  &= \exp(M\lvert S \rvert/p)\int_1^\infty u^{M\lvert S \rvert/p+\lvert S \rvert-1} \exp(-uM\lvert S \rvert/p) \d u \\
  &\le \exp(M\lvert S \rvert/p)\int_0^\infty u^{M\lvert S \rvert/p+\lvert S \rvert-1} \exp(-uM\lvert S \rvert/p) \d u \\
  &= \frac{\exp(M\lvert S \rvert/p)}{(M\lvert S \rvert/p)^{M\lvert S \rvert/p+\lvert S \rvert}} \Gamma(M\lvert S \rvert/p+\lvert S \rvert) \\
  &= \frac{\exp(M\lvert S \rvert/p)}{(M\lvert S \rvert/p)^{M\lvert S \rvert/p+\lvert S \rvert}} (M\lvert S \rvert/p+1)^{\lvert S \rvert-1} \Gamma(M\lvert S \rvert/p+1) (1+o(1)) \\
  &= \frac{\exp(M\lvert S \rvert/p)}{(M\lvert S \rvert/p)^{M\lvert S \rvert/p+\lvert S \rvert}} (M\lvert S \rvert/p+1)^{\lvert S \rvert-1} \sqrt{2\pi M\lvert S \rvert/p} \Big(\frac{M\lvert S \rvert/p}{e}\Big)^{M\lvert S \rvert/p} (1+o(1)) \\
  &= \sqrt{2\pi} (M\lvert S \rvert/p)^{-1/2} \Big(1+\frac{p}{M\lvert S \rvert}\Big)^{\lvert S \rvert-1} (1+o(1)),
\end{align*}
where the last three steps are from Stirling's approximation using $M \to \infty$.

By $M \to \infty$ and Assumption~\ref{asp:proof1}, we have
\begin{align*}
  &\lim_{N_0 \to \infty} \Big(\frac{N_0}{M}\Big)^{p-\lvert S \rvert} \int_1^\infty u^{\lvert S \rvert-1} \E\Big[\ind\Big(V_k > \frac{M}{N_0} u \Big) \Big([(1+\epsilon) f_1(\phi(x)) V_m(4/\epsilon')^d \sup_{z \in \cZ^p} \lVert \hat\phi_{(Z_1,\ldots,Z_p) \to z}-\phi \rVert_\infty^d]\wedge 1\Big)^{p-\lvert S \rvert} \Big] \d u \\
  &= 0.
  \yestag\label{eq:catch7}
\end{align*}

Combining \eqref{eq:catch4}, \eqref{eq:catch5}, \eqref{eq:catch6}, \eqref{eq:catch7} by \eqref{eq:catch3} yields
\begin{align*}
  &\limsup_{N_0 \to \infty} \Big(\frac{N_0}{M}\Big)^p \P\Big(W_{\max} < V_{\max} \le \eta_N \frac{M}{N_0}, \max_{k \notin S} \lVert \phi(x) - \phi(Z_k) \rVert < 4\lVert \hat\phi-\phi \rVert_\infty/ \epsilon' \le \delta \Big) \\
  &\le  \frac{1}{p} \lvert S \rvert \Big[\Big(\frac{1-\epsilon}{1+\epsilon} - \frac{1+\epsilon}{1-\epsilon} d \epsilon' \Big)^{-1} \frac{f_1(\phi(x))}{f_0(\phi(x))} \Big]^{\lvert S \rvert} \ind(\lvert S \rvert = p) = \Big[\Big(\frac{1-\epsilon}{1+\epsilon} - \frac{1+\epsilon}{1-\epsilon} d \epsilon' \Big)^{-1} \frac{f_1(\phi(x))}{f_0(\phi(x))} \Big]^{p} \ind(\lvert S \rvert = p).
\end{align*}

For the second term in \eqref{eq:catch2}, by the Chernoff bound,
\begin{align*}
  \limsup_{N_0 \to \infty} \Big(\frac{N_0}{M}\Big)^p \P\Big(V_{\max} > \eta_N \frac{M}{N_0}\Big) \le \lvert S \rvert \limsup_{N_0 \to \infty} \Big(\frac{N_0}{M}\Big)^p \P\Big(U_{(M)} > \eta_N \frac{M}{N_0}\Big) = 0.
\end{align*}

For the third term in \eqref{eq:catch2},
\begin{align*}
  &\P\Big(4\lVert \hat\phi-\phi \rVert_\infty/ \epsilon' > \delta\Big) \le \delta^{-pd} (4/\epsilon')^{pd} \E[\lVert \hat\phi-\phi \rVert_\infty^{pd}] \\
  &\le \delta^{-pd} (4/\epsilon')^{pd} \E[\sup_{z \in \cZ^p} \lVert \hat\phi_{(Z_1,\ldots,Z_p) \to z}-\phi \rVert_\infty^{pd}].
\end{align*}

By Assumption~\ref{asp:proof1},
\begin{align*}
  \limsup_{N_0 \to \infty} \Big(\frac{N_0}{M}\Big)^p \P\Big(4\lVert \hat\phi-\phi \rVert_\infty/ \epsilon' > \delta\Big) = 0.
\end{align*}

Then we obtain for any $S \subset \zahl{p}$,
\begin{align*}
  &\limsup_{N_0 \to \infty} \Big(\frac{N_0}{M}\Big)^p \P\Big(\nu_0(B_{\phi(Z_k),(1-\epsilon')\lVert \phi(x) - \phi(Z_k) \rVert }) \le \nu_0(B_{\phi(Z_k),\Phi_M(Z_k) }){\rm ~for~} k \in S,\\
  & 4\lVert \hat\phi-\phi \rVert_\infty > \epsilon' \max_{k \notin S}\lVert \phi(x) - \phi(Z_k) \rVert\Big)\\
  &\le \Big[\Big(\frac{1-\epsilon}{1+\epsilon} - \frac{1+\epsilon}{1-\epsilon} d \epsilon' \Big)^{-1} \frac{f_1(\phi(x))}{f_0(\phi(x))}\Big]^p \ind(\lvert S \rvert = p),
\end{align*}
and then by \eqref{eq:catch1},
\begin{align*}
  \limsup_{N_0 \to \infty} \Big(\frac{N_0}{M}\Big)^p \P\Big(Z_1,\ldots,Z_p \in A_\phi(x)\Big) \le \Big[\Big(\frac{1-\epsilon}{1+\epsilon} - \frac{1+\epsilon}{1-\epsilon} d \epsilon' \Big)^{-1} \frac{f_1(\phi(x))}{f_0(\phi(x))}\Big]^p.
\end{align*}

By $\epsilon,\epsilon'$ are arbitrary, we obtain
\begin{align*}
  \limsup_{N_0 \to \infty} \Big(\frac{N_0}{M}\Big)^p \P\Big(Z_1,\ldots,Z_p \in A_\phi(x)\Big) \le \Big(\frac{f_1(\phi(x))}{f_0(\phi(x))}\Big)^p.
\end{align*}

A matched lower bound is directly from the H\"older inequality.

Then we obtain
\begin{align*}
  \lim_{N_0 \to \infty} \Big(\frac{N_0}{M}\Big)^p \P\Big(Z_1,\ldots,Z_p \in A_\phi(x)\Big) = \Big(\frac{f_1(\phi(x))}{f_0(\phi(x))}\Big)^p.
\end{align*}

{\bf Part III.} We consider the simple case where $p=1$ and Assumption~\ref{asp:proof2} holds. We only consider the case where $f_1(\phi(x))>0$, while the case where $f_1(\phi(x))=0$ is similar.

{\bf Upper bound.} Note that
\begin{align*}
  &\P\Big(Z_1 \in A_\phi(x)\Big) = \P\Big(\lVert \hat\phi(x) - \hat\phi(Z_1) \rVert \le \hat\Phi_M(Z_1)\Big) \\
  &\le \P\Big(\lVert \phi(x) - \phi(Z_1) \rVert - 2\sup_{\lVert \phi(s) - \phi(t)\rVert \le \lVert \phi(x) - \phi(Z_1)\rVert}\lVert (\hat\phi-\phi)(s) - (\hat\phi-\phi)(t) \rVert \le \Phi_M(Z_1)\Big).
\end{align*}

The inequality is from the fact that under the event $\{\lVert \hat\phi(x) - \hat\phi(Z_1) \rVert \le \hat\Phi_M(Z_1)\}$, if $\lVert \phi(x) - \phi(Z_1) \rVert - 2\sup_{\lVert \phi(s) - \phi(t)\rVert \le \lVert \phi(x) - \phi(Z_1)\rVert}\lVert (\hat\phi-\phi)(s) - (\hat\phi-\phi)(t) \rVert > \Phi_M(Z_1)$, then there exists a set $S \subset \zahl{N_0}$ such that $\lvert S \rvert \ge M$ and for any $i \in S$, $\lVert \phi(x) - \phi(Z_1) \rVert - 2\sup_{\lVert \phi(s) - \phi(t)\rVert \le \lVert \phi(x) - \phi(Z_1)\rVert}\lVert (\hat\phi-\phi)(s) - (\hat\phi-\phi)(t) \rVert > \lVert \phi(X_i) - \phi(Z_1) \rVert$. For these $i \in S$, we then have $\lVert \hat\phi(x) - \hat\phi(Z_1) \rVert > \lVert \hat\phi(X_i) - \hat\phi(Z_1) \rVert$ since $\lVert \phi(x) - \phi(Z_1) \rVert > \lVert \phi(X_i) - \phi(Z_1) \rVert$. Then $\lVert \hat\phi(x) - \hat\phi(Z_1) \rVert > \hat\Phi_M(Z_1)$ using $\lvert S \rvert \ge M$, which contradicts the event we assume.

For any $\epsilon \in (0,1)$, we decompose as
\begin{align*}
  & \P\Big(\lVert \phi(x) - \phi(Z_1) \rVert - 2\sup_{\lVert \phi(s) - \phi(t)\rVert \le \lVert \phi(x) - \phi(Z_1)\rVert}\lVert (\hat\phi-\phi)(s) - (\hat\phi-\phi)(t) \rVert \le \Phi_M(Z_1)\Big)\\
  &= \P\Big(\lVert \phi(x) - \phi(Z_1) \rVert - 2\sup_{\lVert \phi(s) - \phi(t)\rVert \le \lVert \phi(x) - \phi(Z_1)\rVert}\lVert (\hat\phi-\phi)(s) - (\hat\phi-\phi)(t) \rVert \le \Phi_M(Z_1),\\
  &2\sup_{\lVert \phi(s) - \phi(t)\rVert \le \lVert \phi(x) - \phi(Z_1)\rVert}\lVert (\hat\phi-\phi)(s) - (\hat\phi-\phi)(t) \rVert \le \epsilon \lVert \phi(x) - \phi(Z_1) \rVert\Big)\\
  &+ \P\Big(\lVert \phi(x) - \phi(Z_1) \rVert - 2\sup_{\lVert \phi(s) - \phi(t)\rVert \le \lVert \phi(x) - \phi(Z_1)\rVert}\lVert (\hat\phi-\phi)(s) - (\hat\phi-\phi)(t) \rVert \le \Phi_M(Z_1),\\
  &2\sup_{\lVert \phi(s) - \phi(t)\rVert \le \lVert \phi(x) - \phi(Z_1)\rVert}\lVert (\hat\phi-\phi)(s) - (\hat\phi-\phi)(t) \rVert >  \epsilon \lVert \phi(x) - \phi(Z_1) \rVert\Big).
  \yestag\label{eq:catch11}
\end{align*}

For the first term in \eqref{eq:catch11},
\begin{align*}
  &\P\Big(\lVert \phi(x) - \phi(Z_1) \rVert - 2\sup_{\lVert \phi(s) - \phi(t)\rVert \le \lVert \phi(x) - \phi(Z_1)\rVert}\lVert (\hat\phi-\phi)(s) - (\hat\phi-\phi)(t) \rVert \le \Phi_M(Z_1),\\
  &2\sup_{\lVert \phi(s) - \phi(t)\rVert \le \lVert \phi(x) - \phi(Z_1)\rVert}\lVert (\hat\phi-\phi)(s) - (\hat\phi-\phi)(t) \rVert \le \epsilon \lVert \phi(x) - \phi(Z_1) \rVert\Big)\\
  &\le \P\Big((1-\epsilon)\lVert \phi(x) - \phi(Z_1) \rVert \le \Phi_M(Z_1)\Big),
\end{align*}
and then can be handled in the same way as the upper bound part in Part I. 

For the second term in \eqref{eq:catch11}, note that for any $u \in (0,1)$, conditional on $\nu_1(B_{\phi(x),\lVert \phi(x) - \phi(Z_1) \rVert }) = u$ and under $\lVert \phi(x) - \phi(Z_1) \rVert \le \delta$, we have $u \le (1+\epsilon)f_1(\phi(x)) V_m \lVert \phi(x) - \phi(Z_1) \rVert^m$, and then $\lVert \phi(x) - \phi(Z_1) \rVert \ge [u/((1+\epsilon)f_1(\phi(x)) V_m)]^{1/m}$. Then for any $u \in (0,1)$,
\begin{align*}
  &\P\Big(\sup_{\lVert \phi(s) - \phi(t)\rVert \le \lVert \phi(x) - \phi(Z_1)\rVert}\lVert (\hat\phi-\phi)(s) - (\hat\phi-\phi)(t) \rVert >  \epsilon \lVert \phi(x) - \phi(Z_1) \rVert, \lVert \phi(x) - \phi(Z_1) \rVert \le \delta\\
  &\Biggiven \nu_1(B_{\phi(x),\lVert \phi(x) - \phi(Z_1) \rVert }) = u \Big)\\
  &\le \P\Big(\sup_{\lVert \phi(s) - \phi(t)\rVert \le \lVert \phi(x) - \phi(Z_1)\rVert}\sup_{z \in \cZ}\lVert (\hat\phi_{Z_1 \to z}-\phi)(s) - (\hat\phi_{Z_1 \to z}-\phi)(t) \rVert >  \epsilon \lVert \phi(x) - \phi(Z_1) \rVert, \lVert \phi(x) - \phi(Z_1) \rVert \le \delta \\
  &\Biggiven \nu_1(B_{\phi(x),\lVert \phi(x) - \phi(Z_1) \rVert }) = u \Big)\\
  &\le T_{\epsilon}([u/((1+\epsilon)f_1(\phi(x)) V_m)]^{1/m}),
\end{align*}
where the last step is from the fact that $\sup_{z \in \cZ}\lVert (\hat\phi_{Z_1 \to z}-\phi)(s) - (\hat\phi_{Z_1 \to z}-\phi)(t) \rVert$ does not depend on $Z_1$ together with Assumption~\ref{asp:proof2}.

Then
\begin{align*}
  &\P\Big(\lVert \phi(x) - \phi(Z_1) \rVert - 2\sup_{\lVert \phi(s) - \phi(t)\rVert \le \lVert \phi(x) - \phi(Z_1)\rVert}\lVert (\hat\phi-\phi)(s) - (\hat\phi-\phi)(t) \rVert \le \Phi_M(Z_1),\\
  &2\sup_{\lVert \phi(s) - \phi(t)\rVert \le \lVert \phi(x) - \phi(Z_1)\rVert}\lVert (\hat\phi-\phi)(s) - (\hat\phi-\phi)(t) \rVert >  \epsilon \lVert \phi(x) - \phi(Z_1) \rVert\Big)\\
  &\le\P\Big(2\sup_{\lVert \phi(s) - \phi(t)\rVert \le \lVert \phi(x) - \phi(Z_1)\rVert}\lVert (\hat\phi-\phi)(s) - (\hat\phi-\phi)(t) \rVert >  \epsilon \lVert \phi(x) - \phi(Z_1) \rVert\Big)\\
  &\le\P\Big(2\sup_{\lVert \phi(s) - \phi(t)\rVert \le \lVert \phi(x) - \phi(Z_1)\rVert}\lVert (\hat\phi-\phi)(s) - (\hat\phi-\phi)(t) \rVert >  \epsilon \lVert \phi(x) - \phi(Z_1) \rVert, \lVert \phi(x) - \phi(Z_1) \rVert \le \delta\Big)\\
  &+\P\Big(2\sup_{\lVert \phi(s) - \phi(t)\rVert \le \lVert \phi(x) - \phi(Z_1)\rVert}\lVert (\hat\phi-\phi)(s) - (\hat\phi-\phi)(t) \rVert >  \epsilon \lVert \phi(x) - \phi(Z_1) \rVert, \lVert \phi(x) - \phi(Z_1) \rVert > \delta\Big)\\
  &\le \int_0^1 T_{\epsilon/2}([u/((1+\epsilon)f_1(\phi(x)) V_m)]^{1/m}) \d u + \P\Big(\lVert \hat\phi-\phi\rVert_\infty >  \epsilon \delta/4\Big)\\
  &= (1+\epsilon)f_1(\phi(x)) V_m \int_0^{1/((1+\epsilon)f_1(\phi(x)) V_m)} T_{\epsilon/2}(u^{1/m}) \d u + \P\Big(\lVert \hat\phi-\phi\rVert_\infty >  \epsilon \delta/4\Big).
\end{align*}

By Assumption~\ref{asp:proof2} and $\epsilon$ is arbitrary, using \eqref{eq:catch11}, we obtain
\begin{align*}
  \limsup_{N_0 \to \infty} \frac{N_0}{M} \P\Big(Z_1 \in A_\phi(x)\Big) \le \frac{f_1(\phi(x))}{f_0(\phi(x))}.
\end{align*}

{\bf Lower bound.} For any $\epsilon \in (0,1)$,
\begin{align*}
  &\P\Big(Z_1 \in A_\phi(x)\Big) = \P\Big(\lVert \hat\phi(x) - \hat\phi(Z_1) \rVert \le \hat\Phi_M(Z_1)\Big) \\
  &\ge \P\Big(\lVert \hat\phi(x) - \hat\phi(Z_1) \rVert \le \hat\Phi_M(Z_1), \sup_{\delta \ge \lVert \phi(x) - \phi(Z_1)\rVert} \delta^{-1} \sup_{\lVert \phi(s) - \phi(t)\rVert \le \delta}\lVert (\hat\phi-\phi)(s) - (\hat\phi-\phi)(t) \rVert \le \epsilon \Big) \\
  &\ge \P\Big((1+\epsilon)\lVert \phi(x) - \phi(Z_1) \rVert \le (1-\epsilon)\Phi_M(Z_1), \sup_{\delta \ge \lVert \phi(x) - \phi(Z_1)\rVert} \delta^{-1} \sup_{\lVert \phi(s) - \phi(t)\rVert \le \delta}\lVert (\hat\phi-\phi)(s) - (\hat\phi-\phi)(t) \rVert \le \epsilon\Big).
\end{align*}
The last inequality is from the fact that under the event $\{(1+\epsilon)\lVert \phi(x) - \phi(Z_1) \rVert \le (1-\epsilon)\Phi_M(Z_1)\}$, there exists a set $S \subset \zahl{N_0}$ such that $\lvert S \rvert \ge N_0-M$ and for any $i \in S$, $(1+\epsilon)\lVert \phi(x) - \phi(Z_1) \rVert \le (1-\epsilon)\lVert \phi(X_i) - \phi(Z_1) \rVert$. Under the event that $\{\sup_{\delta \ge \lVert \phi(x) - \phi(Z_1)\rVert} \delta^{-1} \sup_{\lVert \phi(s) - \phi(t)\rVert \le \delta}\lVert (\hat\phi-\phi)(s) - (\hat\phi-\phi)(t) \rVert \le \epsilon\}$, for these $i \in S$, we then have $\lVert \hat\phi(x) - \hat\phi(Z_1) \rVert \le \lVert \hat\phi(X_i) - \hat\phi(Z_1) \rVert$ since $\lVert \phi(x) - \phi(Z_1) \rVert \le \lVert \phi(X_i) - \phi(Z_1) \rVert$ for $i \in S$. Then $\lVert \hat\phi(x) - \hat\phi(Z_1) \rVert \le \hat\Phi_M(Z_1)$.

Then
\begin{align*}
  \P\Big(Z_1 \in A_\phi(x)\Big)
  &\ge \P\Big((1+\epsilon)\lVert \phi(x) - \phi(Z_1) \rVert \le (1-\epsilon)\Phi_M(Z_1)\Big)\\
  & \qquad - \P\Big(\sup_{\delta \ge \lVert \phi(x) - \phi(Z_1)\rVert} \delta^{-1} \sup_{\lVert \phi(s) - \phi(t)\rVert \le \delta}\lVert (\hat\phi-\phi)(s) - (\hat\phi-\phi)(t) \rVert > \epsilon\Big).
\end{align*}

The first term can be handled in the same way as the lower bound part in Part I. The second term can be handled in the same way as the second term in \eqref{eq:catch11} in this proof. Then we can obtain a matched lower bound.

{\bf Part IV.} We then consider the general case where $p$ is a fixed positive integer and Assumption~\ref{asp:proof2} holds. We only consider the case where $f_1(\phi(x))>0$, while the case where $f_1(\phi(x))=0$ is similar.

For any $\epsilon' \in (0,1)$, we have
\begin{align*}
  &\P\Big(Z_1,\ldots,Z_p \in A_\phi(x)\Big) = \P\Big(\lVert \hat\phi(x) - \hat\phi(Z_k) \rVert \le \hat\Phi_M(Z_k), \forall k \in \zahl{p} \Big) \\
  &\le \P\Big(\lVert \phi(x) - \phi(Z_k) \rVert - 2\sup_{\lVert \phi(s) - \phi(t)\rVert \le \lVert \phi(x) - \phi(Z_k)\rVert}\lVert (\hat\phi-\phi)(s) - (\hat\phi-\phi)(t) \rVert \le \Phi_M(Z_k), \forall k \in \zahl{p} \Big)\\
  &= \sum_{S \subset \zahl{p}}\P\Big(\lVert \phi(x) - \phi(Z_k) \rVert - 2\sup_{\lVert \phi(s) - \phi(t)\rVert \le \lVert \phi(x) - \phi(Z_k)\rVert}\lVert (\hat\phi-\phi)(s) - (\hat\phi-\phi)(t) \rVert \le \Phi_M(Z_k),\\
  &2\sup_{\lVert \phi(s) - \phi(t)\rVert \le \lVert \phi(x) - \phi(Z_k)\rVert}\lVert (\hat\phi-\phi)(s) - (\hat\phi-\phi)(t) \rVert \le \epsilon' \lVert \phi(x) - \phi(Z_k) \rVert {\rm ~for~} k \in S,\\
  &2\sup_{\lVert \phi(s) - \phi(t)\rVert \le \lVert \phi(x) - \phi(Z_k)\rVert}\lVert (\hat\phi-\phi)(s) - (\hat\phi-\phi)(t) \rVert > \epsilon' \lVert \phi(x) - \phi(Z_k) \rVert {\rm ~for~} k \notin S\Big)\\
  &\le \sum_{S \subset \zahl{p}}\P\Big((1-\epsilon')\lVert \phi(x) - \phi(Z_k) \rVert \le \Phi_M(Z_k){\rm ~for~} k \in S,\\
  &2\lVert \phi(x) - \phi(Z_k) \rVert^{-1}\sup_{\lVert \phi(s) - \phi(t)\rVert \le \lVert \phi(x) - \phi(Z_k)\rVert}\lVert (\hat\phi-\phi)(s) - (\hat\phi-\phi)(t) \rVert > \epsilon' {\rm ~for~} k \notin S\Big).
\end{align*}

If $\lvert S \rvert =p$, we have in the same way as the upper bound part in Part III that for any $\epsilon \in (0,1)$,
\begin{align*}
  &\limsup_{N_0 \to \infty} \Big(\frac{N_0}{M}\Big)^p \P\Big((1-\epsilon')\lVert \phi(x) - \phi(Z_k) \rVert \le \Phi_M(Z_k){\rm ~for~} k \in S,\\
  &2\lVert \phi(x) - \phi(Z_k) \rVert^{-1}\sup_{\lVert \phi(s) - \phi(t)\rVert \le \lVert \phi(x) - \phi(Z_k)\rVert}\lVert (\hat\phi-\phi)(s) - (\hat\phi-\phi)(t) \rVert > \epsilon' {\rm ~for~} k \notin S\Big) \\
  &\le  \Big[\Big(\frac{1-\epsilon}{1+\epsilon} - \frac{1+\epsilon}{1-\epsilon} d \epsilon' \Big)^{-1} \frac{f_1(\phi(x))}{f_0(\phi(x))} \Big]^{p}.
\end{align*}

Now we consider $\lvert S \rvert < p$. Recall that $W_k = \nu_0(B_{\phi(Z_k),(1-\epsilon')\lVert \phi(x) - \phi(Z_k) \rVert })$ and $V_k = \nu_0(B_{\phi(Z_k),\Phi_M(Z_k) })$ for any $k \in \zahl{p}$. Fix $S \subset \zahl{p}$. Recall that $W_{\max} = \max_{k \in S} W_k$ and $V_{\max} = \max_{k \in S} V_k$. We have
\begin{align*}
  & \P\Big((1-\epsilon')\lVert \phi(x) - \phi(Z_k) \rVert \le \Phi_M(Z_k){\rm ~for~} k \in S,\\
  &2\lVert \phi(x) - \phi(Z_k) \rVert^{-1}\sup_{\lVert \phi(s) - \phi(t)\rVert \le \lVert \phi(x) - \phi(Z_k)\rVert}\lVert (\hat\phi-\phi)(s) - (\hat\phi-\phi)(t) \rVert > \epsilon' {\rm ~for~} k \notin S\Big)\\
  &\le \P\Big(W_{\max} < V_{\max}, 2\min_{k \notin S}\lVert \phi(x) - \phi(Z_k) \rVert^{-1}\sup_{\lVert \phi(s) - \phi(t)\rVert \le \lVert \phi(x) - \phi(Z_k)\rVert}\lVert (\hat\phi-\phi)(s) - (\hat\phi-\phi)(t) \rVert > \epsilon'\Big)\\
  &\le \P\Big(W_{\max} < V_{\max} \le \eta_N \frac{M}{N_0}, 2\min_{k \notin S}\lVert \phi(x) - \phi(Z_k) \rVert^{-1}\sup_{\lVert \phi(s) - \phi(t)\rVert \le \lVert \phi(x) - \phi(Z_k)\rVert}\lVert (\hat\phi-\phi)(s) - (\hat\phi-\phi)(t) \rVert > \epsilon' \Big) \\
  &+ \P\Big(V_{\max} > \eta_N \frac{M}{N_0}\Big).
\end{align*}

Note that
$$2\min_{k \notin S}\lVert \phi(x) - \phi(Z_k) \rVert^{-1}\sup_{\lVert \phi(s) - \phi(t)\rVert \le \lVert \phi(x) - \phi(Z_k)\rVert}\lVert (\hat\phi-\phi)(s) - (\hat\phi-\phi)(t) \rVert > \epsilon'$$ implies that $$(\max_{k \notin S}\lVert \phi(x) - \phi(Z_k) \rVert)^{-1}\sup_{\lVert \phi(s) - \phi(t)\rVert \le \max_{k \notin S}\lVert \phi(x) - \phi(Z_k) \rVert}\lVert (\hat\phi-\phi)(s) - (\hat\phi-\phi)(t) \rVert >  \epsilon'/2.$$ Then
\begin{align*}
  & \P\Big(W_{\max} < V_{\max} \le \eta_N \frac{M}{N_0}, 2\min_{k \notin S}\lVert \phi(x) - \phi(Z_k) \rVert^{-1}\sup_{\lVert \phi(s) - \phi(t)\rVert \le \lVert \phi(x) - \phi(Z_k)\rVert}\lVert (\hat\phi-\phi)(s) - (\hat\phi-\phi)(t) \rVert > \epsilon' \Big)\\
  &\le \P\Big(\Big(\frac{1-\epsilon}{1+\epsilon} - \frac{1+\epsilon}{1-\epsilon} d \epsilon' \Big) \frac{f_0(\phi(x))}{f_1(\phi(x))}\max_{k \in S}\nu_1(B_{\phi(x),\lVert \phi(x) - \phi(Z_k) \rVert }) < V_{\max}, \\
  & (\max_{k \notin S}\lVert \phi(x) - \phi(Z_k) \rVert)^{-1}\sup_{\lVert \phi(s) - \phi(t)\rVert \le \max_{k \notin S}\lVert \phi(x) - \phi(Z_k) \rVert}\lVert (\hat\phi-\phi)(s) - (\hat\phi-\phi)(t) \rVert >  \epsilon'/2\Big).
\end{align*}

In the same way as the upper bound part in Part II, it suffices to show that
\begin{align*}
  \lim_{N_0 \to \infty} \Big(\frac{N_0}{M}\Big)^{p-\lvert S \rvert} \P\Big( (\max_{k \notin S}\lVert \phi(x) - \phi(Z_k) \rVert)^{-1}\sup_{\lVert \phi(s) - \phi(t)\rVert \le \max_{k \notin S}\lVert \phi(x) - \phi(Z_k) \rVert}\lVert (\hat\phi-\phi)(s) - (\hat\phi-\phi)(t) \rVert >  \epsilon'/2\Big) = 0.
\end{align*}

For any $u \in (0,1)$, conditional on $\max_{k \notin S} \nu_1(B_{\phi(x),\lVert \phi(x) - \phi(Z_k) \rVert }) = u$ and under $\max_{k \notin S} \lVert \phi(x) - \phi(Z_k) \rVert \le \delta$, we have 
\[
\max_{k \notin S} \lVert \phi(x) - \phi(Z_k) \rVert \ge [u/((1+\epsilon)f_1(\phi(x)) V_m)]^{1/m}. 
\]
Then by Assumption~\ref{asp:proof2},
\begin{align*}
  & \Big(\frac{N_0}{M}\Big)^{p-\lvert S \rvert} \P\Big( (\max_{k \notin S}\lVert \phi(x) - \phi(Z_k) \rVert)^{-1}\sup_{\lVert \phi(s) - \phi(t)\rVert \le \max_{k \notin S}\lVert \phi(x) - \phi(Z_k) \rVert}\lVert (\hat\phi-\phi)(s) - (\hat\phi-\phi)(t) \rVert >  \epsilon'/2\Big)\\
  &= \Big(\frac{N_0}{M}\Big)^{p-\lvert S \rvert} \int_0^1 (p-\lvert S \rvert) u^{p-\lvert S \rvert-1} \P\Big( (\max_{k \notin S}\lVert \phi(x) - \phi(Z_k) \rVert)^{-1} \sup_{\lVert \phi(s) - \phi(t)\rVert \le \max_{k \notin S} \lVert \phi(x) - \phi(Z_k) \rVert}\\
  &\lVert (\hat\phi-\phi)(s) - (\hat\phi-\phi)(t) \rVert  >  \epsilon'/2, \max_{k \notin S} \lVert \phi(x) - \phi(Z_k) \rVert \le \delta \Biggiven \max_{k \notin S} \nu_1(B_{\phi(x),\lVert \phi(x) - \phi(Z_k) \rVert }) = u \Big) \d u\\
  &+ \Big(\frac{N_0}{M}\Big)^{p-\lvert S \rvert} \P\Big(\lVert \hat\phi-\phi\rVert_\infty >  \epsilon' \delta/4\Big)\\
  &\le \Big(\frac{N_0}{M}\Big)^{p-\lvert S \rvert} \Big[\int_0^1 (p-\lvert S \rvert) u^{p-\lvert S \rvert-1} T_{\epsilon'}([u/(\lVert f_1 \rVert_\infty V_m)]^{1/m}) \d u + \P\Big(\lVert \hat\phi-\phi\rVert_\infty >  \epsilon' \delta/4\Big) \Big] = o(1).
\end{align*}

Then
\begin{align*}
  \limsup_{N_0 \to \infty} \Big(\frac{N_0}{M}\Big)^p \P\Big(Z_1,\ldots,Z_p \in A_\phi(x)\Big) \le \Big[\Big(\frac{1-\epsilon}{1+\epsilon} - \frac{1+\epsilon}{1-\epsilon} d \epsilon' \Big)^{-1} \frac{f_1(\phi(x))}{f_0(\phi(x))}\Big]^p.
\end{align*}

By $\epsilon,\epsilon'$ are arbitrary, we obtain
\begin{align*}
  \limsup_{N_0 \to \infty} \Big(\frac{N_0}{M}\Big)^p \P\Big(Z_1,\ldots,Z_p \in A_\phi(x)\Big) \le \Big(\frac{f_1(\phi(x))}{f_0(\phi(x))}\Big)^p.
\end{align*}

A matched lower bound is directly from the H\"older inequality.

\subsection{Proof of Theorem~\ref{thm:cons}}

Consider any $\epsilon \in (0,1)$ be given. From Assumption~\ref{asp:proof3}, $\cX$ is compact, and then $\cZ$ is also compact. Since $f_0, f_1$ are continuous over their compact supports, they are uniformly continuous, that is, there exists $\delta > 0$ such that for any $x,z \in \cZ$ with $\lVert \phi(z)-\phi(x) \rVert \le 3\delta$, we have $\lvert f_1(\phi(z)) - f_1(\phi(x)) \rvert \le \epsilon^2$, and for any $x,z \in \cX$ with $\lVert \phi(z)-\phi(x) \rVert \le 3\delta$, we have $\lvert f_0(\phi(z)) - f_0(\phi(x)) \rvert \le \epsilon^2$.

Let $\cE_1  =  \{x: f_1(\phi(x)) \le \epsilon\}$, $\cE_2  =  \{x: f_0(\phi(x)) \le \epsilon~{\rm or}~{\rm dist}(x,\partial \cZ) \vee {\rm dist}(x,\partial \cX) \le 3\delta\}$. We then seperate the proof into three cases. In the following, it suffices to consider $x$ such that $f_0(\phi(x))>0$ since we are considering $L_p$ risk.

{\bf Case I.} $x \notin \cE_1 \cup \cE_2$. In this case we have $f_0(\phi(x)), f_1(\phi(x)) > \epsilon$. Then for any $z \in \cX$ with $\lVert \phi(z) - \phi(x) \rVert \le 3\delta$, we have $z \in \cZ$ by the definition of $\cE_2$ and then $\lvert f_0(\phi(z)) - f_0(\phi(x)) \rvert \le \epsilon f_0(\phi(x))$ and $\lvert f_1(\phi(z)) - f_1(\phi(x)) \rvert \le \epsilon f_1(\phi(x))$.

Proceeding as in the proof of Case I in Lemma~\ref{lemma:moment,p}, we obtain
\begin{align*}
  \lim_{N_0\to\infty} \sup_{x \notin \cE_1 \cup \cE_2} \E \Big[ \Big\lvert \hat{r}_\phi(x) - r(\phi(x))\Big\rvert^p \Big] = 0,
\end{align*}
and then
\begin{align*}
  &\lim_{N_0\to\infty} \E \Big[ \Big\lvert \hat{r}_\phi(X) - r(\phi(X))\Big\rvert^p \ind\Big(X \notin \cE_1 \cup \cE_2\Big) \Big] = \lim_{N_0\to\infty} \E \Big[ \E\Big[\Big\lvert \hat{r}_\phi(X) - r(\phi(X))\Big\rvert^p\Big] \ind\Big(X \notin \cE_1 \cup \cE_2\Big) \Big]\\
  &\le \lim_{N_0\to\infty} \E \Big[ \sup_{x \notin \cE_1 \cup \cE_2} \E\Big[\Big\lvert \hat{r}_\phi(x) - r(\phi(x))\Big\rvert^p\Big] \ind\Big(X \notin \cE_1 \cup \cE_2\Big) \Big]
  = 0.
\end{align*}

{\bf Case II.} $x \in \cE_1 \setminus \cE_2$. In this case we have $f_0(\phi(x)) > \epsilon$. Then for any $z \in \cX$ with $\lVert \phi(z) - \phi(x) \rVert \le 3\delta$, we have $z \in \cZ$ by the definition of $\cE_2$ and then $\lvert f_0(\phi(z)) - f_0(\phi(x)) \rvert \le \epsilon f_0(\phi(x))$ and $f_1(\phi(z)) \le \epsilon + \epsilon^2$.

Proceeding as in the proof of Case II in Lemma~\ref{lemma:moment,p}, we obtain
\begin{align*}
  \lim_{N_0\to\infty} \sup_{x \in \cE_1 \setminus \cE_2} \E \Big[ \Big\lvert \hat{r}_\phi(x) - r(\phi(x))\Big\rvert^p \Big] = 0,
\end{align*}
and then
\begin{align*}
  \lim_{N_0\to\infty} \E \Big[ \Big\lvert \hat{r}_\phi(X) - r(\phi(X))\Big\rvert^p \ind\Big(X \notin \cE_1 \setminus \cE_2\Big) \Big] = 0.
\end{align*}

{\bf Case III.} $x \in \cE_2$. In this case we have $f_0(\phi(x)) \le \epsilon$ or ${\rm dist}(x,\partial \cZ) \vee {\rm dist}(x,\partial \cX) \le 3\delta$. Since $\cX$ is compact, the surface areas of $\cX$ and $\cZ$ are bounded, and $r$ is bounded uniformly, we have
\begin{align*}
  \limsup_{N_0\to\infty} \E \Big[ \Big\lvert \hat{r}_\phi(X) - r(\phi(X))\Big\rvert^p \ind\Big(X \in \cE_2\Big) \Big] \lesssim \P\Big(X \in \cE_2\Big) \lesssim \epsilon + \delta.
\end{align*}

Since $\epsilon$ is arbitrary and $\delta$ can be taken arbitrary small, we obtain
\begin{align*}
  \lim_{N_0\to\infty} \E \Big[ \Big\lvert \hat{r}_\phi(X) - r(\phi(X))\Big\rvert^p \ind\Big(X \in \cE_2\Big) \Big] = 0.
\end{align*}

Combining the above three cases completes the proof.

\subsection{Proof of Lemma~\ref{lemma:drequal}}

Note that for any $x \in \cX$,
\begin{align*}
  \frac{f_{\phi,X \given D=1}(\phi(x))}{f_{\phi,X \given D=0}(\phi(x))} = \frac{\P(D=1\given\phi(X)=\phi(x))\P(D=0)}{\P(D=0\given\phi(X)=\phi(x))\P(D=1)},
\end{align*}
and
\begin{align*}
  \frac{f_{X \given D=1}(x)}{f_{X \given D=0}(x)} = \frac{\P(D=1\given X=x)\P(D=0)}{\P(D=0\given X=x)\P(D=1)}.
\end{align*}
This completes the proof.

{
\bibliographystyle{apalike}
\bibliography{AMS}
}

\end{document}